\documentclass[aop,MSNbibl,seceqn,dvips]{arximspdf}
\usepackage{mathrsfs}
\usepackage{graphicx}

\doi{10.1214/12-AOP828} 
\volume{41}
\issue{6}
\pubyear{2013}
\firstpage{4162}
\lastpage{4213}

\makeatletter

\newcommand{\al}{\alpha} 
\newcommand{\be}{\beta} 
\newcommand{\gam}{\gamma} \newcommand{\Gam}{\Gamma}
\newcommand{\ep}{\varepsilon} 
\newcommand{\de}{\delta} \newcommand{\De}{\Delta}
 
\newcommand{\lm}{\lambda} \newcommand{\Lm}{\Lambda}
\newcommand{\ka}{\kappa} 
\newcommand{\vpi}{\varpi} \newcommand{\vph}{\varphi}
 
\newcommand{\si}{\sigma} 
 \newcommand{\Thet}{\Theta}

\newcommand{\E}{\mathbb{E}}

\renewcommand{\P}{\mathbb{P}}

\newcommand{\R}{\mathbb{R}}

\newcommand{\Z}{\mathbb{Z}}

\newcommand{\cG}{\mathcal{G}}
\newcommand{\cH}{\mathcal{H}}

\newcommand{\cN}{\mathcal{N}}

\newcommand{\cR}{\mathcal{R}}

\newcommand{\cT}{\mathcal{T}}

\newcommand{\wt}{\widetilde}
\newcommand{\wh}{\widehat}

\renewcommand{\vec}[1]{\underline{#1}}

\newcommand{\decto}{\downarrow}
\newcommand{\incto}{\uparrow}

\newcommand{\im}{\operatorname{im}} 

\newcommand{\loc}{\operatorname{loc}}

\newcommand{\supp}{\operatorname{supp}}

\newcommand{\abs}[1]{|#1|}

\newcommand{\angl}[1]{\langle#1\rangle}
\newcommand{\anglb}[1]{\langle#1\rangle}

\newcommand{\ball}[2]{B_{#1}(#2)}
\renewcommand{\d}[1]{\, d #1}
\newcommand{\f}{\frac}

\newcommand{\I}{\mathbf{1}}
\newcommand{\Ind}[1]{\I\{#1\}}
\newcommand{\pd}{\partial}
\newcommand{\set}[1]{\{#1\}}

\newcommand{\eqref}[1]{(\ref{#1})}
\newcommand{\Pois}{\operatorname{Pois}}
\newcommand{\relent}[2]{H(#1 \Vert  #2)}

\newcommand{\hypreg}{\textup{(H1)}}
\newcommand{\hypint}{\textup{(H2)}}
\newcommand{\hypintbeta}{\textup{(H2$^\be$)}}
\newcommand{\hypintB}{\textup{(H2$^B$)}}
\newcommand{\hypdiffbeta}{\textup{(H3$^\be$)}}
\newcommand{\hypdiffB}{\textup{(H3$^B$)}}

\newcommand{\rt}{o}
\newcommand{\rtsymbol}{{\bullet}}
\newcommand{\graphs}{\cG}
\newcommand{\graphsrt}{{\graphs_{\rtsymbol}}}
\newcommand{\graphspr}{{\graphs_{\rtsymbol\rtsymbol}}}
\newcommand{\graphsrtmark}{\cG_\rtsymbol^\spins}
\newcommand{\graphsprmark}{\cG_{\rtsymbol\rtsymbol}^\spins}
\newcommand{\trees}{\cT}
\newcommand{\treesrt}{{\trees_\rtsymbol}}
\newcommand{\edge}{\mathrm{e}}
\newcommand{\treesed}{{\trees_\rtsymbol^+}}
\newcommand{\treespr}{{\trees_{\edge}}}
\newcommand{\treereg}[1]{{\acrsf{T}_{#1}}}
\newcommand{\subt}[3][T]{{#1}_{{#2}}^{#3}}
\newcommand{\lwc}{\to_{\mathit{lwc}}}
\newcommand{\graphsrtz}{\cN_\rtsymbol}
\newcommand{\br}{\operatorname{br}}

\newcommand{\acrsf}[1]{{\textsf{#1}}}
\newcommand{\indir}{\uparrow}
\newcommand{\outdir}{\downarrow}

\newcommand{\Bethe}{{\mathrm{Bethe}}}
\newcommand{\BP}{\acrsf{BP}}
\newcommand{\danglb}[1]{[\![ #1 ]\!]}
\newcommand{\free}{\mathrm{f}}
\newcommand{\gibbs}{\mathscr{G}}
\newcommand{\perm}{\mathrm{p}} 
\newcommand{\simplex}{\De_\spins}
\newcommand{\simplexpr}{\De_{\spins^2}}
\newcommand{\spins}{\mathscr{X}}
\newcommand{\unif}{\vu}
\newcommand{\unimsr}{\mathscr{U}}
\newcommand{\vertex}{\mathrm{vx}}
\newcommand{\vpsi}{{\bar\psi}}
\newcommand{\vxi}{{\bar\xi}}
\newcommand{\hstar}{\cH^\star}
\newcommand{\hloc}{\cH_{\loc}}
\newcommand{\hlocf}{\cH_{\loc}^{\mathrm{fin}}}
\newcommand{\hlocopsi}{\cH_{\loc}^\circ[\psi]}
\newcommand{\hlocpmpsi}{\cH_{\loc}^\pm[\psi]}
\newcommand{\pr}[1]{\mathbf{#1}}
\newcommand{\prr}[1]{\bolds{#1}}
\newcommand{\bh}{\mathbf{h}}
\newcommand{\bde}{\bolds{\de}}
\newcommand{\vh}{\bar h}
\newcommand{\vm}{\bar m}
\newcommand{\vu}{\bar u}
\newcommand{\vde}{\bar{\de}}
\newcommand{\vka}{\bar{\ka}}

\newcommand{\ghost}{\star}
\newcommand{\vghost}{{v^\star}}
\newcommand{\es}{\vpi} 
\newcommand{\rc}{\pi} 
\newcommand{\conn}{\leftrightsquigarrow}
\newtheorem{thmm}{Theorem}[section]
\newtheorem{cor}[thmm]{Corollary}
\newtheorem{lem}[thmm]{Lemma}
\newtheorem{ppn}[thmm]{Proposition}

\newproclaim{dfn}[thmm]{Definition}
\newproclaim{ex}[thmm]{Example}
\newproclaim{rmk}[thmm]{Remark}

\makeatother

\begin{document}
\begin{frontmatter}

\title{Factor models on locally tree-like graphs}
\runtitle{Factor models on locally tree-like graphs}

\begin{aug}
\author[A]{\fnms{Amir} \snm{Dembo}\thanksref{t1}},
\author[B]{\fnms{Andrea} \snm{Montanari}\thanksref{t1,t2}}
\and
\author[C]{\fnms{Nike} \snm{Sun}\corref{}\ead[label=e3]{nikesun@stanford.edu}\thanksref{t1,t3}}
\runauthor{A. Dembo, A. Montanari and N. Sun}
\thankstext{t1}{Supported in part by NSF Grant DMS-11-06627.}
\thankstext{t2}{Supported in part by NSF Grant CCF-0743978.}
\thankstext{t3}{Supported in part by Dept. of Defense NDSEG Fellowship.}
\affiliation{Stanford University}
\address[A]{A. Dembo\\
Department of Mathematics\\
Stanford University\\
Building 380, Sloan Hall\\
Stanford, California 94305\\
USA\\
and\\
Department of Statistics\\
Stanford University\\
Sequoia Hall, 390 Serra Mall\\
Stanford, California 94305\\
USA} 
\address[B]{A. Montanari\\
Department of Electrical Engineering\\
Stanford University\\
350 Serra Mall\\
Stanford, California 94305\\
USA\\
and\\
Department of Statistics\\
Stanford University\\
Sequoia Hall, 390 Serra Mall\\
Stanford, California 94305\\
USA}
\address[C]{N. Sun\\
Department of Statistics\\
Stanford University\\
Sequoia Hall, 390 Serra Mall\\
Stanford, California 94305\\
USA}
\end{aug}

\received{\smonth{10} \syear{2011}}
\revised{\smonth{11} \syear{2012}}

%
\begin{abstract}
We consider homogeneous factor models on uniformly sparse graph
sequences converging locally to a (unimodular) random tree $T$, and
study the existence of the free energy density $\phi$, the limit of the
log-partition function divided by the number of vertices $n$ as $n$
tends to infinity. We provide a new interpolation scheme and use it to
prove existence of, and to explicitly compute, the quantity $\phi$
subject to uniqueness of a relevant Gibbs measure for the factor model
on $T$. By way of example we compute $\phi$ for the independent set (or
hard-core) model at low fugacity, for the ferromagnetic Ising model at
all parameter values, and for the ferromagnetic Potts model with both
weak enough and strong enough interactions. Even beyond uniqueness
regimes our interpolation provides useful explicit bounds on $\phi$.

In the regimes in which we establish existence of the limit, we show
that it coincides with the Bethe free energy functional evaluated at a
suitable fixed point of the belief propagation (Bethe) recursions on
$T$. In the special case that $T$ has a Galton--Watson law, this
formula coincides with the nonrigorous ``Bethe prediction'' obtained
by statistical physicists using the ``replica'' or ``cavity'' methods.
Thus our work is a rigorous generalization of these heuristic
calculations to the broader class of sparse graph sequences converging
locally to trees. We also provide a variational characterization for
the Bethe prediction in this general setting, which is of independent interest.
\end{abstract}

%

\end{frontmatter}
%
\section{Introduction}

Let $G=(V,E)$ be a finite undirected graph, and $\spins$ a finite
alphabet of \emph{spins}. A \emph{factor model} on $G$ is a probability
measure on the space of \emph{(spin) configurations} $\vec{\si}\in
\spins
^V$ of form
%
\begin{equation}
\label{efm} \nu^{\be,B}_{G,\vec{\psi}}(\vec{\si}) = \f{1}
{Z_{G,\vec{\psi
}}(\be,B)} \prod_{(ij)\in E}
\psi^\be(\si_i,\si_j) \prod
_{i\in V} \vpsi^B(\si_i),
\end{equation}
where $\psi\equiv\psi^\be$ is a symmetric function $\spins^2\to\R
_{\ge
0}$ parametrized by $\be\in\R$, $\vpsi\equiv\vpsi^B$ is a positive
function $\spins\to\R_{\ge0}$ parametrized by $B\in\R$ and
$Z_{G,\vec
\psi}(\be,B)$ is the normalizing constant, called the \emph{partition
function} (with its logarithm called the \emph{free energy}). The pair
$\vec{\psi}\equiv(\psi,\vpsi)$ is called a \emph{specification}
for the
factor model~\eqref{efm}.

In this paper we study the asymptotics of the free energy for sequences
of (random) graphs $G_n=(V_n=[n],E_n)$ in the \emph{thermodynamic
limit} $n\to\infty$. More precisely, with $Z_n(\be,B) \equiv
Z_{G_n,\vec
\psi}(\be,B)$ and $\E_n$ denoting expectation with respect to the law
of $G_n$, we seek to establish the existence of the \emph{free energy density}
%
\begin{equation}
\label{elim}\qquad \phi(\be,B) \equiv\lim_{n\to\infty}\phi_n(
\be,B), \qquad\mbox{where } \phi_n(\be,B)\equiv \f{1} {n}
\E_n \bigl[\log Z_n (\be,B) \bigr],
\end{equation}
and to determine its value. [In the literature, $\phi(\be,B)$ is also
referred to as the ``free entropy density'' or ``pressure.'']

The primary example we consider is the Potts model for a system of
interacting spins on a graph. Formally, the \emph{$q$-Potts model} on
$G$ with inverse temperature $\be$ and magnetic field $B$ is the
probability measure on $\spins^V=[q]^V$ (with $[q]\equiv\set
{1,\ldots
,q}$) given by
%
\begin{equation}
\label{epotts} \nu^{\be,B}_G(\vec{\si}) = \f{1}
{Z_G(\be,B)} \exp \biggl\{ \be\sum_{(ij)\in E}
\Ind{\si_i=\si_j} +B\sum_{i\in V}
\Ind{\si_i=1} \biggr\}.
\end{equation}
For $\be>0$ the system favors monochromatic edges and is said to be
\emph{ferromagnetic}, while for $\be<0$ the system favors edge
disagreements and is said to be \emph{anti-ferromagnetic}; the magnetic
field $B$ biases vertices toward the distinguished spin $1$. The
$q$-Potts model generalizes the Ising model which corresponds to the
case $q=2$. In analogy with the Potts model, in the general factor
model setting we continue to refer to $\be$ as the interaction or
temperature parameter and to $B$ as the magnetic field.

Potts models have been intensively studied in statistical mechanics
because of their key role in the theory of phase transitions \cite{MR641370}, critical phenomena \cite{MR1227790} and conformally invariant
scaling limits \cite{MR2280251}. As demonstrated, for instance, in
\cite{MR2875752} for the Ising model, determining the limit \eqref{elim}
plays a key role in characterizing the asymptotic structure of the
measures $\nu^{\be,B}_{G_n}$ in the thermodynamic limit. Potts models
are also of great interest in combinatorics: recall in fact that the
partition function admits a random-cluster representation (\cite{MR0359655,MR2243761}; see also Section~\ref{ssecpotts}), which at $B=0$ reads
\[
Z_G(\be,0)=\sum_{F\subseteq E}
\bigl(e^\be-1\bigr)^{\abs{F}} q^{k(F)},
\]
with $k(F)$ denoting the number of connected components induced by the
subset of edges $F\subseteq E$; cf. \eqref{erc}. Up to a
multiplicative constant this coincides with the Tutte polynomial
$T_G(x,y)$ of $G$ evaluated at $x = 1+q(e^{\beta}-1)^{-1}$,
$y=e^{\beta
}$; see, for example,  \cite{MR2187739}.

Mathematical statistical mechanics has focused so far on specific graph
sequences $G_n$, for example, on finite exhaustions of the rectangular
grid or other regular lattices in $d$ dimensions with $d$ fixed. Under
mild conditions on the sequence, existence of the free energy density
is a consequence of the following well-known argument (see, e.g., \cite{MR0289084}, Proposition~2.3.2): each graph $G_n$ can be decomposed into
smaller blocks by deleting a collection of edges whose number is
negligible in comparison with the volume. Consequently the sequence
$\log Z_{G_n}$ is approximately sub-additive in $n$, implying existence
of the limit; see \cite{MR0137800}.

In this paper we consider sparse graphs with a locally tree-like
structure---formally, graph sequences $G_n$ converging locally weakly
to (random) trees; see Definition~\ref{dlwc} below; see also \cite
{MR2354165,MR1873300}. Although the study of statistical
mechanics ``beyond $\Z^d$'' is not directly motivated by physics
considerations, physicists have been interested in models on
alternative graph structures for a long time (an early example being
\cite{MR0436850}). Moreover, the study of factor models on \emph{sparse}
graphs has many motivations coming from computer science and
statistical inference; see \cite{MR2643563,MR2518205}. Indeed, another
example we
will consider is the hard-core model for random independent sets on a
graph. In this model the configuration space is $\spins^V=\set{0,1}^V$,
where $0$ means unoccupied, and $1$ means occupied, and the only
configurations receiving positive measure are those for which no two
neighboring vertices are occupied, that is, so that the occupied
vertices form an independent set in the graph. Formally, the \emph
{independent set} or \emph{hard-core} model on $G$ with fugacity $\lm
>0$ is the probability measure on $\set{0,1}^V$ given by
%
\begin{equation}
\label{eis} \nu^\lm_G(\vec{\si}) = \f{1}
{Z_G(\lm)} \prod_{(ij)\in E} \Ind{
\si_i\si_j\ne1} \prod_{i\in V}
\lm^{\si_i},
\end{equation}
so that as $\lm$ increases the measure becomes more biased toward the
larger independent sets (and we write $B\equiv\log\lm$ for the magnetic
field). Due to the hard constraint preventing neighboring $1$s, this
system always has anti-ferromagnetic interactions and is of significant
interest in computer science. The independent set decision problem is
\textsc{np}-complete (via the clique decision problem \cite
{Cook1971CTP800157805047,MR0378476}). As $\lm$ increases the measure $\nu^\lm_G$ becomes increasingly
concentrated on the maximal independent sets; the optimization problem
of finding such sets is \textsc{np}-hard~\cite{MR584512} and hard to
approximate (\cite{Zuckerman2006LDE11325161132612} and references therein). The problem of
counting independent sets [i.e., computing $Z_G(1)$] for graphs of
maximum degree $\De$ is \#\textsc{p}-complete for $\De\ge3$ (\cite
{MR1791090} and references therein). Although there exists a \textsc
{ptas} (polynomial-time approximation scheme) for $Z_G(\lm)$ for $\lm$
below a certain ``uniqueness threshold'' \cite{MR2277139}, a series of
previous works (see \cite{MR2475668,101109FOCS201034,springerlink101007978364222935048} and references therein)
gave strong evidence that computation is hard for any $\lm$ above this
threshold. This question was resolved simultaneously in the subsequent
works \cite{arXiv12032226,arXiv12032602}, with \cite{arXiv12032602}
building on methods from
this paper.

Since infinite trees are nonamenable, $G_n$ cannot be decomposed by
removing a vanishing fraction of edges, so the preceding argument no
longer applies: in physics terms, \emph{surface effects are
nonnegligible even in the thermodynamic limit.} Despite this,
statistical physicists expect the free energy density \eqref{elim} to
exist on a large class of locally tree-like graphs. Even more
surprisingly, employing nonrigorous but mathematically sophisticated
heuristics such as the ``replica'' or ``cavity'' methods, they derive
exact formulas for this limit for a number of statistical mechanics
models on locally tree-like graphs; see, for example, \cite
{MR2518205} and
the references therein. The primary example considered in these works
is the graph chosen uniformly at random from those with $n$ vertices
and $m=m(n)$ edges, with $m/n\to\gamma\in\R$; such graphs converge
locally to the Galton--Watson tree with $\Pois(2\gamma)$ offspring
distribution. The Galton--Watson tree with general offspring
distribution can be obtained as the local weak limit of random graphs
with specified degree profile corresponding to the offspring
distribution; the physics heuristics extend to this and even more
general settings.

There is no good argument for why the limit \eqref{elim} exists; the
heuristic replica or cavity methods compute this limit starting from
the postulate that it exists. A significant breakthrough was achieved
by the interpolation method first developed by Guerra and Toninelli
\cite{MR1930572} for the Sherrington--Kirkpatrick model from
spin-glass theory, and then generalized to a number of statistical
physics models on sparse graphs \cite{MR1972121,MR2025238,MR2095932} and related constraint satisfaction problems
\cite
{MR2743259}. This method establishes super-additivity of
$\log Z_{G_n}$ which implies existence of the limit \eqref{elim}.
Unfortunately, this approach appears limited to models with \emph
{repulsive} interactions, that is, in which higher weight is given to
configurations in which neighboring vertices take different values. In
particular, it does not apply to the ferromagnetic Potts model. This is
especially puzzling because the heuristic physics predictions do not
distinguish between the two cases, and there is no fundamental reason
why the limit should be computable in one case and not in the other.
Further, this interpolation method only applies to very restricted
classes of graph sequences (typically, uniformly random given the
degree sequence); notably, existence of the limit is not proved for
\emph{deterministic} graph sequences. Finally, the method gives no way
to actually compute the limit, although interpolation has been used to
prove upper bounds \cite{MR1972121,MR2025238,MR2095932}.

In this paper we follow a different approach relying only on local weak
convergence of the graph sequence $(G_{n})_{n\ge1}$ to some limiting
(random) tree. The general idea is that the corresponding factor models
\eqref{efm} must converge (passing to a subsequence as needed), to a
Gibbs measure on the limiting tree; the task then ``reduces'' to the
one of identifying the correct limit. This is still a substantial
challenge because, in general, there is an uncountable number of
``candidate'' Gibbs measures for the limit. Nevertheless, this program
was carried through in \cite{MR2650042} for Ising models on graphs
converging locally to a Galton--Watson tree, under a ``uniform
sparsity'' assumption (Definition~\ref{dunifsparse}), on the degree
distribution. (It is further assumed in \cite{MR2650042} that the
distribution has finite second moment; this condition was relaxed in
\cite{MR2733399}, thereby handling the case of power law graphs.)
The result of \cite{MR2650042,MR2733399} provides also a fairly explicit
expression $\Phi(\be,B)$ for the free energy density, defined solely in
terms of the limiting tree. This expression coincides with the
so-called ``Bethe prediction'' of statistical physics, derived earlier
for random graphs with given degree distribution using the ``replica''
or ``cavity'' methods.

We develop this approach here in more generality. Rather than
considering a specific model such as the Ising, we establish results
for general abstract factor models satisfying mild regularity
conditions [see \hypreg\ below], covering in particular the Potts and
independent set models. We also make no distributional assumptions on
the graphs $G_n$ or the limiting random tree, other than some
integrability conditions [see Definition~\ref{dunifsparse} and
\hypint\ below]. In this setting we develop a general interpolation scheme
(Theorem~\ref{tfm}) which, under appropriate assumptions, bounds
differences $\phi_n(\be,B)-\phi_n(\be_0,B_0)$ in the limit $n\to
\infty$
by differences $\Phi(\be,B)-\Phi(\be_0,B_0)$ for $\Phi$ a functional
defined solely in terms of the limiting tree; see \eqref{ebethepred}.
We refer the reader to \cite{MR2023650} for a discussion of the
computation of limits of finite large random structures through
optimization procedures on the limiting infinite structure. Although we
continue to refer to this $\Phi(\be,B)$ as the ``Bethe prediction,'' we
remark that it is a considerable generalization of earlier formulas
obtained in the special case of Galton--Watson trees by statistical
physics methods. It is defined as the evaluation of the ``Bethe free
energy functional'' \eqref{ebethehstar} at a specific Gibbs measure
on the limiting tree, and corresponds to what physicists call the
``replica symmetric solution'': whereas it is expected to hold in the
high-temperature regime (i.e., with small enough interactions), for
many factor models it is incorrect at low temperature. However, we will
show that in ``uniqueness regimes,'' where the set of Gibbs measures on
the limiting tree corresponding to the factor model specification $\vec
\psi$ is a singleton, the upper and lower bounds of Theorem~\ref{tfm}
match to completely verify the Bethe prediction (Theorem~\ref{tfmunique}).

We then apply our interpolation scheme to compute the free energy
density in specific models. We verify the Bethe prediction for the
independent set model with low fugacity (Theorem~\ref{tis}) as a
consequence of Theorem~\ref{tfmunique}. Further, by using
monotonicity properties to restrict the set of relevant Gibbs measures,
we obtain results for the Potts model going beyond the implications of
Theorem~\ref{tfmunique}: for $q=2$ (Ising), we verify the Bethe
prediction for all $\be\ge0$, $B\in\R$ (Theorem~\ref{tising}),
extending the results of \cite{MR2650042,MR2733399} to general locally
tree-like graph sequences. For general $q$, we verify the prediction in
regimes of nonnegative $(\be,B)$ in which two specific Gibbs measures
on the limiting tree coincide, namely, the Gibbs measures arising from
free and 1 boundary conditions coincide, see Definition~\ref{dbd}
below. This condition is satisfied throughout the range $\set{\be\ge
0,B>0}$ for $q=2$; when $q\ge3$ there are regimes of nonuniqueness in
which it fails, but we will show that it is satisfied both at $\be$
sufficiently small and sufficiently large, that is, at high and low
temperatures.

Theorem~\ref{tfm} can give useful bounds even beyond uniqueness
regimes. As an illustration, we study the Potts model in the case that
$G_n$ converges locally to the $d$-regular tree $\treereg{d}$. In
Theorem~\ref{tpottsreg} we explicitly characterize the nonuniqueness
regime of this model and use Theorem~\ref{tfm} to give bounds for
$\phi
_n(\be,B)$ within this regime. In a subsequent work \cite
{arXiv12075500} we prove that in this setting, $\phi(\be,B)$ exists
and matches the lower bound of Theorem~\ref{tpottsreg}. We also
compute there the asymptotic free energy $\phi(\lm)$ (all $\lm\ge0$)
for the independent set model on $d$-regular bipartite graphs. In
contrast, for generic nonbipartite $G_n$ the consensus in physics is
for a full replica symmetry breaking for large enough $\lm$, and
consequently there does not exist even a heuristic prediction for the
free energy density in this regime.

As mentioned above, the Bethe prediction $\Phi(\be,B)$ is the
evaluation of the Bethe free energy functional at a specific Gibbs
measure on the limiting tree. This Gibbs measure has a characterization
in terms of ``messages'' $h_{x\to y} \equiv h_{(T,x\to y)}$ defined on
the directed edges $x\to y$ of each tree $T$, such that the entire
collection of messages is a fixed point of a certain ``belief
propagation'' or ``Bethe recursion''~\eqref{ebp}. Motivated by the
finite-graph optimization of \cite{MR2246363}, we provide a variational
characterization of the Bethe prediction (Theorem~\ref{tbetheopt})
which is of independent interest. In particular, this formulation
suggests nontrivial connections with large deviation principles.

\subsection{Local weak convergence and the Bethe prediction}
\label{ssecintrolwcbethe}

We study factor models on graphs which are ``locally tree-like'' in a
sense which we now formalize, starting with a few notation and
conventions. All graphs are taken to be \emph{undirected and locally
finite}. In a graph $G=(V,E)$, let $d$ denote graph distance, and for
$v\in V$ write $\ball{t}{v}$ for the sub-graph of $G$ induced by $\set
{w\in V \dvtx d(v,w)\le t}$. Write $v\sim w$ if $v,w$ are neighbors in $G$,
and write $\pd v$ for the set of neighbors of $v$ and $D_v\equiv|\pd
v|$. Let $\graphsrt$ denote the space of isomorphism classes of (finite
or infinite) rooted, connected graphs $(G,\rt)$. A metric on this space
is given by defining the distance between $(G_1,\rt_1)$ and $(G_2,\rt
_2)$ in $\graphsrt$ to be $1/(1+R)$ where $R$ is the maximal $r\in\Z
_{\ge0} \cup\set{\infty}$ such that $\ball{R}{\rt_1}\cong\ball
{R}{\rt
_2}$; with this definition $\graphsrt$ is a complete separable metric
space; see, for example,  \cite{MR2354165}. Let $\treesrt\subset
\graphsrt$ denote the closed subspace of \emph{\textup{(}rooted\textup{)} trees}
$T\equiv
(T,\rt)$, the acyclic elements of $\graphsrt$. We write $\subt{}{t}$
for $\ball{t}{\rt}$ in $T$, and in particular we use $\subt{}{0}$ to
denote the single-vertex tree. We now define the precise notion of
graph limits considered throughout this paper.

\begin{dfn}\label{dlwc}
Let $G_n = (V_n,E_n)$ ($n\ge1$) be a sequence of random graphs, and let
$I_n$ be a vertex chosen uniformly at random from $V_n$. We say $G_n$
\emph{converges locally \textup{(}weakly\textup{)}} to the random tree $T$ if for each
$t\ge0$, $\ball{t}{I_n}$ converges in law to $\subt{}{t}$ in the space
$\graphsrt$. We say in this case that the $G_n$ are \emph{locally tree-like}.
\end{dfn}

We will make repeated use of the fact that \emph{any} local weak limit
of graph sequences satisfies the ``unimodularity'' or
``mass-transport'' property whose definition we recall here; for a
detailed account, see \cite{MR2354165}. Let $\graphspr$ denote the
space of isomorphism classes of \emph{bi-rooted}, connected graphs with
a distinguished ordered pair, denoted $(G,i,j)$ (we do not require
$i\sim j$); $\graphspr$ is metrizable in a similar manner as
$\graphsrt$.

\begin{dfn}\label{dunim}
A Borel probability measure $\mu$ on $\graphsrt$ is said to be \emph
{unimodular} if it obeys the \emph{mass-transport principle},
%
\begin{eqnarray}
\label{eunim} \E_\mu \biggl[ \sum_{x\in V(G)}
f(G,\rt,x) \biggr] = \E_\mu \biggl[ \sum
_{x\in
V(G)} f(G,x,\rt) \biggr]
\nonumber
\\[-8pt]
\\[-8pt]
 \eqntext{\forall f \dvtx \graphspr\to
\R_{\ge0} \mbox{ Borel}.}
\end{eqnarray}
We say that $\mu$ is \emph{involution invariant} if \eqref{eunim}
holds when restricted to $f$
supported only on those $(G,x,y)$ with $x \sim y$.
\end{dfn}

A measure $\mu$ on $\graphs$ is involution invariant if and only if it
is unimodular (\cite{MR2354165}, Proposition~2.2). Unimodularity
corresponds to ``indistinguishability of the root;'' the concept first
appeared in \cite{MR1873300} where it was observed that local weak
limits of graph sequences must be unimodular (\cite{MR1873300}, Section~3.2). The converse of this implication remains a
well-known open question; see~\cite{MR2354165}.

\begin{dfn}\label{dunifsparse}
The graph sequence $G_n$ is \emph{uniformly sparse} if the $D_{I_n}$
are uniformly integrable, that is, if
\[
\lim_{L\to\infty} \Bigl( \limsup_{n\to\infty}
\E_n \bigl[D_{I_n} \Ind{D_{I_n} \ge L}\bigr] \Bigr)= 0
\]
(where $\E_n$ denotes expectation over the law of $G_n$ and $I_n$).
\end{dfn}

We assume throughout that $G_n$ ($n\ge1$) is a uniformly sparse graph
sequence converging locally weakly to the random tree $T$ of
(unimodular) law $\mu$ such that the root degree $D_\rt$ is nonzero
with positive $\mu$-probability; this entire setting is hereafter
denoted $G_n\lwc\mu$. In this setting we will describe general
conditions under which the asymptotic free energy $\phi(\be,B)$ for the
factor model \eqref{efm} exists and agrees with the ``Bethe energy
prediction,'' which we now describe. [If the sequence of random graphs
$G_n$ is such that $G_n\lwc\mu$ for almost every realization of the
sequence---as is the case for Erd\"{o}s--R\'{e}nyi random
graphs or random graphs with given degree distribution (see, e.g., \cite{MR2643563}, Propositions~2.5~and~2.6)---then our results apply
instead to the a.s. limit of $n^{-1}\log Z_n(\be,B)$.]

Let $\simplex$ denote the $(\abs{\spins}-1)$-dimensional simplex of
probability measures on the finite alphabet of spins $\spins$. Let
$\treesed$ denote $\treesrt$ without the single-vertex tree $\subt
{}{0}$, and let $\treespr\subset\graphspr$ be the space of isomorphism
classes of trees $T\in\treesed$ rooted at a directed edge $x\to y$,
written $(T,x\to y)$ or simply $x\to y$ for short. If $T$ has law $\mu$
for $\mu$ a unimodular measure on $\treesrt$, we let $\mu^\indir$ and
$\mu^\outdir$ denote the laws of $(T,J\to\rt)$ and $(T,\rt\to J)$,
respectively, for $J$ chosen uniformly at random from $\pd\rt$
conditioned on the event $\set{T\in\treesed}$. Involution invariance of
$\mu$ is then equivalent to
\begin{eqnarray*}
\E_{\mu^\outdir}\bigl[D_x f(T,x\to y)\bigr]& =&\f{
\E_\mu [\sum_{j\in\pd\rt} f(T,\rt\to j)
]} {\mu(D_\rt>0)} =\f{ \E_\mu [\sum
_{j\in\pd\rt}f(T,j\to\rt) ]} {\mu(D_\rt>0)}\\
& =&
\E_{\mu^\indir}\bigl[D_y f(T,x\to y)\bigr]
\end{eqnarray*}
(where $\rt$ corresponds to $x$ on the left-hand side and to $y$ on the
right-hand side), so in particular $\mu^\indir$ and $\mu^\outdir$ are
mutually absolutely continuous.

\begin{dfn}\label{dmsg}
The \emph{message space} is the space $\cH\equiv\cH_\mu$ of
measurable functions
\[
h \dvtx \treespr\times\R^2 \to\simplex, \qquad\bigl((T,x\to y),\be,B\bigr)
\mapsto\bigl(h^{\be,B}_{x\to y}(\si)\bigr)_{\si\in\spins},
\]
taken up to $\mu^\indir$-equivalence.
\end{dfn}

\begin{rmk}\label{rTij}
For $(T,x\to y)\in\treespr$ let $T_{x\to y}$ denote the component
sub-tree rooted at $x$ which results from deleting edge $(x,y)$ from
$T$. The interpretation of $h_{x\to y}$ is that it is a message from
$x$ to $y$ on the tree $T$, giving the distribution of $\si_x$ for the
factor model \eqref{efm} on $T_{x\to y}$. Indeed, although we do not
require it in general, in our concrete examples $h_{x\to y}$ depends
only on this component sub-tree.
\end{rmk}

\begin{figure}
\centering
\begin{tabular}{@{}cc@{}}

\includegraphics{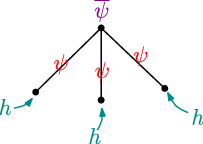}
 & \includegraphics{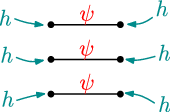}\\
\footnotesize{(a) Star graph $\subt{}{1}$} & \footnotesize{(b) Edge graph}
\end{tabular}
\caption{$\Phi^\vertex_T$ and $2\Phi^\edge_T$ are log-partition
functions of star and edge graphs.}\label{fstaredge}
\end{figure}

For $T\in\treesrt$ and $h\in\cH$, let
%
\begin{eqnarray}
\label{ePhiT} \Phi_T(\be,B,h) &\equiv&\Phi^\vertex_T(
\be,B,h)-\Phi^\edge_T(\be,B,h)
\nonumber
\\[-8pt]
\\[-8pt]
\nonumber
&\equiv&\Phi^\vertex_T(
\be,B,h)-\f{1} {2} \sum_{j\in\pd\rt} \Phi
^{(\rt
j)}_T(\be,B,h),
\end{eqnarray}
where ``$\vertex$'' and ``$\edge$'' indicate vertex and edge terms,
respectively:
%
\begin{equation}
\label{ePhivx} \Phi^\vertex_T(\be,B,h) \equiv \log \biggl\{
\sum_\si\vpsi(\si) \prod
_{j\in\pd\rt} \biggl( \sum_{\si_j} \psi(
\si,\si_j) h_{j\to\rt}(\si_j) \biggr) \biggr\},
\end{equation}
the log-partition function of the star graph $\subt{}{1}$ with boundary
conditions $h$ [see Figure~\ref{fstaredge}(a)] and
%
\begin{eqnarray}
\Phi^\edge_T(\be,B,h) &\equiv& \f12\sum
_{j\in\pd\rt}\Phi^{(\rt j)}_T(\be,B,h)
\nonumber
\\[-8pt]
\\[-8pt]
\nonumber
&=&\f12\sum
_{j\in\pd\rt} \log \biggl\{ \sum
_{\si,\si_j} \psi(\si,\si_j) h_{j\to\rt}(
\si_j) h_{\rt\to j}(\si) \biggr\},
\end{eqnarray}
half the log-partition function on $D_\rt$ disjoint edges with boundary
conditions $h$; see Figure~\ref{fstaredge}(b). (See Definition~\ref{dbd} below
for a detailed discussion of boundary conditions.)

We take the usual convention that the empty sum is zero, and the empty
product is one, so $\Phi_T=\log(\sum_\si\vpsi(\si))$ in case
$T=\subt{}{0}$.
Although we suppress it from the notation, in the above equations $\vec
\psi$ and $h$ are taken to be evaluated at $(\be,B)$.
The \emph{Bethe free energy functional} on $\cH$ for the factor model
\eqref{efm} on $G_n\lwc\mu$ is defined by
%
\begin{equation}
\label{ebethehstar} \Phi_\mu(\be,B,h) \equiv\E_\mu\bigl[
\Phi_T(\be,B,h)\bigr],
\end{equation}
provided the expectation exists; see Lemma~\ref{lbethehbd}.

\begin{dfn}
The \emph{belief propagation} or \emph{Bethe recursion} is the mapping
$\BP\equiv\BP^{\be,B}\dvtx \cH\to\cH$,
%
\begin{eqnarray}
\label{ebp}&& \bigl(\BP^{\be,B} h\bigr)_{x\to y}(\si)
\nonumber
\\[-8pt]
\\[-8pt]
\nonumber
&&\qquad \equiv\f{1}
{z_{x\to y}(\be,B)} \vpsi^B(\si) \prod
_{v \in\pd x \setminus y} \biggl( \sum_{\si_v}
\psi^\be(\si,\si_v) h_{v\to x}(
\si_v) \biggr),
\end{eqnarray}
with $z_{x\to y}(\be,B)$ normalizing constants.
For $\mu$ a measure on $\treesrt$ and fixed $(\be,B)$, let $\hstar
_\mu
(\be,B)$
denote the space of measurable\vadjust{\goodbreak} functions $h\dvtx \treespr\to\simplex$, again
taken up to $\mu^\indir$-equivalence, which are fixed points of the
Bethe recursion: that is, satisfying
%
\begin{equation}
\label{ehrecurs} h = \BP^{\be,B} h,\qquad \mu^\indir\mbox{-a.s.}
\end{equation}
The \emph{Bethe prediction} is that
the asymptotic free energy $\phi(\be,B)$ of \eqref{elim} exists and equals
%
\begin{equation}
\label{ebethepred} \Phi^\Bethe_\mu(\be,B) \equiv
\Phi_\mu\bigl(\be,B,h^\star\bigr)
\end{equation}
for $h^\star$ a certain element of $\cH^\star_\mu(\be,B)$. We often
drop the subscript $\mu$ when it is clear from context.
\end{dfn}

\begin{rmk}\label{rdefbethe}
In the case that the recursion \eqref{ehrecurs} has multiple
solutions ($\abs{\cH^\star_\mu(\be,B)}>1$), the \emph{Bethe~prediction}
is defined to be the supremum of $\Phi(\be,B,h^\star)$ over admissible
fixed points $h^\star$. While in the abstract factor model setting all
fixed points are admissible, in specific models typically there are
``natural'' criteria restricting the set of admissible fixed points. We
will demonstrate this in the Ising and Potts models where restrictions
are imposed by monotonicity and symmetry considerations.
\end{rmk}

The rationale behind the Bethe recursions and Bethe prediction is
explained in detail in \cite{MR2643563}, Section~3; see also \cite
{MR2518205}. In
brief, solutions to the Bethe recursions correspond to consistent
``boundary laws'' for the factor model on tree-like graphs; for further
details, see Remark~\ref{rbp} below. When $G$ is a finite tree, and
$\mu_G$ is the law of $(G,I)$ for $I$ a uniform element of $V$ (here
$\mu_G$ is a measure on $\treesrt$, but not necessarily unimodular),
the Bethe recursions have a unique solution, given by the so-called
``standard message set;'' see \cite{MR2643563}, Remark~3.5. In this setting
it holds exactly (see \cite{MR2643563}, Proposition~3.7) that
\[
\abs{V}^{-1} \log Z_G = \Phi_{\mu_G} =
\abs{V}^{-1} \sum_{v\in G}
\Phi_{(G,v)},
\]
where $\Phi_{(G,v)}$ is as defined by \eqref{ePhiT} with $T=(G,v)$.
The heuristic then is that for $G_n$ locally like the random tree
$T\sim
\mu$, the (normalized) free energy $\phi_n$ is approximated by $\Phi
_\mu
\equiv\E_\mu[\Phi_T]$ for $n$ large. We emphasize that \emph{no}
averaging over the vertices of the tree $T$ takes place in the
definition of $\Phi_T$; indeed for $T\in\treesrt$ the sub-trees
$\subt
{}{t}$ typically do \emph{not} converge locally weakly to $T$. For
example, when $T$ is the $d$-regular tree $\treereg{d}$, the subtrees
$\subt{}{t}$ converge locally weakly to the so-called $d$-canopy tree;
see, for example, \cite{MR2643563}, Lemma~2.8. Instead the averaging of
$\Phi
_{(G,v)}$ over the vertices $v\in G$ in the evaluation of $\Phi_{\mu
_G}$ corresponds to the averaging with respect to the law $\mu$ in the
evaluation of the Bethe prediction $\Phi_\mu$.

The following is a terminology which we adopt throughout the paper:

\begin{dfn}\label{dbd}
If $G$ is any graph and $U$ a sub-graph, the \emph{external boundary}
$\pd U$ of $U$\vadjust{\goodbreak} is the set of vertices of $G\setminus U$ adjacent to
$U$. Let $U^+$ denote the sub-graph of $G$ induced by the vertices in
$V_U\cup\pd U$.
For $U$ finite (so $U^+$ is finite, since $G$ is locally finite), and
$\nu^\ddagger$ a measure on $\spins^{\pd U}$, the factor model
on $U$ with \emph{$\nu^\ddagger$ boundary conditions} is the
probability measure on configurations $\vec{\si}_U\in\spins^{V_U}$
given by
%
\begin{equation}
\label{ebd-cond} \nu^\ddagger_{U,G,\vec{\psi}}(\vec{\si}_U)
\cong\int\prod_{(ij)\in E_{U^+}} \psi(\si_i,
\si_j) \prod_{i\in U} \vpsi(
\si_i) \d\nu^\ddagger(\vec{\si}_{\pd U}).
\end{equation}
(Throughout, $\cong$ indicates equivalence up to a positive normalizing
constant.)
The case in which $\nu^\ddagger$ gives probability one to the
identically-$\si_0$ spin configuration on $\pd U$ ($\si_0\in\spins
$) is
referred to as \emph{$\si_0$ boundary conditions} and denoted $\nu
^\ddagger=\nu^{\si_0}$, while the case in which $\nu^\ddagger$ is
uniform measure on $\spins^{\pd U}$ is referred to as \emph{free
boundary conditions} and denoted $\nu^\ddagger=\nu^\free$.
\end{dfn}
\subsection{Application to Ising, Potts and independent set}\label
{ssecintroapp}

Before formally stating our main theorem for general factor models, we
mention its consequences in some models of interest: we verify the
Bethe prediction for the ferromagnetic Ising model at all temperatures,
the ferromagnetic Potts model with field $B\ge0$ in uniqueness regimes,
and the independent set model with low fugacity $\lm$.

\subsubsection{Ising model}\label{sssecintroising}

The \emph{Ising model} is the Potts model \eqref{epotts} with $q=2$.
For convenience we use the equivalent formulation which takes $\spins
=\set{\pm1}$ and defines the probability measure on $\spins^V$
%
\begin{equation}
\label{eising} \nu^{\be,B}_G(\vec{\si}) = \f{1}
{Z_G(\be,B)} \exp \biggl\{ \be\sum_{(ij)\in E}
\si_i \si_j + B \sum_{i\in V}
\si_i \biggr\}.
\end{equation}
For $T\in\treesrt$ let $\vh^{t,+}_T\equiv\vh^{t,+,\be,B}_T$
denote the
root marginal for the Ising model of parameters $(\be,B)$ on $\subt
{}{t}$ with $+$ boundary conditions (i.e., with $\si_v$ conditioned to
be $+1$ for all $v$ at level $t+1$), and similarly define $\vh^\free_T$
corresponding to free boundary conditions. For $\ddagger\in\set
{\free
,+}$ let $\vh^\ddagger_T\equiv\vh^{\ddagger,\be,B}_T\equiv\lim_{t\to
\infty} \vh^{t,\ddagger,\be,B}_T$. (Existence of the limits $\vh
^\free
_T,\vh^+_T$ for the Ising model is an easy consequence of Griffiths's
inequality; see Lemma~\ref{lemgrif}.) We then define messages
$h^\ddagger\in\cH_\mu$ by
\[
h^\ddagger_{x\to y} = \vh^\ddagger_{ \subt{x\to y}{} }
\]
for $\subt{x\to y}{}$ as defined in Remark~\ref{rTij}. For $G_n\lwc
\mu$,
the Bethe free energy prediction for the Ising model with $\be\ge0$,
$B>0$ is $\phi(\be,B)=\Phi_\mu(\be,B,h^+)$.
This prediction was verified in \cite{MR2650042}, Theorem~2.4, for uniformly
sparse graph sequences converging locally weakly to Galton--Watson trees
subject to the second-moment condition $\E_\mu[D_\rt^2]<\infty$, which
was relaxed in \cite{MR2733399} to a $(1+\ep)$-moment condition. We
have the following generalization of this result to an arbitrary
limiting law.

\begin{thmm}\label{tising}
For the Ising model \eqref{eising} on $G_n\lwc\mu$,
\[
\phi(\be,B) = \Phi_\mu\bigl(\be,B,h^\free\bigr) =
\Phi_\mu\bigl(\be,B,h^+\bigr)
\]
for $\be\ge0$, $B>0$. Also $\phi(\be,B)=\phi(\be,-B)$ and $\phi
(\be,0)
= \lim_{B\to0} \phi(\be,B)$.
\end{thmm}

Note that in the Ising model we are able to characterize the free
energy density for \emph{all} $\be\ge0$. The underlying reason is that
for $B>0$, all boundary conditions dominating the free boundary
condition give rise to the same Gibbs measure on the limiting tree,
that is, $\vh^\free=\vh^+$. This phenomenon appears to be in line with
physicists' intuition that the Ising model always undergoes a
second-order phase transition. The physics argument suggests therefore
that the zero-magnetization phase becomes unstable below the critical
temperature. In other words, even with free boundary conditions, an
arbitrarily small external field $B>0$ is sufficient to drive the
system into the ``plus'' phase.

\subsubsection{Potts model}\label{sssecintropotts}

Throughout the remainder
let $(\be_0,B_0)\le(\be_1,B_1)$, where~$\le$ means coordinate-wise
less than or equal to. An \emph{interpolation path} is a piecewise
linear path, with each piece parallel to a coordinate axis, increasing
from $(\be_0,B_0)$ to $(\be_1,B_1)$ with respect to the partial order
$\le$.

We restrict our attention to the Potts model with $\be,B\ge0$. In this
regime we are able to use a random-cluster representation to extract
important monotonicity properties. For $T\in\treesrt$ and $\ddagger
\in
\set{\free,1}$ let $\vh^{t,\ddagger}_T\equiv\vh^{t,\ddagger,\be,B}_T$
denote the root marginal for the Potts model on $\subt{}{t}$ with
$\ddagger$ boundary conditions. Let $\vh^\ddagger_T\equiv\vh
^{\ddagger
,\be,B}_T \equiv\lim_{t\to\infty} \vh^{t,\ddagger,\be,B}_T$ (existence
of the limits $\vh^\free_T,\vh^1_T$ for the Potts model follows from
monotonicity properties of the random-cluster representation; see
Corollary~\ref{corpottsmonot}). We then define messages $h^\ddagger\in
\cH
_\mu$ by $h^\ddagger_{x\to y}=\vh^\ddagger_{\subt{x\to y}{}}$, and let
\[
\cR_\mu\equiv\bigl\{(\be,B) \dvtx h^\free=
h^1, \mu^\indir\mbox{-a.s.}\bigr\}.
\]
We also define
\[
\cR_\infty\equiv \bigl(\set{0} \times\R_{\ge0}\bigr) \cup\bigl(
\R_{\ge0}\times\set{\infty}\bigr) \cup \bigl(\set{\infty} \times
\R_{>0}\bigr).
\]

\begin{thmm}\label{tpotts}
For the Potts model \eqref{epotts} with $q>2$ and $\be,B\ge0$ on
$G_n\lwc\mu$, the following hold (with $\Phi\equiv\Phi_\mu$,
$\cR
\equiv\cR_\mu$):
\begin{longlist}[(a)]
\item[(a)]
If there exists an interpolation path contained in $\cR$ joining
$(\be,B)$ and
$\cR_\infty$,
then
\[
\phi(\be,B) = \Phi\bigl(\be,B,h^\free\bigr) = \Phi\bigl(
\be,B,h^1\bigr).
\]

\item[(b)]
If there exists an interpolation path from $(\be_0,B_0)$ to $(\be
_1,B_1)$ along which $h^\free$ is continuous (in the interpolation
parameter), then
\[
\liminf_{n\to\infty} \bigl[\phi_n(\be_1,B_1)-
\phi_n(\be_0,B_0)\bigr] \ge\Phi\bigl(
\be_1,B_1,h^\free\bigr)-\Phi\bigl(
\be_0,B_0,h^\free\bigr).
\]\eject\noindent
If $h^\free$ is replaced with $h^1$, then we have instead
\[
\limsup_{n\to\infty} \bigl[\phi_n(\be_1,B_1)-
\phi_n(\be_0,B_0)\bigr] \le\Phi\bigl(
\be_1,B_1,h^1\bigr)-\Phi\bigl(
\be_0,B_0,h^1\bigr).
\]
\end{longlist}
\end{thmm}

We obtain more explicit results when the limiting tree is the
$d$-regular tree $\treereg{d}$.

\begin{thmm}\label{tpottsreg}
For the Potts model \eqref{epotts} with $q>2$ and $\be,B\ge0$ on
$G_n\lwc\treereg{d}$,
the following hold (with $\Phi\equiv\Phi_{\treereg{d}}$, $\cR
\equiv
\cR_{\treereg{d}}$, and $\cR_{\ne}\equiv\set{\be,B\ge
0}\setminus
\cR$):

\begin{longlist}[(a)]
\item[(a)]
If $d=2$, $\cR_{\ne}=\varnothing$.
If $d>2$ and $q=2$, there exists $0<\be_-<\infty$ such that $\cR
_{\ne}=\set{B=0,\be>\be_-}$.
If $d>2$ and $q>2$, there exists $0<B_+<\infty$ and
smooth curves $\be_\free(B)\le\be_+(B)$ defined on $[0,B_+]$ with
$\be
_\free(B_+)=\be_+(B_+)$ such that
\[
\cR_{\ne} =\bigl\{B=0,\be\ge\be_\free(0)\bigr\} \cup\bigl\{0<B<
B_+, \be\in\bigl[\be_\free(B),\be_+(B)\bigr]\bigr\}.
\]

\item[(b)]
For $(\be,B)\notin\cR_{\ne}$, $\phi(\be,B)=\Phi(\be
,B,h^\free)=\Phi
(\be,B,h^1)$. If $(\be,B)\in\pd\cR_{\ne}$ with $\be=\be
_\free(B)$,
then $\phi(\be,B)=\Phi(\be,B,h^\free)$. If $(\be,B)\in\pd
\cR_{\ne}$
with $\be\ge\be_+(B)$, then $\phi(\be,B)=\Phi(\be,B,h^1)$. For
$(\be
,B)$ in the interior $\cR_{\ne}^\circ$ of $\cR_{\ne}$,
\begin{eqnarray*}
\max\bigl\{ \Phi\bigl(\be,B,h^1\bigr),\Phi\bigl(\be,B,h^\free
\bigr) \bigr\} &\le&\liminf_{n\to\infty} \phi_n(\be,B)
\\
&\le&\limsup_{n\to\infty}\phi_n(\be,B) \le\min\bigl\{ \wt
\Phi^\free(\be,B),\wt\Phi^1(\be,B) \bigr\},
\end{eqnarray*}
where
\begin{eqnarray*}
\wt\Phi^1(\be,B) &\equiv&\Phi\bigl(\be_\free(B),B,h^\free
\bigr) + \bigl[\Phi \bigl(\be,B,h^1\bigr) - \Phi\bigl(
\be_\free(B),B,h^1\bigr)\bigr],
\\
\wt\Phi^\free(\be,B) &\equiv&\Phi\bigl(\be_+(B),B,h^1\bigr)
- \bigl[\Phi\bigl(\be _+(B),B,h^\free\bigr) - \Phi\bigl(
\be,B,h^\free\bigr)\bigr].
\end{eqnarray*}
\end{longlist}
\end{thmm}

\begin{figure}[b]
\centering
\begin{tabular}{@{}c@{\hspace*{2pt}}c@{}}

\includegraphics{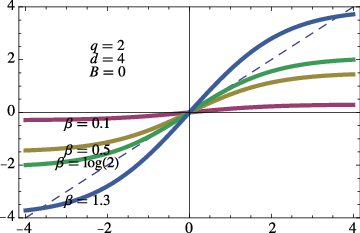}
 & \includegraphics{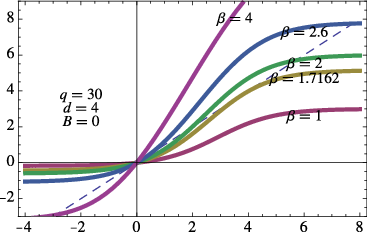}\\
\footnotesize{(a) Ising} & \footnotesize{(b) Potts}
\end{tabular}
\caption{Ising and Potts Bethe recursions.}\label{fbp}
\end{figure}

%

\begin{figure}
\centering
\begin{tabular}{@{}cc@{}}

\includegraphics{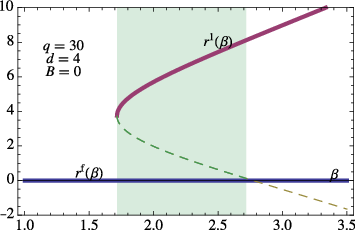}
 & \includegraphics{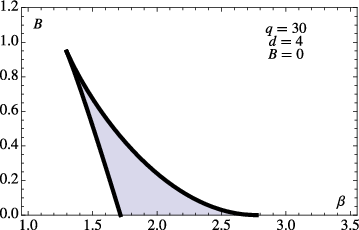}\\
\footnotesize{(a) BP fixed points} & \footnotesize{(b) Regime $\cR_{\ne}$ (shaded)}
\end{tabular}
\caption{Potts Bethe fixed points and the intermediate regime
$\cR
_{\ne}$.}\label{fpotts}
\end{figure}

\begin{figure}[b]

\includegraphics{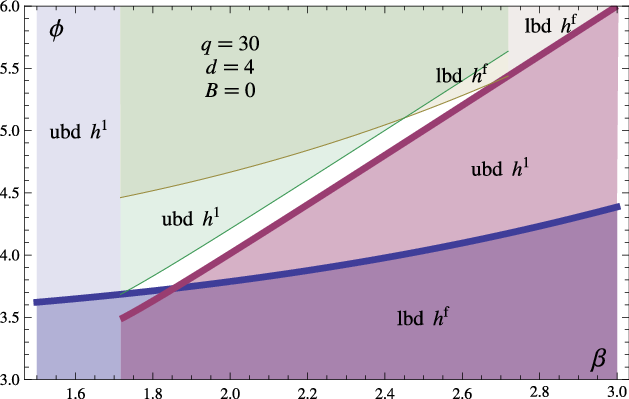}

\caption{Potts Bethe interpolation: the heavy (light) shaded regions are the
asymptotic lower (upper) bounds on $\phi_n$ given by Theorem~\protect\ref
{tpottsreg}; the bounds fail to match when $(\be,B)\in\cR_{\ne}$.
The Bethe prediction is the upper envelope of the thick lines.
In the figure, a shaded region marked ``lbd $h^\dagger$'' (resp., ``ubd
$h^\dagger$'') means an asymptotic bound on $\phi_n$ obtained from
interpolation using the asymptotic lower (resp., upper) bound on
$a^\edge_n(\be,B)$ by $a^\edge(\be,B,h^\dagger)$, in the notation of
Theorem~\protect\ref{tfm}. For example, the shaded region labeled ``ubd
$h^1$'' is an asymptotic (lower) bound on $\phi_n$ obtained by
interpolating from $\be=\infty$ using the asymptotic \emph{upper bound}
$\limsup a^\edge_n(\be,B)\le a^\edge(\be,B,h)^1$.}
\label{fpottsbethe}
\end{figure}

Figures~\ref{fbp}--\ref{fpottsbethe} highlight the difficulty in
analyzing the Potts model ($q>2$) as opposed to the Ising
model.\vadjust{\goodbreak}
Figure~\ref{fbp}(a) shows the Ising Bethe recursion parametrized
in terms of the log-likelihood ratio $r\equiv\log h(+)-\log h(-)$. For
sufficiently large $\be$ the recursion has three fixed points, but in
this case the $r=0$ fixed point is unstable, and we will see in the
proof of Theorem~\ref{tising} that adding a small magnetic field
resolves the nonuniqueness. The remaining plots were computed for the
Potts model with $q=30$ and $d=4$. Figure~\ref{fbp}(b) shows the
Potts Bethe recursion at $B=0$ restricted to those $h$ which are
symmetric among the spins $\ne1$, and parametrized by $r \equiv\log
h(1) - \log h(2)$. The fixed point at $r=0$ corresponds to $h^\free$
while the uppermost fixed point corresponds to $h^1$; Figure~\ref{fpotts}(a) shows how the fixed points vary with $\be$. In an
intermediate regime of $\be$-values [shaded in Figure~\ref{fpotts}(a)] both fixed points are stable, and perturbing by a
magnetic field does not resolve the nonuniqueness: indeed, Figure~\ref{fpotts}(b) shows that there is a two-dimensional region $\cR
_{\ne}$ of $(\be,B)$ values for which $h^\free\ne h^1$, making the
exact Bethe prediction inaccessible via our current interpolation
scheme. Figure~\ref{fpottsbethe} shows the discrepancy between the
upper and lower bounds of Theorem~\ref{tpottsreg}(b) inside $\cR_{\ne}$.

\subsubsection{Independent set model}\label{sssecintrois}

We consider the independent set model \eqref{eis} in the regime of low
fugacity. For $\ddagger\in\set{0,1}$ let $\vh^{t,\ddagger}_T\equiv
\vh
^{t,\ddagger,\lm}_T$ denote the root marginal on $\subt{}{t}$ with
$\ddagger$ boundary conditions on $\pd\subt{}{t}$: that is, $\vh
^{t,1}_T$ (resp., $\vh^{t,0}_T$) is calculated conditional on the
event of being fully occupied (unoccupied) at level $t+1$ of $T$. Let
$\vh^\ddagger_T\equiv\lim_{t\to\infty} \vh^{2t-1,\ddagger}_T$
(existence of the limits $\vh^0_T,\vh^1_T$ for the independent set
model follows from anti-monotonicity; see Section~\ref{ssecis}). We
then define messages $h^\ddagger\in\cH_\mu$ by $h^\ddagger_{x\to
y}=\vh
^\ddagger_{\subt{x\to y}{}}$, and let
\[
\lm_c\equiv\lm_{c,\mu} \equiv\inf\bigl\{\lm\ge0\dvtx
\mu^\indir\bigl(h^{0,\lm}_{x\to y}= h^{1,\lm}_{x\to y}
\bigr)<1\bigr\}
\]
denote the uniqueness threshold. For $T\in\treesrt$ we write
%
\begin{eqnarray}
\label{ebranch} \br T&\equiv&\inf \biggl\{y>0\dvtx \liminf_{\abs{\Pi}\to\infty}
\sum_{v\in\Pi} y^{-d(\rt,v)}=0 \biggr\}
\nonumber
\\[-8pt]
\\[-8pt]
\nonumber
& =&\sup \biggl
\{y>0\dvtx \liminf_{\abs{\Pi}\to\infty} \sum_{v\in\Pi}
y^{-d(\rt,v)}=\infty \biggr\}
\end{eqnarray}
(where the limit is taken over cutsets $\Pi$ of $T$ with distance
$\abs
{\Pi}$ from the root tending to infinity) for the \emph{branching
number} of $T$; see \cite{MR1062053}, Section~2.

\begin{thmm}\label{tis}
Consider the independent set model \eqref{eis} on $G_n\lwc\mu$, and
write $\lm_c\equiv\lm_{c,\mu}$.
\begin{longlist}[(a)]
\item[(a)]
If $\lm<\lm_c$ and the function $\lm\mapsto h^{0,\lm}=h^{1,\lm}$ has
total variation bounded by a deterministic constant on $[0,\log\lm]$, then
%
\begin{equation}
\label{eisbethe} \phi(\lm)=\Phi_\mu\bigl(\lm,h^0\bigr)=
\Phi_\mu\bigl(\lm,h^1\bigr),
\end{equation}
which converges to $\phi(\lm_c)$ as $\lm\incto\lm_c$.
\item[(b)]
If $\br T_{x\to y}\le\De-1$ $\mu^\indir$-a.s. for $\De$ a
deterministic constant, then \eqref{eisbethe} holds for $\lm<\lm_c$
with $\lm(\De-2)<1$.
\item[(c)]
If $\mu=\de_{\treereg{d}}$, then \eqref{eisbethe} holds for $\lm
\le\lm_c$.
\end{longlist}
\end{thmm}

For the $d$-regular tree $\treereg{d}$, the uniqueness threshold $\lm
_c(d)$ is $(d-1)^{d-1}/(d-2)^d$ (see \cite{MR844469}, Section~2), and
\cite{MR2277139}, Theorem~2.3, shows that $\treereg{d}$ has the lowest
value of $\lm_c$ among trees with maximum degree at most $d$. The
identity \eqref{eisbethe} has been proved in the case that the $G_n$
are random $d$-regular graphs \cite{MR2373815,MR2462251}. It
is also suggested by Weitz's \textsc{ptas} for $Z_G(\lm)$ on a finite
graph $G$ of maximum degree $\De$ and with $\lm<\lm_c(\De)$
(\cite{MR2277139}, Corollary~2.8). For $\mu$ a unimodular measure on
$\treesrt$
giving a local tree approximation to $G$ (in the sense of
Definition~\ref{dlwc}), $\lm_{c,\mu}$ is often an improvement over
$\lm
_c(\De)$, making it possible to compute $\phi(\lm)$ above $\lm
_c(\De)$
provided \hypdiffB\ can be verified. In \cite{arXiv12032602} the
interpolation scheme of Theorem~\ref{tis} is refined to give a
verification of the Bethe prediction on locally tree-like $d$-regular
\emph{bipartite} graphs for all $\lm>0$; this result is then leveraged
to show inapproximability of the hard-core partition function on
$d$-regular graphs above $\lm_c(d)$.

\subsection{Results for general factor models}\label{ssecintrofm}

We now state our results for the factor model \eqref{efm}.
With the convention $\log0\equiv-\infty$, let $\log\psi\equiv\xi$
and $\log\vpsi\equiv\vxi$, and impose the following regularity condition:
\begin{longlist}
\item[\hypreg] The specification is \emph{permissive}, that is,
$\vpsi
(\si)>0$ for all $\si\in\spins$,
and there exists a ``permitted state'' $\si^\perm\in\spins$ such that
$\min_\si\psi(\si,\si^\perm)>0$.

For any $\si\in\spins$, $\vxi^B(\si)$ is continuously
differentiable in
$B$. For any $\si,\si'\in\spins$, $\xi^\be(\si,\si')$ is either
identically $-\infty$ over all $\be$, or finite and continuously
differentiable in $\be$.
\end{longlist}
Recalling Definition~\ref{dmsg} of the message space $\cH\equiv\cH
_\mu
$, for $h\in\cH_\mu$ we can define $\vh\dvtx \treesrt\to\simplex$ up
to $\mu
$-equivalence by
%
\begin{equation}
\label{evh} \vh_T(\si)\cong \vpsi(\si) \prod
_{j\in\pd\rt} \biggl( \sum_{\si_j} \psi(
\si,\si_j) h_{j\to\rt}(\si_j) \biggr).
\end{equation}
In particular, if $h\in\hstar_\mu(\be,B)$ and $T\in\treesed$, then
comparing \eqref{evh} with \eqref{ebp} gives
%
\begin{equation}
\label{ehrecrt} \vh_T(\si)\cong \sum_{\si_j}
\psi(\si,\si_j) h_{\rt\to j} (\si) h_{j \to\rt} (
\si_j)
\end{equation}
independently of the choice of $j\in\pd\rt$. From now on, for $h\in
\cH
_\mu$, we will write $h\in\hstar$ to indicate that $h^{\be,B}\in
\hstar
_\mu(\be,B)$ for $(\be,B)$ in the range being considered.

\begin{figure}
\centering
\begin{tabular}{@{}cc@{}}

\includegraphics{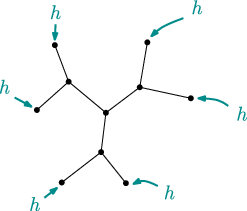}
 & \includegraphics{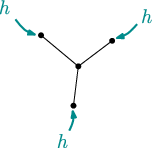}\\
\footnotesize{(a) $U=\subt{}{2}$} & \footnotesize{(b) $U=\subt{}{1}$}
\end{tabular}
\caption{A Bethe fixed point defines a consistent family of f.d.d. $\nu
^h_{U,T}$ (Remark~\protect\ref{rbp}).}\label{fconsistent}
\end{figure}

\begin{rmk}\label{rbp}
The elements of $\hstar$ are consistent with the recursion structure of
the tree in the following precise sense: for $T\in\treesrt$ and $U$ a
finite connected sub-graph of $T$, consider the factor model $\nu
^h_{U,T}$ on $U$ with boundary conditions $\si_v\sim h_{v\to p(v)}$
independently for $v\in\pd U$, where $p(v)$ denotes the (necessarily
unique) neighbor of $v$ inside $U$. Then the marginal of $\nu^h_{\subt
{}{t},T}$ on $\subt{}{t-1}$ is exactly the factor model
$\nu^{\scriptsize{\BP} h}_{\subt{}{t-1},T}$ on $\subt{}{t-1}$ with boundary
conditions $\si_u\sim(\BP h)_{u\to p(u)}$ independently for $u\in\pd
\subt{}{t-1}$, including any $u$ which are leaves of $\subt{}{t}$. This
statement remains valid if $\pd\subt{}{t}$ or even $\pd\subt
{}{t-1}$ is
empty, since if $\pd\subt{}{t}=\varnothing$ then $\nu^h_{\subt{}{t},T}$
is simply $\nu_T$ as defined by \eqref{efm}. Continuing the recursion
up the tree, we see that $h\in\hstar$ implies that the marginal law of
$\si_\rt$ will be $\vh_T$ as defined by \eqref{evh}. From this it is
easy to see that the measures $\nu^h_{U,T}$ form a consistent family of
finite-dimensional marginals (see Figure~\ref{fconsistent}), so by the
Kolmogorov consistency theorem they uniquely determine a probability
measure $\nu_T\equiv\nu^h_T$ belonging to $\gibbs_T$, the set of Gibbs
measures (or Markov random fields) associated to the specification
$\vec
\psi\equiv(\psi,\vpsi)$ on $T$.\setcounter{footnote}{3}\footnote{Strictly speaking the term
``Gibbs measures'' refers to the case $\psi>0$, but we will follow
common practice and say Gibbs measures also for the general case. For
the general theory of Gibbs measures see, for example, \cite{MR956646}.}
(In fact this mapping is one-to-one, e.g., by Remark~\ref
{rembed} below.)
Each $\nu_T$ belongs to a special class of measures in $\gibbs_T$ which
are called \emph{Markov chains} or \emph{splitting Gibbs measures} in
the literature, and the entire collection $(\nu_T)_{T\in\treesrt}$
arising from $h\in\hstar_\mu$ has a consistency property which leads us
to term them ``unimodular Markov chains'' or ``Bethe Gibbs measures;''
see Section~\ref{ssecdiscuss}.
\end{rmk}

In this general setting, the \emph{Bethe prediction} is the supremum of
$\Phi_\mu(\be,B,h)$ over $\hstar_\mu(\be,B)$; cf. Remark~\ref
{rdefbethe}.
(It will be shown in Lemma~\ref{lbethehbd} that $\Phi_\mu$ is
uniformly bounded on $\hstar_\mu(\be,B)$
subject to $\E_\mu[D_\rt^2]<\infty$; if further $\psi>0$, then
$\Phi_\mu
$ is in fact uniformly bounded on $\cH$ subject only to $\E_\mu
[D_\rt
]<\infty$.)
We define the following integrability condition for unimodular measures
$\mu$ on $\treesrt$ (not necessarily arising from a graph sequence):
\begin{longlist}[\hypint]
\item[\hypint] The probability measure $\mu$ on $\treesrt$
satisfies $\E
_\mu[D_\rt]<\infty$. If $\psi$ is not everywhere positive, then
furthermore $\E_\mu[e^{c D_\rt}]<\infty$ for all $c\in\R$.
\end{longlist}
Note that if $G_n\lwc\mu$ and $\psi>0$, then \hypint\ holds trivially
by the assumption of uniform sparsity. We will in fact justify our
interpolation scheme under a weaker assumption than \hypint; for the
exact condition see \hypintbeta, \hypintB\ in
Section~\ref{ssecinterp}.

\subsubsection{Bethe interpolation}
\label{sssecintrobethe}

We will deduce the results of Section~\ref{ssecintroapp} from the
abstract interpolation method given by Theorem~\ref{tfm} below, which
bounds differences of $\phi(\be,B)$ by differences of $\Phi(\be,B,h)$
($h\in\hstar$)
when the limiting expectation of a certain edge or vertex functional in
the finite graph (capturing resp. $\pd_\be\phi_n$ or $\pd
_B\phi
_n$) is bounded by the expectation of an analogous functional on the
infinite tree.

To be more precise, recall that $I_n$ denotes a uniformly random vertex
of $V_n$. Let $\anglb{\cdot}^{\be,B}_n$ denote expectation with respect
to $\nu_{G_n,\vec{\psi}}$, conditioned on $G_n$. For $h\in\hstar
_\mu(\be
,B)$ and $T\in\treesrt$, let $\danglb{}^{h,\be,B}_T$ denote expectation
with respect to $\nu^h_T$
(as defined in Remark~\ref{rbp}), conditioned on $T$, and define
\begin{eqnarray*}
a^\edge_n(\be,B)&\equiv&\f{1} {2} \E_n \biggl[
\sum_{j\in\pd I_n} \bigl\langle\pd_\be\xi(
\si_{I_n},\si_j)\bigr\rangle^{\be
,B}_n \biggr],\\
a^\edge(\be,B,h) &\equiv&\f{1} {2} \E_\mu \biggl[ \sum
_{j\in\pd\rt} \bigl[\!\bigl[\pd_\be\xi(
\si_\rt,\si_j)\bigr]\!\bigr]^{h,\be,B}_T \biggr],
\\
a^\vertex_n(\be,B)& \equiv&\E_n \bigl[ \bigl\langle
\pd_B\vxi(\si _{I_n})\bigr\rangle^{\be,B}_n \bigr],\\
a^\vertex(\be,B,h)& \equiv& \E_\mu \bigl[ \bigl[\!\bigl[
\pd_B\vxi(\si_\rt)\bigr]\!\bigr]^{h,\be,B}_T \bigr].
\end{eqnarray*}
The left-hand side expressions are the derivatives $\pd_\be\phi_n$, $\pd
_B\phi_n$ (Lemma~\ref{lemlogZderiv}). The right-hand side expressions are
the infinite-tree analogues, which, as we will show in Proposition~\ref
{pbethedbeta}, may be thought of as derivatives in $\be$ and $B$ of
$\Phi_\mu$.

\begin{ex}
For example in the Potts model \eqref{epotts} we have $\pd_B\vxi
(\si
)=\Ind{\si=1}$, so $a^\vertex_n(\be,B)$ is the expected density of $1$s
in the graph while $a^\edge_n(\be,B)$ is $1/n$ times the expected
number of edge agreements, both with respect to the Potts measure on
$G_n$. The infinite tree analogues of $a^\vertex_n$ and $a^\edge_n$
are
\begin{eqnarray*}
a^\vertex(\be,B,h) &=&\E_\mu \biggl[{\biggl(e^B\prod
_{j\in\pd\rt}\bigl[\bigl(e^\be-1\bigr)
h_{j\to\rt}(1)+1 \bigr]\biggr)}\\
&&\hspace*{22pt}{}\bigg/ \biggl(e^B\prod
_{j\in\pd\rt}\bigl[\bigl(e^\be-1\bigr)
h_{j\to\rt}(1)+1 \bigr]\\
&&\hspace*{41pt}{} + \sum_{\si\ne1} \prod
_{j\in\pd\rt}\bigl[\bigl(e^\be-1\bigr)
h_{j\to\rt}(\si)+1 \bigr]\biggr) \biggr],
\end{eqnarray*}
the $\nu^h_T$-probability (averaged over $T\sim\mu$) that the root spin
takes value $1$ and
\[
a^\edge(\be,B,h) =\f12\E_\mu \biggl[ \sum
_{j\in\pd\rt} \f{\sum_\si
e^\be h_{\rt\to j}(\si) h_{j\to\rt}(\si)} {1 + \sum
_\si(e^\be-1)
h_{\rt\to j}(\si) h_{j\to\rt}(\si)} \biggr],
\]
the $\nu^h_T$-expectation (averaged over $T\sim\mu$) of half the number
of edge agreements incident to the root.
\end{ex}

For interpolation in $\be$ on a compact interval $[\be_0,\be_1]$ using
some particular $h\in\hstar$, we require the following regularity
condition on $h$:
\begin{longlist}
\item[\hypdiffbeta] On $[\be_0,\be_1]$, for all $\si\in\spins$
it holds
$\mu^\indir$-a.s. that the function $\be\mapsto h^\be_{x\to y}(\si)$
is continuous with total variation in $\be$ bounded by a deterministic
constant depending only on $\be_0,\be_1$.
\end{longlist}
Likewise for interpolation in $B$ on a compact interval $[B_0,B_1]$
using $h\in\hstar$ we require
\begin{longlist}
\item[\hypdiffB] On $[B_0,B_1]$, for all $\si\in\spins$ it holds
$\mu
^\indir$-a.s. that the function $B\mapsto h^B_{x\to y}(\si)$ is
continuous with total variation in $B$ bounded by a deterministic
constant depending only on $B_0,B_1$.
\end{longlist}
The condition of boundedness in total
variation is implied for example whenever the functions $h$ are
(anti-)monotone in the interpolation parameter.

\begin{thmm}\label{tfm}
Let $\vec{\psi}\equiv(\psi,\vpsi)$ specify a factor model \eqref
{efm} on
$G_n\lwc\mu$ such that \hypreg\ and \hypint\ are satisfied.
\begin{longlist}[(a)]
\item[(a)]\label{tfmbeta}
If on $[\be_0,\be_1]$ we have $h\in\hstar$ satisfying \hypdiffbeta, and
%
\begin{equation}
\label{edbetaassump} \limsup_{n\to\infty} a^\edge_n(
\be,B) \le a^\edge(\be,B,h),
\end{equation}
then $\limsup_{n\to\infty} [\phi_n(\be_1,B)-\phi_n(\be_0,B)]
\le\Phi(\be_1,B,h)-\Phi(\be_0,B,h)$.

\item[(b)]\label{tfmB}
If on $[B_0,B_1]$, we have $h\in\hstar$ satisfying \hypdiffB, and
%
\begin{equation}
\label{edBassump} \limsup_{n\to\infty}a^\vertex_n(
\be,B) \le a^\vertex(\be,B,h),
\end{equation}
then
$\limsup_{n\to\infty} [\phi_n(\be,B_1)-\phi_n(\be,B_0)]
\le\Phi(\be,B_1,h)-\Phi(\be,B_0,h)$.
\end{longlist}
The same results hold if all inequalities are reversed, replacing limit
superior with inferior.
\end{thmm}

Conditions \eqref{edbetaassump}, \eqref{edBassump} (and their
reverses) are automatically verified in the following special case,
where we recall that $\gibbs_T$ denotes the set of Gibbs measures
associated to the specification $\vec{\psi}$ on $T$; cf. Remark~\ref{rbp}:

\begin{thmm}\label{tfmunique}
Let $\vec{\psi}\equiv(\psi,\vpsi)$ specify a factor model \eqref
{efm} on
$G_n\lwc\mu$ satisfying \hypreg\ and \hypint. We say that \emph
{uniqueness} holds if $\gibbs_T$ at $(\be,B)$ consists of a single
measure $\nu_T$, $\mu$-a.s. In this case, $\hstar_\mu(\be,B)$ is a
singleton.

\begin{longlist}[(a)]
\item[(a)]\label{tfmuniquebeta}
If on $[\be_0,\be_1]\times\set{B}$ uniqueness holds and the unique
element $h\in\hstar$ satisfies \hypdiffbeta, then
\[
\lim_{n\to\infty} \bigl[\phi_n(\be_1,B)-
\phi_n(\be_0,B)\bigr] =\Phi(\be _1,B,h)-\Phi(
\be_0,B,h).
\]
\item[(b)]\label{tfmuniqueB}
If on $\set{\be}\times[B_0,B_1]$ uniqueness holds and the unique
element $h\in\hstar$ satisfies \hypdiffB, then
\[
\lim_{n\to\infty} \bigl[\phi_n(\be,B_1)-
\phi_n(\be,B_0)\bigr] = \Phi(\be,B_1,h)-
\Phi(\be,B_0,h).
\]
\end{longlist}
\end{thmm}

Uniqueness for $\gibbs_T$ corresponds to the vanishing effect of
boundary conditions on $\pd\subt{}{t}$ as $t \to\infty$ (\cite{MR956646}, Chapter~7).
Dobrushin's uniqueness theorem (see, e.g., \cite{MR543198}) gives a sufficient condition for uniqueness to hold,
together with a bound on the rate of convergence of the root marginal
in $\subt{}{t}$ to the limit as $t\to\infty$. Note that if the
convergence rate is uniform in $\be,B$ then the continuity required in
\hypdiffbeta\ and \hypdiffB\ immediately follows. We will obtain
continuity in uniqueness regimes via a different route, making use of
certain monotonicity properties; see the proof of Theorem~\ref{tising}.

\subsubsection{Variational principle} \label{sssecintroopt}

We further develop the theory by providing a variational principle for
the Bethe prediction: we express $\Phi_\mu(\be,B)$ as an optimum of a
function $\Phi_\mu(\be,B,\bh)$ defined for $\bh$ in a larger space
$\hloc$ which, unlike $\hstar_\mu(\be,B)$, is \emph{independent}
of $\be
,B$. This alternative characterization of $\Phi_\mu$ is the
infinite-tree analogue of the finite-graph optimization problem that is
considered in~\cite{MR2246363}. Recall from Section~\ref{ssecintrolwcbethe} that $\treespr$ denotes the space of trees
rooted at a \emph{directed} edge.

\begin{dfn}\label{dhloc}
The \emph{local polytope} $\hloc\equiv\cH_{\loc,\mu}$ is the
space of
measurable functions
\[
\bh\dvtx \treespr\to\simplexpr, \qquad(T,x\to y) \mapsto\bh_{(T,x\to y)}\equiv
\bh_{xy},
\]
taken up to $\mu^\indir$-equivalence, such that:
\begin{longlist}[(ii)]
\item[(i)]$\bh_{xy}(\si,\si') = \bh_{yx}(\si',\si)$ for all $\si,\si
'\in
\spins$, and
\item[(ii)] for $T\in\treesed$, the one-point marginal $\vh_x(\si)\equiv
\vh
_{(T,x)}(\si) \equiv\sum_{\si_y} \bh_{xy}(\si,\si_y)$
is well-defined, that is, does not depend on the choice of $y\in\pd x$.
\end{longlist}
We also define
\begin{eqnarray*}
\hloc[\psi]&\equiv& \bigl\{\bh\in\hloc\dvtx \mu^\indir(\supp\bh\subseteq
\supp\psi)=1\bigr\},
\\
\hlocopsi&\equiv& \bigl\{\bh\in\hloc\dvtx \mu^\indir(\supp\bh=\supp
\psi)=1\bigr\}.
\end{eqnarray*}
In accordance with \eqref{evh}, we set
%
\begin{equation}
\label{evhdegzero} \vh_T(\si)\cong\vpsi(\si) \qquad\mbox{if } T=\subt{}
{0}.
\end{equation}
\end{dfn}

For \emph{fixed} $(\be,B)$, by symmetry of $\psi^\be$ and \eqref
{ehrecrt}, the space $\hstar_\mu(\be,B)$ has a natural mapping into
$\hloc$ given by
%
\begin{equation}
\label{eembed} h \mapsto\bh,\qquad \bh_{xy}\bigl(\si,\si'
\bigr) \cong\psi\bigl(\si,\si'\bigr) h_{x\to y}(\si)
h_{y\to x}\bigl(\si'\bigr).
\end{equation}
With $\vec{\psi}$ permissive this is in fact an embedding; see
Remark~\ref
{rembed}. We define the \emph{Bethe free energy functional on $\hloc
$} by
%
\begin{eqnarray}
\label{ebethehlocfirst} \Phi_\mu(\bh)& \equiv&\E_\mu
\biggl[\bigl\langle\vxi(\si_\rt) \bigr\rangle_{\vh_\rt} - (D_\rt-1)
H(\vh_\rt)
\nonumber
\\[-8pt]
\\[-8pt]
\nonumber
&&\hspace*{18pt}{}+ \f{1} {2} \sum_{j\in\pd\rt} \bigl
\{ \bigl\langle\xi(\si_\rt,\si_j)\bigr\rangle_{\bh_{\rt j}} + H(
\bh_{\rt j}) \bigr\} \biggr],
\end{eqnarray}
where $H(p)$ denotes the Shannon entropy $-\sum_k p_k \log p_k$ for $p$
a probability measure on a finite space. This is an infinite-tree
analogue of the definition of \cite{MR2246363}, (37)--(38), for finite
graphs. With the usual conventions $\log0\equiv-\infty$, $0\log0
\equiv0$ and $0 \log(0/0) \equiv0$, $\Phi_\mu$ is bounded above on
$\hloc$ whenever $\E_\mu[D_\rt]<\infty$, and we show in Lemma~\ref
{lembetheext} that for $\mu$ unimodular, this $\Phi_\mu$ extends the
previous definition~\eqref{ebethehstar} on $\hstar$ [under the
embedding \eqref{eembed}], provided the latter is finite. Furthermore,
writing $\relent{q}{p}$ for the relative entropy $\sum_k q_k \log
(q_k/p_k)$ between $q$ and $p$ (well defined for any nonnegative
reference measure $p$), for $\mu$ unimodular we can alternatively express
%
%
\begin{eqnarray}
\qquad\Phi_\mu(\bh) &= &-\E_\mu\bigl[\relent{\vh_\rt} {
\vpsi}\bigr] -\f{1} {2} \E_\mu \biggl[\sum_{j\in\pd\rt}
\relent{\bh_{\rt j}} {\psi} \biggr] - \E_\mu \bigl[
D_\rt H(\vh_\rt) \bigr]\label{ebethehloc}
\\
\nonumber
&=& -\E_\mu\bigl[\relent{\vh_\rt} {\vpsi}\bigr]-\f{1} {2}
\E_\mu \biggl[\sum_{j\in\pd\rt} \bigl\{ \relent{
\bh_{\rt j}} {\psi} + H(\vh_\rt) + H(\vh_j) \bigr
\} \biggr]
\\
\label{ebethehlocrelent} &= &- \E_\mu\bigl[ \relent{\vh_\rt} {
\vpsi} \bigr] - \f{1} {2} \E_\mu \biggl[ \sum_{j \in\pd\rt}
\relent{\bh_{\rt j}} {\vh_\rt\times_\psi
\vh_j} \biggr],
\end{eqnarray}
where $(\vh_\rt\times_\psi\vh_j)(\si_\rt,\si_j)\equiv\vh_\rt
(\si_\rt
) \psi(\si_\rt,\si_j) \vh_j(\si_j)$, and unimodularity is used in the
second identity.

This extended definition of $\Phi_\mu$ provides the following
variational principle for the Bethe free energy:

\begin{thmm}\label{tbetheopt}
Let $\vec{\psi}\equiv(\psi,\vpsi)$ specify a factor model \eqref{efm}
satisfying \hypreg, and let $\mu$ be a unimodular measure on
$\treesrt$
with $\E_\mu[D_\rt]<\infty$.
\begin{longlist}[(a)]
\item[(a)]
$\wt\Phi_\mu(\be,B)\equiv\sup_{\bh\in\hloc}\Phi_\mu(\be
,B,\bh)$ is
continuous in $(\be,B)$.
\item[(b)]
Any local maximizer of $\Phi_\mu(\be,B)$ belongs to $\hlocopsi$. Any
stationary point of $\Phi_\mu(\be,B)$ belonging to $\hlocopsi$ is the
image under \eqref{eembed} of an element of $\hstar_\mu(\be,B)$. In
particular, if $\Phi_\mu$ attains its supremum on $\hloc$, then
\[
\wt\Phi_\mu(\be,B) =\max_{h\in\hstar_\mu(\be,B)}
\Phi_\mu(\be,B,h)\equiv\Phi ^\Bethe_\mu(\be,B),
\]
so that the Bethe free energy is also continuous in $(\be,B)$.
\end{longlist}
\end{thmm}

Although we do not pursue this point, we mention that even in specific
models where the abstract definition of $\Phi^\Bethe$
is supplanted by $\Phi(\be,B,h)$ for some ``naturally'' distinguished
$h$, an adaptation of Theorem~\ref{tbetheopt}
[involving a restricted subspace of $\hloc$ which is independent of
$(\be,B)$], may be relevant.

\begin{rmk}
In case $G_n\lwc\treereg{d}$ the $d$-regular tree,
$\hloc$ is parame\-trized by a single measure $\bh_{xy}$ on $\spins^2$
whose one-point marginals are required to agree,
and the formula \eqref{ebethehlocrelent} simplifies to
%
\begin{equation}
\label{ebethereg} -\Phi_\mu(\bh) = \relent{\vh_0} {
\vpsi} + \f{d} {2} \relent{\bh _{0 1} } {\vh_0
\times_\psi\vh_1}.
\end{equation}
For $\vec{\si}\in\spins^{V_n}$ let $L^\vertex_n\equiv n^{-1}\sum_{i\in
V_n} \de_{\si_i}$ and $L^\edge_n \equiv(2\abs{E_n})^{-1}\sum_{(ij)\in
E_n} [\de_{(\si_i,\si_j)}+\de_{(\si_j,\si_i)}]$ denote the induced
empirical and pair empirical measures, respectively. If $G_n$ is
$d$-regular, then the one-point marginals of $L^\edge_n$ coincide with
$L^\vertex_n$, and
\[
\phi_n = \log\abs{\spins} + \f{1} {n} \E_n \biggl[ \log
\E_{\unif_n} \exp \biggl\{ n\angl{\vxi}_{L^\vertex_n} + \f{nd} {2} \angl{
\xi}_{L^\edge_n} \biggr\} \biggr],
\]
where the law of $\vec{\si}$ is the uniform measure $\unif_n$ on
$\spins
^{[n]}$ and $\E_{\unif_n}$ denotes expectation with respect to $\unif
_n$ (with $G_n$ fixed).

If $(G_n)$ is an independent sequence of uniformly random $d$-regular
graphs and $\vec{\si}_n\sim\unif_n$, one might guess that for a.e. $(G_n)$ the induced sequence $L^\edge_n$ satisfies a large deviation
principle (\textsc{ldp}) with good rate function
%
\begin{equation}
\label{ebetheline} I(\bh_{01}) = \relent{\vh_0} {\unif} +
\f{d} {2} \relent{\bh_{0 1} } {\vh _0 \times
\vh_1},
\end{equation}
where $\unif\equiv\unif_1$.
If this were the case, it would be an immediate consequence of
Varadhan's lemma (see \cite{MR1619036}, Section~4.3.1) that $\phi_n
\to
\wt
\Phi_\mu(\be,B)$ (as defined in Theorem~\ref{tbetheopt}) for
\emph
{any} factor model satisfying \hypreg. However, for many of these
models the Bethe prediction is known to fail at low temperature for
$d\ge3$. So, while Theorem~\ref{tbetheopt} suggests a potential
connection to large deviations theory, such a connection would be
highly nontrivial and applicable only in certain regimes of $(\be,B)$.

One special case in which everything trivializes is the (rooted)
infinite line $\treereg{2}$, the local weak limit of the simple
path\vadjust{\goodbreak}
$G_n$ on $n$ vertices. In this case $\unif_n$ may be viewed as the law
of a stationary reversible Markov chain on $\spins$ with transitions
$q(\si,\si')=\unif(\si')$ and reversing measure $\unif$, and it is
well-known (see, e.g., \cite{MR1619036}, Theorem~3.1.13) that the associated
pair empirical measure $L^\edge_n$ satisfies an \textsc{ldp} with good
rate function $I(\bh_{01}) = \relent{\bh_{01}(\si,\si')}{\vh
_0(\si)
q(\si,\break\si')}$ which matches~\eqref{ebetheline}. The implication of
Varadhan's lemma is also easy to see: a factor model on the
simple path $G_n$ with general positive specification $\vec{\psi}$
corresponds in the limit
$n\to\infty$ to a reversible Markov chain with transition kernel $p$
and positive
reversing measure $\pi$ given by
\[
p\bigl(\si,\si'\bigr) = \f{1} {\rho} \wt\psi\bigl(\si,
\si'\bigr) \f{m(\si')} {m(\si )},\qquad \pi(
\si)=m(\si)^2,
\]
where $\rho$ and $m$ are the Perron--Frobenius eigenvalue and
eigenvector of the symmetric positive $\abs{\spins}$-dimensional
matrix with entries $\wt\psi(\si,\si')\equiv \psi(\si,\si')\times\vpsi
(\si
)^{1/2} \vpsi(\si')^{1/2}$.
The Bethe free energy functional \eqref{ebethereg} is then maximized
at $\bh_{01}(\si,\si') = \wt{\psi}(\si,\si') m(\si) m(\si
')/\rho$,
where it takes the value $\Phi_\mu(\bh)=\log\rho$ which coincides with
$\phi$ by the Perron--Frobenius theorem; see, for example,
\cite{MR1619036}, Theorem~3.1.1.
\end{rmk}

\subsection*{Outline of the paper}

\begin{itemize}
\item In Section~\ref{secfm} we prove the abstract interpolation
results. Section~\ref{ssecprelim} presents some preliminary lemmas
which will be useful in our proofs. Our main result for abstract factor
models, Theorem~\ref{tfm}, is proved in Section~\ref{ssecinterp}.
Section~\ref{ssecdiscuss} contains the specialization of this theorem
to the uniqueness case (Theorem~\ref{tfmunique}) and also contains
discussion on unimodular Markov chains (or Bethe Gibbs measures).
Section~\ref{ssecis} shows how to deduce our result for independent
set (Theorem~\ref{tis}) from Theorem~\ref{tfm}.

\item In Section~\ref{sechloc} we prove the variational
characterization Theorem~\ref{tbetheopt} for the Bethe free energy
prediction, establishing in particular the correspondence between
interior stationary points $\bh\in\hlocopsi$ of $\Phi_\mu$ and fixed
points $h\in\cH^\star$ of the Bethe recursion. We further provide in
Proposition~\ref{plocmax} a simple criterion for such stationary
points to be local maximizers.

\item Section~\ref{secisingpotts} contains applications of our
abstract results to the Ising and Potts models. In Section~\ref{ssecising} we prove Theorem~\ref{tising}, generalizing the results
of \cite{MR2650042,MR2733399}. In Section~\ref{ssecpotts} we prove
Theorem~\ref{tpotts} by appealing to a random-cluster representation.
Finally, Section~\ref{ssecpottsreg} analyzes the $d$-regular case and
proves Theorem~\ref{tpottsreg}.\
\end{itemize}

\section{Bethe interpolation for general factor models}\label{secfm}

\subsection{Preliminaries}\label{ssecprelim}

We begin with some straightforward observations on the boundedness of
the free energy $\phi_n$ and the Bethe free
energy $\Phi_\mu$ as defined on $\cH$, and we prove that the mapping
\eqref{eembed} of $\hstar$ into $\hloc$ is in fact an embedding for
permissive specifications.

\begin{lem}\label{lemlogZderiv}
For the factor model \eqref{efm} satisfying \hypreg\ on $G_n\lwc\mu$,
the functions $\phi_n(\be,B)$ are uniformly bounded and equicontinuous
on compact regions of $(\be,B)$, with
\begin{eqnarray}
\label{elogZderivexch} \pd_\be\phi_n(\be,B) &=& \f{1} {n}
\E_n\bigl[\pd_\be\log Z_n(\be,B)\bigr],
\nonumber
\\[-8pt]
\\[-8pt]
\nonumber
\pd_B\phi_n(\be,B) &=& \f{1} {n} \E_n\bigl[
\pd_B\log Z_n(\be,B)\bigr].
\end{eqnarray}
Further,
\begin{eqnarray*}
\f{1} {n} \pd_\be\log Z_n(\be,B)& =&\f{1} {2} \sum
_{j\in\pd I_n} \bigl\langle \pd_\be\xi(
\si_{I_n},\si_j)\bigr\rangle ^{\be,B}_n, \\
\f{1} {n}
\pd_B\log Z_n(\be,B) &=& \bigl\langle\pd_B \vxi(
\si_{I_n})\bigr\rangle^{\be,B}_n,
\end{eqnarray*}
with the convention $\pd_\be\xi(\si,\si')\equiv0$ in case $\xi
(\si,\si
')\equiv-\infty$.
\end{lem}

\begin{pf}
The expressions for $n^{-1} \pd_\be\log Z_n(\be,B)$ and $n^{-1}\pd
_B\log Z_n(\be,B)$ are obtained by a straightforward computation. Now
note that if $G_n\lwc\mu$, then the uniform sparsity assumption gives
%
\begin{equation}
\label{eedgefreq} \f{1} {n} \E_n\bigl[\abs{E_n}\bigr] =\f{1} {2}
\E_n[ D_{I_n}] \to\f{1} {2} \E_\mu[D_\rt]<
\infty.
\end{equation}
Let $(\be,B)$ vary within a given compact region. By \hypreg\ we have
$\xi,\vxi\le\xi_{\max}$ as well as
$\xi(\si^\perm,\cdot),\vxi\ge\xi_{\min}$. Therefore,
\[
\bigl(1 + \abs{E_n}/n\bigr) \xi_{\min} \le n^{-1} \log
Z_n(\be,B) \le\log\abs{\spins} + \bigl(1 + \abs{E_n}/n\bigr)
\xi_{\max}
\]
so $\phi_n = n^{-1} \E_n[\log Z_n(\be,B)]$ is uniformly bounded by
uniform sparsity. The exchange of differentiation and integration in
\eqref{elogZderivexch} is justified by
Vitali's convergence theorem, in view of the boundedness of $\pd_\be
\xi
$, $\pd_B\vxi$ and the uniform integrability of $|E_n|/n$.
It follows furthermore that $\pd_\be\phi_n(\be,B)$ and $\pd_B \phi
_n(\be,B)$ are bounded uniformly in $n$, from which equicontinuity follows.
\end{pf}

\begin{lem}\label{lbethehbd}
Let $\vec{\psi}\equiv(\psi,\vpsi)$ specify a
factor model \eqref{efm} satisfying \hypreg, and let $\mu$ be a
unimodular measure on $\treesrt$. For any compact region of $(\be,B)$
there exists a deterministic constant $C<\infty$ such that:
\begin{longlist}[(a)]
\item[(a)]
$\abs{\Phi_T(\be,B,h)} \le C (D_\rt^2+1)$ for any $h\in\hstar$, and

\item[(b)]
if further $\psi>0$, then $\abs{\Phi_T(\be,B,h)} \le C(D_\rt+1)$ for
any $h\in\cH$.
\end{longlist}
Hence, on any compact region of $(\be,B)$, $\Phi_\mu$
is uniformly bounded on $\hstar_\mu$ provided $\E_\mu[D_\rt
^2]<\infty$,
and if $\psi>0$, uniformly bounded on $\cH_\mu$
subject only to $\E_\mu[D_\rt]<\infty$.
\end{lem}

\begin{pf}
Let $\xi_{\min},\xi_{\max}$ be as in the proof of Lemma~\ref
{lemlogZderiv}. Then, for any $h\in\cH$,
\[
\log\abs{\spins} + (D_\rt+1) \xi_{\max} \ge
\Phi^\vertex_T(h) \ge(D_\rt +1)
\xi_{\min}.
\]
If $\psi>0$, then we also have
\[
D_\rt\xi_{\max} \ge2 \Phi^\edge_T(h)
\ge D_\rt\xi_{\min},
\]
so $\abs{\Phi_T(\be,B,h)} \le C(D_\rt+1)$ on $\cH$, which proves
(b). For general permissive $\psi$, the preceding lower
bound on $\Phi^\edge_T(h)$ may fail, but
\eqref{ehrecurs} implies that for $h\in\hstar$,
%
\begin{equation}
\label{ehpermlbd} \log h_{\rt\to j}\bigl(\si^\perm\bigr) \ge
D_\rt(\xi_{\min} - \xi_{\max}) - \log\abs{\spins}\qquad
\forall j\in\pd\rt.
\end{equation}
Therefore,
\begin{eqnarray*}
D_\rt\xi_{\max} \ge2 \Phi^\edge_T(h)
&\ge& \sum_{j\in\pd\rt} \bigl( \xi_{\min} + \log
h_{\rt\to j}\bigl(\si ^\perm\bigr) \bigr) \\
&\ge& D_\rt\bigl(
\xi_{\min} - \log\abs{\spins}\bigr) + D_\rt^2 (
\xi_{\min} - \xi _{\max}),
\end{eqnarray*}
which proves (a).
\end{pf}

\begin{rmk}\label{rembed}
It is now easy to see that the mapping \eqref{eembed} of $\hstar_\mu
(\be,B)$ into $\hloc$ is injective: if $h,h'\in\cH$ give rise to the
same $\bh$, then
\[
h_{x\to y}(\si) h_{y\to x}\bigl(\si^\perm\bigr) =
z_{x,y} h'_{x\to y}(\si) h'_{y\to x}
\bigl(\si^\perm\bigr)\qquad \forall\si \in \spins
\]
for $z_{x,y}$ a positive scaling factor. If $h,h'\in\hstar_\mu(\be,B)$,
then \eqref{ehpermlbd} implies that $\mu^\indir$-a.s. both
$h_{y\to
x}$ and $h'_{y\to x}$ give positive measure to $\si^\perm$.
Therefore, $\mu^\indir$-a.s. the $\abs{\spins}$-dimensional vectors
$h_{x\to y}$ and $h'_{x\to y}$ are equivalent up to scaling, and since
both are probability measures on $\spins$, we must have $h=h'$ $\mu
^\indir$-a.s. as claimed.
\end{rmk}

\subsection{Bethe interpolation}\label{ssecinterp}

We now prove Theorem~\ref{tfm}(a). The result is for
fixed $B$, so we suppress it from the notation.
The proof of Theorem~\ref{tfm}(b) is very similar and will
be given in brief at the end of this section.

Our interpolation procedure relies on the proposition below which
expresses $\Phi_\mu$ as the integral of its partial derivative with
respect to $\be$ only,
ignoring the dependence on $\be$ through the function $h$. Recall that
although it is suppressed from the notation, $\vec{\psi}$ and $h$ depend
on $\be$,
and are taken to be evaluated at $\be$ in expressions such as $\Phi
_T(\be,B)$. We will prove our result under the following integrability
condition, which by \eqref{ehpermlbd} is a relaxation of~\hypint:
\begin{longlist}
\item[\hypintbeta]
The probability measure $\mu$ on $\treesrt$ satisfies $\E_\mu[D_\rt
]<\infty$. If $\psi$ is not everywhere positive, then furthermore,
\[
\E_\mu \biggl[ \sum_{j\in\pd\rt} \sup
_{\be\in[\be_0,\be_1]} \f {1} {h^\be _{j\to\rt}(
\si^\perm)} \biggr]<\infty.
\]
\end{longlist}
We define the analogous condition \hypintB\ on an interval $[B_0,B_1]$.

\begin{ppn}\label{pbethedbeta}
Let $\vec{\psi}\equiv(\psi,\vpsi)$ be a specification satisfying
\hypreg,
and $\mu$ a unimodular measure on $\treesrt$. If on $[\be_0,\be_1]$ we
have $h\in\hstar$ satisfying \hypintbeta\ and~\hypdiffbeta, then
\[
\Phi(\be_1,h)-\Phi(\be_0,h) = \int_{\be_0}^{\be_1}
a^\edge(\be,h) \d\be.
\]
\end{ppn}

\begin{pf}
For fixed $T\in\treesrt$ we shall regard $\Phi_T$ simply as a function
of a vector $(\be,h_{x\to y})_{x\to y\in\subt{}{1}}$ in $(1+2\abs
{\spins
}D_\rt)$-dimensional euclidean space (with $h$ depending on $\be$).
We begin by computing the partial derivatives of this function with
respect to $\be$ and $h$.
We abbreviate $\wh h^\be_{\rt\to j}(\si) \equiv(\BP^\be h)_{\rt
\to
j}(\si)$ for the
belief propagation mapping of \eqref{ebp}, which for fixed $T$ and
each $j\in\pd\rt$ is a well-defined function on the same euclidean
space as $\Phi_T$. Making use of \hypreg\ we find
%
%
\begin{eqnarray}
\label{edPhivdbe} \f{\pd\Phi_T^\vertex} {\pd\be}(\be,h) &=&
\sum_{j\in\pd\rt} \f{\sum_{\si,\si_j}
\pd_\be\xi(\si,\si_j) \psi(\si,\si_j)
h_{j\to\rt}(\si_j) \wh h^\be_{\rt\to j}(\si)}
{\sum_{\si,\si_j} \psi(\si,\si_j)
h_{j\to\rt}(\si_j) \wh h^\be_{\rt\to j}(\si)}
,
\\
\f{\pd\Phi_T^{(\rt j)}} {\pd\be}(\be,h) &= &\f{ \sum
_{\si,\si_j} \pd_\be\xi(\si,\si_j) \psi(\si,
\si_j) h_{j\to\rt}(\si_j) h_{\rt\to j}(\si)} {
\sum_{\si,\si_j} \psi(\si,\si_j)
h_{j\to\rt}(\si_j) h_{\rt\to j}(\si) }.\label{edPhiedbe}
\end{eqnarray}
If $h\in\hstar$, then $\wh h^\be=h$, therefore
(recalling the notation $\danglb{}^{h,\be}_T$ from Section~\ref{sssecintrobethe}) we re-express the above as
\[
\f{\pd\Phi_T^\vertex} {\pd\be}(\be,h) = \sum
_{j\in\pd\rt} \bigl[\!\bigl[\pd_\be\xi(\si_\rt,\si
_j)\bigr]\!\bigr]^{h,\be}_T, \qquad \f{\pd\Phi_T^{(\rt j)}}
{\pd\be}(\be,h) = \bigl[\!\bigl[\pd_\be\xi(\si_\rt,
\si_j)\bigr]\!\bigr]^{h,\be}_T,
\]
and combining gives
%
\begin{equation}
\label{ebethedbeta} \f{\pd\Phi_T} {\pd\be} (\be,h) = \f{1} {2} \sum
_{j\in\pd\rt} \bigl[\!\bigl[\pd_\be\xi(
\si_\rt,\si _j)\bigr]\!\bigr]^{h,\be}_T \equiv
a^\edge_T(\be,h).
\end{equation}
Likewise we compute that for $T\in\treesed$,
\begin{eqnarray*}
\f{\pd\Phi^\vertex_T(\be,h)} {\pd h_{\rt\to j}(
\si)}&=&0,\\
 \f{\pd\Phi^\vertex_T(\be,h)} {\pd
h_{j\to\rt}(\si_j)}& =& \wh g^\be_{\si_j}(j
\to\rt;h) \equiv\f{\sum_\si\psi(\si,
\si_j) \wh h^\be_{\rt\to j}(\si)} { \sum
_{\si',\si_j'} \psi(\si',\si_j'
) h_{j\to\rt}(\si_j') \wh
h^\be_{\rt\to j}(\si')},
\\
\f{\pd\Phi^\edge_T(\be,h)} {\pd h_{j\to\rt}(
\si_j)} &=& \f{1} {2} g^\be_{\si_j}(j\to\rt;h),\qquad \f{
\pd\Phi^\edge_T(\be,h)} {\pd h_{\rt\to j}(\si)} =
\f{1} {2} g^\be_\si(\rt\to j;h),
\end{eqnarray*}
where $g^\be_\si$ is the same as $\wh g^\be_\si$ but with $h$ in place
of $\wh h$. Note that for permissive $\psi$ and any $\si\in\spins$,
%
\begin{equation}
\label{egubd} \qquad\wh g^\be_\si(x\to y;h) \le\f{ \sum
_{\si'} \psi(\si',\si) \wh
h^\be_{y\to x}(\si') } { \sum
_{\si'} \psi(\si',\si^\perm)
h_{x\to y}(\si^\perm) \wh h^\be
_{y\to x}(\si') } \le\f{\psi_{\max}^\be}
{\psi_{\min}^\be h_{x\to y}(\si^\perm
)}.
\end{equation}
If further $\psi>0$ everywhere, then $\wh g^\be_\si(x\to y;h) \le
\psi
^\be_{\max}/\psi^\be_{\min}$ is uniformly bounded on $[\be_0,\be_1]$.

Consider now a small sub-interval $[\be,\be+\de]$ of $[\be_0,\be_1]$.
Writing $\De_{\be,\de} h \equiv h^{\be+\de}-h^\be$ and applying
the mean
value theorem to the differentiable function $t \mapsto\Phi_T(\be+ t
\de, h+ t \De_{\be,\de} h)$ for $t\in[0,1]$ gives
%
\begin{eqnarray}
\label{ephiincrement} &&\Phi_T(\be+\de,h)-\Phi_T(\be,h)
\nonumber
\\[-8pt]
\\[-8pt]
\nonumber
&&\qquad= \f{\pd\Phi_T} {\pd\be}\bigl(\be+ t\de, h^\be+ t
\De_{\be,\de} h\bigr) \de+ \Gam_T(\be,\de) +
E_T(\be,\de)
\end{eqnarray}
for some $t\equiv t_{\be,\de}\in[0,1]$, where
\begin{eqnarray*}
\Gam_T(\be,\de) &\equiv&\sum_\si\sum
_{j\in\pd\rt} \biggl\{\f{\pd\Phi_T} {\pd
h_{j\to\rt}(\si)}(\be,h) \De_{\be,\de} h_{j\to\rt}(\si)\\
&&\hspace*{40pt}{} +\f{\pd
\Phi_T} {\pd h_{\rt\to j}(\si)}(\be,h) \De_{\be,\de}
h_{\rt\to j}(\si) \biggr\},
\\
E_T(\be,\de) &\equiv&\sum_{\si} {\sum
_{x\to y}}^* \biggl\{\f{\pd\Phi_T} {\pd
h_{x\to y}(\si)}\bigl(\be+t \de, h^\be+ t \De_{\be
,\de}
h\bigr)\\
&&\hspace*{102pt}{} -\f{\pd\Phi_T} {\pd h_{x\to y}(\si)}(\be,h) \biggr\}
\De_{\be,\de} h_{x\to y}(\si),
\end{eqnarray*}
and $\sum^*_{x\to y}$ indicates the sum over the $2 D_\rt$ directed
edges $x\to y$ within $\subt{}{1}$.

Setting $\de\equiv\de_m \equiv(\be_1-\be_0)/m$, we now sum $\Phi
(\be
+\de_m,h)-\Phi(\be,h)$ over
$\be\in\Pi_m \equiv\set{\be_0 + k \de_m \dvtx 0\le k < m}$ and analyze
separately the contribution of
each term on the right-hand side of \eqref{ephiincrement}:

\begin{longlist}[(a)]
\item[(a)]
First we show that
$\E_\mu[\Gam_T(\be,\de)]=0$ for any $[\be,\be+\de]\subseteq
[\be_0,\be
_1]$. Indeed, since $h\in\hstar$
we have $\wh h^\be=h^\be$ and $\wh g^\be=g^\be$. Therefore,
\[
\Gam_T(\be,\de) = \f{1} {2} \sum_\si
\sum_{j\in\pd\rt} \bigl\{g^\be_\si(j
\to\rt;h) \De_{\be,\de} h_{j\to\rt}(\si) - g^\be_\si(
\rt\to j;h) \De_{\be,\de} h_{\rt\to j}(\si) \bigr\}.
\]
The result then follows from unimodularity of $\mu$, subject to $\mu
$-integrability of
\[
\sum_\si\sum_{j\in\pd\rt}
\bigl|g^\be_\si(j\to\rt;h) \De _{\be,\de}
h_{j\to\rt}(\si)\bigr|.
\]
Clearly $|\De_{\be,\de} h_{x\to y} (\si)| \le2$ so integrability
certainly holds when $\psi>0$,
since $\E_\mu[D_\rt]<\infty$ and $g^\be_\si$ is deterministically
uniformly bounded on $[\be_0,\be_1]$ as noted above. More generally,
for permissive $\psi$ the required $\mu$-integrability follows from~\eqref{egubd} and \hypintbeta.

\item[(b)]
The total contribution of the first term on the right-hand side of
\eqref{ephiincrement} is
\[
A_m \equiv\de\E_\mu \biggl[ \sum
_{\be\in\Pi_m} \f{\pd\Phi_T} {\pd\be}\bigl(
\be+t_{\be,\de} \de,h^\be+t_{\be,\de} \De_{\be
,\de}h
\bigr) \biggr].
\]
Observe that $A_m=\int Y_m \d(\lm\times\mu)$ where $\lm$ is Lebesgue
measure on $[\be_0,\be_1]$ and
\[
Y_m\bigl(\be',T\bigr) \equiv\sum
_{\be\in\Pi_m} \mathbf{1}\bigl\{\be\le\be'<\be+\de\bigr\} \f{\pd
\Phi_T} {\pd\be}\bigl(\be+t_{\be,\de} \de,h^\be+t_{\be,\de}
\De_{\be
,\de}h\bigr).
\]
For $(\lm\times\mu)$-a.e. $(\be',T)$, this sum has at most one
nonzero term, in which the argument of $\pd_\be\Phi_T$ converges by
\hypdiffbeta\
to $(\be',h^{\be'})$ as $m\to\infty$. From \hypreg, \eqref{ebp} and
the computation of $\pd_\be\Phi_T$
in \eqref{edPhivdbe}--\eqref{edPhiedbe}, we see that
$\pd_\be\Phi_T(\be,h)$ is continuous in $(\be,h)$. Therefore,
$Y_m(\be
',T)\to a^\edge_T(\be',h)$, $(\lm\times\mu)$-a.e.
Furthermore, \hypreg\ implies that $\abs{\pd_\be\xi}\le C$
uniformly on
$[\be_0,\be_1]$ for
some deterministic constant $C$, so
$\abs{Y_m}\le2 C D_\rt$ for all $m$, $(\lm\times\mu)$-a.e.  see
\eqref
{edPhivdbe} and \eqref{edPhiedbe}. Dominated convergence then gives
\[
\lim_{m\to\infty} A_m = \int a^\edge_T
\bigl(\be',h\bigr) \d(\lm\times\mu) = \int_{\be_0}^{\be_1}
a^\edge\bigl(\be',h\bigr) \d\be'.
\]

\item[(c)]
The contribution of the final term in \eqref{ephiincrement} is $\E
_\mu
[E_{T,m}]$ where
\[
E_{T,m}\equiv\sum_{\be\in\Pi_m}
E_T(\be,\de),
\]
and we conclude the proof by showing that $\lim_{m\to\infty}\E_\mu
[ E_{T,m}]=0$.

Indeed, it is not hard to see that $\lim_{m\to\infty}E_{T,m}=0$ $\mu
$-a.s.: by the uniform bound on total variation assumed in \hypdiffbeta,
there exists deterministic $C$ such that
\begin{eqnarray*}
&&\abs{E_{T,m}} \le C {\sum_{x\to y}}^* \max
_\si \sup_{\be\in[\be_0,\be_1]} \sup_{t\in[0,1]}
\biggl|\f{\pd\Phi_T} {\pd h_{x\to y}(\si)}\bigl(\be+t
\de_m, h^\be+ t \De _{\be,\de_m} h\bigr) \\
&&\hspace*{232pt}{}-\f{\pd
\Phi_T} {\pd h_{x\to y}(\si)}(\be,h)\biggr|
\end{eqnarray*}
$\mu$-a.s., uniformly in $m$.
It also follows from \hypdiffbeta\ that $\mu$-a.e. $h^\be$ is
uniformly continuous on $[\be_0,\be_1]$. Using \hypreg, the partials
$\pd_h \Phi_T$ computed above are uniformly continuous in $(\be,h)$ for
$\be\in[\be_0,\be_1]$ and $h^\be_{j\to\rt}(\si^\perm)$ uniformly
bounded away from zero. By \eqref{ehpermlbd} there exists
deterministic $c$ such that
\[
\inf_{\be\in[\be_0,\be_1]} h^\be_{j\to\rt}\bigl(
\si^\perm\bigr)\ge e^{-c(D_j+1)} \qquad\forall j \in\pd\rt, \mu\mbox{-a.s.}
\]
Combining these observations gives $\lim_{m\to\infty}E_{T,m}=0$ $\mu$-a.s.

To take the limit in $\mu$-expectation, we argue similarly as in part
(a): by \eqref{egubd} and~\hypreg\ there exists
deterministic $C'$ such that
\begin{eqnarray*}
\biggl| \f{\pd\Phi_T} {\pd h_{x\to y}(\si)}\bigl(\be+t
\de,h^\be+t\De_{\be
,\de} h\bigr) \biggr| &\le&\f{C'}
{h^\be_{x\to y}(\si^\perm) + t
\De_{\be,\de} h_{x\to
y}(\si ^\perm)}\\
&\le&\sup
_{\be'\in[\be_0,\be_1]} \f{C'} {h^{\be'}_{x\to
y}
(\si ^\perm)}
\end{eqnarray*}
for all $\be\in[\be_0,\be_1-\de]$, $x\to y\in\subt{}{1}$, $\si
\in\spins
$ and $t\in[0,1]$, hence
\[
\abs{E_{T,m}} \le C C' {\sum
_{x\to y}}^* \sup_{\be\in[\be_0,\be_1]} \f{1}
{h^\be_{x\to y}(\si^\perm)}.
\]
This is integrable by \hypintbeta\ and unimodularity of $\mu$, so
dominated convergence implies that
$\lim_{m\to\infty}\E_\mu[E_{T,m}]=0$
as claimed.
\end{longlist}
Combining (a)--(c) gives the
result of the proposition.
\end{pf}

\begin{pf*}{Proof of Theorem~\ref{tfm}(a)}
Recalling Lemma~\ref{lemlogZderiv},
\begin{eqnarray*}
\limsup_{n\to\infty} \bigl[\phi_n(\be_1)-
\phi_n(\be_0)\bigr] &=& \limsup_{n\to\infty}
\int_{\be_0}^{\be_1} a^\edge_n(
\be,B) \d \be
\\
&\le&\int_{\be_0}^{\be_1} \limsup_{n\to\infty}
a^\edge_n(\be,B) \d\be \le\int_{\be_0}^{\be_1}
a^\edge(\be,h)\d\be,
\end{eqnarray*}
where the first inequality follows by (the reversed) Fatou's lemma and
the second one by the hypothesis \eqref{edbetaassump}. By
Proposition~\ref{pbethedbeta} the right-most expression equals to
$\Phi(\be_1,h)-\Phi(\be_0,h)$, so the theorem is proved.
\end{pf*}

The justification for interpolation in $B$ is entirely similar:

\begin{pf*}{Proof of Theorem~\ref{tfm}(b)}
Now $\be$ is fixed, so we suppress it from the notation. For $h\in\cH$
and $T\in\treesed$, then
\begin{eqnarray}
\f{\pd\Phi_T(B,h)} {\pd B} =\f{\pd\Phi_T^\vertex(B,h)}
{\pd B} =\f{ \sum_{\si,\si_j} \pd_B\vxi(\si)
\psi(\si,\si_j) h_{j\to\rt}(\si_j) \wh
h_{\rt\to j}(\si) } { \sum_{\si,\si_j} \psi(\si,
\si_j) h_{j\to\rt}(\si_j) \wh h_{\rt\to j}(\si)
}
\nonumber\\
\eqntext{ \forall j\in\pd\rt,}
\end{eqnarray}
while if $T=\subt{}{0}$, then $\pd_B\Phi_T= \sum_\si\pd_B \vxi
(\si)
\vpsi(\si)/ \sum_\si\vpsi(\si)$. If $h\in\hstar$, then $\wh
h^B=h^B$, so
\[
\E_\mu\bigl[\pd_B\Phi_T(B,h)\bigr] =
\E_\mu\bigl[\bigl[\!\bigl[\pd_B\vxi(\si_\rt
)\bigr]\!\bigr]^{h,B}_T\bigr] \equiv a^\vertex(B,h).
\]
The result now follows by adapting the proofs of Proposition~\ref
{pbethedbeta} and Theorem~\ref{tfm}(a).
\end{pf*}

\subsection{Discussion and first consequences}\label{ssecdiscuss}

We now prove Theorem~\ref{tfmunique} by considering an extended
notion of local weak convergence. As discussed in \cite{MR2354165},
a graph $G=(V,E)$ together with a spin configuration $\vec{\si}\in
\spins
^V$ on the graph can be regarded as a graph with \emph{marks} in
$\spins
$. Let $\graphsrtmark$ and $\graphsprmark$ denote the spaces of marked
isomorphism classes of connected, rooted and bi-rooted graphs,
respectively, with marks in $\spins$. These spaces are metrizable by
the obvious generalizations of the metrics on $\graphsrt,\graphspr$
defined in Section~\ref{ssecprelim}, giving rise to the notion of
local weak convergence for pairs $(G_n,\vec{\si}_n)$ of graphs with spin
configurations. Definition~\ref{dunim} generalizes naturally to this
setting, and we show next that if $\vec{\si}_n$ is a random configuration
on $G_n$ with law $\nu_{G_n,\vec{\psi}}$ [as defined in \eqref{efm}],
then a local weak limit of $(G_n,\vec{\si}_n)$, if it exists, must be
unimodular.

\begin{lem}\label{lemtight}
If $G_n\lwc\mu$ and $\vec{\si}_n\sim\nu_{G_n,\vec{\psi}}$,
then the laws of $(G_n,\vec{\si}_n)$ have subsequential local weak limits
belonging to the space $\unimsr$ of unimodular measures on
$\graphsrtmark$.
\end{lem}

\begin{pf}
For each fixed $t$, the laws of $\ball{t}{I_n}$ are weakly convergent,
hence by Prohorov's theorem
form a uniformly tight sequence. Consequently, for each $\ep>0$ there
exists $K_\ep\subseteq\graphsrt$
compact with $\sup_n \P_n( \ball{t}{I_n} \notin K_\ep) \le\ep$.
Further, $K_\ep$ may be taken to
contain only graphs of depth at most $t$, whereby the minimal distance
between any two graphs in
$K_\ep$ is uniformly bounded below [by $1/(1+t)$], hence the
compactness of $K_\ep$ implies that it
must be a finite set. The collection of all marked graphs in
$\graphsrtmark$ whose underlying graph
is in $K_\ep$ must therefore be finite, hence compact as well. Thus, by
yet another application of
Prohorov's theorem, the joint laws of $(\ball{t}{I_n},\vec{\si
}_{\ball
{t}{I_n}})$ are uniformly tight
in $\graphsrtmark$ and consequently have subsequential weak limits. By
extracting successive
subsequences for increasing $t$ and taking the diagonal subsequence, it
follows that
the sequence $(G_n,\vec{\si}_n)$ admits subsequential local weak limits
$\wh\mu\in\unimsr$.
\end{pf}

For $\wh\mu\in\unimsr$, the marginal $\mu$ of $\wh\mu$ is a unimodular
measure on $\graphsrt$. If it is supported on a single tree $T$
as in the $d$-regular case, then clearly $\wh\mu$
may be represented as $\de_T \times\nu$ where $\nu\in\gibbs_T$,
the space of Gibbs measures on $T$ corresponding to specification $\vec
\psi$.
To make such a statement in the general setting, note that there is a
continuous mapping $\pi$ from $\graphsrt$ to the space $\graphsrtz$
of graphs on $\Z_{\ge0}$ rooted at $0$, taking an isomorphism class
to its
canonical representative (\cite{MR2354165}, page~1461).
Thus $\wh\mu$ may be regarded as a measure on the product space
$\graphsrtz\times\spins^{\Z_{\ge0}}$, and consequently $\wh\mu$
has a representation as the measure
$\mu\otimes\nu$ on pairs $(T,\vec{\si})$
where $T$ has law $\mu$ and $\vec{\si}$ given $T$ has law $\nu_T\in
\gibbs_T$.
In particular, if $\abs{\gibbs_T}=1$ $\mu$-a.s., then $\mu\otimes
\nu$
is uniquely determined.

Let $\mu$ be a unimodular measure on $\treesrt$. It was noted in
Remark~\ref{rbp} that there is a mapping from $\hstar_\mu(\be,B)$ to
collections $(\nu_T\in\gibbs_T)_{T\in\treesrt}$. For such $\nu$,
$\mu
\otimes\nu$ belongs to $\unimsr$:
if $f$ is a nonnegative Borel function on $\graphsprmark$, it follows
from the $\treespr$-measurability of elements of $\hloc$ that
\[
\E_{\mu\otimes\nu} \biggl[ \sum_{j\in\pd\rt} f\bigl((T,
\vec{\si }),\rt,j\bigr) \biggr] =\E_\mu \biggl[ \sum
_{j\in\pd\rt} \bar f(T,\rt,j) \biggr],
\]
where $\bar f$ is a
nonnegative Borel function on $\graphspr$. The unimodularity of the
underlying measure $\mu$ then gives
\[
\E_{\mu\otimes\nu} \biggl[ \sum_{j\in\pd\rt} f\bigl((T,
\vec{\si }),\rt,j\bigr) \biggr] = \E_{\mu\otimes\nu} \biggl[ \sum
_{j\in\pd\rt} f\bigl((T,\vec{\si }),j,\rt\bigr) \biggr],
\]
and therefore $\mu\otimes\nu\in\unimsr$.

\begin{rmk}\label{rmc}
An element $\nu\in\gibbs_T$ is called a \emph{Markov chain} (or
\emph
{splitting Gibbs measure}) if for any finite connected sub-graph
$U\subseteq T$, the marginal of $\nu$ on $U$ is a Markov random field
\cite{MR714953}; see also \cite{MR956646}, Chapter~12, and \cite
{MR0378152}.
A collection $\Lm_T\equiv(\lm^j_i)_{(ij)\in E_T}$ of probability
measures on $\spins$ is called an \emph{entrance law} (or \emph
{boundary law}) for the specification $\vec{\psi}\equiv(\psi,\vpsi
)$ on
$T$ if it satisfies the consistency requirement (\cite{MR714953}, (3.4))
\[
\lm^j_i(\si_i) = \prod
_{k\in\pd i\setminus j} \biggl( \sum_{\si_k}
\phi_{ik}(\si_i,\si_k) \lm^i_k(
\si_k) \biggr),
\]
where $\phi_{ij}(\si_i,\si_j)\equiv\vpsi(\si)^{1/D_i}\psi(\si
,\si')\vpsi
(\si')^{1/D_j}$, the pairwise interaction potential corresponding to
$\vec{\psi}$. It is shown in \cite{MR714953}, Theorem~3.2, that
there is a
one-to-one correspondence between Markov chains $\nu$ and entrance laws
$\Lm_T$, given by
\[
\nu(\vec{\si}_U) \cong\prod_{(ij)\in E_U}
\phi_{ij}(\si_i,\si_j) \prod
_{i\in\pd U} \biggl( \sum_{\si_i}
\phi_{i p(i)}(\si_i,\si_{p(i)}) \lm
^{p(i)}_i(\si _i) \biggr)
\]
for $U$ any finite connected sub-graph of $T$, with $p(i)$ denoting the
unique neighbor of $i$ inside $U$ for $i\in\pd U$. In particular, we
see that the Gibbs measure $\nu_T$ arising from $h\in\hstar_\mu(\be,B)$
is precisely the Markov chain with entrance law $\lm^i_j(\si)\cong
h_{j\to i}(\si)\vpsi(\si)^{-1/D_j}$. Extremal elements of $\gibbs_T$
are Markov chains (\cite{MR714953}, Theorem~2.1), but the converse is
false; for example, the free-boundary Ising Gibbs measure is
nonextremal at low temperature; see \cite{MR1768240,MR0676482}. The
measures $\mu\otimes\nu$ arising from elements of $\hstar_\mu(\be,B)$
might naturally be termed ``unimodular Markov chains'' or
``Bethe Gibbs measures,'' in the sense that the entrance laws for the
entire collection $(\nu_T)_{T\in\treesrt}$ are specified by a single
measurable function $h\dvtx \treespr\to\simplex$ which is a Bethe fixed
point. In the case $\mu=\de_{\treereg{d}}$ these correspond precisely
to the \emph{completely homogeneous Markov chains} studied in
\cite{MR714953}, Section~4.
\end{rmk}

\begin{pf*}{Proof of Theorem~\ref{tfmunique}}
Suppose uniqueness holds at $(\be,B)$, that is, $\gibbs_T=\set{\nu_T}$
$\mu$-a.s. Then $\hstar_\mu(\be,B)$ has size at most one by
Remark~\ref
{rembed}. For $\mu$-a.e. $T$, the measure $\nu_T$ is extremal, and so
specifies a Markov chain on $T$ with entrance law $\Lm_T$; see
Remark~\ref{rmc}. If we define $h_{x\to y}(\si) \equiv h_{(T,x\to
y)}(\si) \cong\lm^y_x(\si) \vpsi(\si)^{1/D_x}$, then $h\in\hstar
_\mu
(\be,B)$, which proves that $\hstar_\mu(\be,B)$ is a singleton.

Now consider interpolation in $\be$ or $B$. All the conditions of
Theorem~\ref{tfm} are satisfied by assumption except \eqref
{edbetaassump} and \eqref{edBassump}. If uniqueness holds at $(\be
,B)$, it follows from the preceding discussion that there is a unique
$\mu\otimes\nu\in\unimsr$ corresponding to the specification
$(\psi^\be
,\vpsi^B)$. Any local weak limit of $(G_n,\vec{\si}_n)$ must be such a
measure, so $(G_n,\vec{\si}_n)\lwc\mu\otimes\nu$; likewise, any element
of $\hstar_\mu(\be,B)$ gives rise to $\mu\otimes\nu$. Therefore,
\[
\lim_{n\to\infty} a^\edge_n(\be,B) = \f{1}
{2} \E_{\mu\otimes\nu} \biggl[ \sum_{j\in\pd\rt}
\pd_\be \xi(\si_\rt,\si_j) \biggr] =
a^\edge(\be,B,h),
\]
where the limit in expectation is justified by the boundedness of $\pd
_\be\xi$ on compacts and uniform sparsity (as in the proof of
Lemma~\ref
{lemlogZderiv}). This verifies \eqref{edbetaassump}, and the
verification of \eqref{edBassump} is entirely similar. The result
therefore follows from Theorem~\ref{tfm}.
\end{pf*}

\begin{rmk}
If uniqueness of Gibbs measures does not hold, one may consider
extremal decomposition of the subsequential local weak limits $\wh\mu$
of $(G_n,\vec{\si}_n)$, either in the spaces $\gibbs_T$ (possibly losing
unimodularity in the decomposition), or in the space $\unimsr$.
Extremal decomposition in $\unimsr$ is discussed in \cite{MR2354165}, Section~4, but it is unclear whether extremal elements would be
unimodular Markov chains in the sense described here. A decomposition
of $\wh\mu=\mu\otimes\nu$ into unimodular Markov chains $\mu
\otimes\nu
'$ would obviously yield a substantial generalization of Theorem~\ref
{tfmunique}.
\end{rmk}

\subsection{Application to independent set}
\label{ssecis}

We now prove Theorem~\ref{tis}, our result for the independent set
model \eqref{eis}, by verifying the conditions of Theorem~\ref
{tfmunique} for the interpolation parameter $B\equiv\log\lm$. In this
setting a convenient parametrization for the messages $h\in\cH$ is
$u\equiv h(0)$, so that the BP mapping \eqref{ebp} becomes
%
\begin{equation}
\label{ebpis} \bigl(\BP^\lm u\bigr)_{x\to y} = \f{1} { 1 + \lm
\prod_{v\in\pd x\setminus y} u_{v\to x}}.
\end{equation}
A single BP iteration is anti-monotone in the messages $u_{v\to x}$, so
a double iteration is monotone. Since the root marginal for an
independent set model in $\subt{}{2t-1}$ is obtained by an even number
of BP iterations starting from level $2t$ (see Remark~\ref{rbp}), it
is monotone in the boundary conditions. Recalling from Section~\ref{sssecintrois} the definition of $\vh^{t,\ddagger}_T\equiv\vh
^{t,\ddagger,\lm}_T$ for $\ddagger\in\set{0,1}$ and writing $\vu
^{t,\ddagger}_T\equiv\vh^{t,\ddagger}_T(0)$, the above implies that
for $1\le s\le t$,
\[
\vu^{2s-1,0}_T \ge\vu^{2t-1,0}_T \ge
\vu^{2t-1,1}_T \ge\vu^{2s-1,1}_T \ge
\vu^{1,1}_T = \f{1} {1+\lm}.
\]
Thus the $t\to\infty$ limits $\vh^0_T,\vh^1_T$ are well-defined with
$\vh^0_T(1)\ge\vh^1_T(1)\ge1/(1+\lm)$, and using these we define
messages $h^\ddagger\in\cH$, $h^\ddagger_{x\to y}=\vh^\ddagger
_{\subt
{x\to y}{}}$. The next lemma gives the boundary values for the interpolation.

\begin{lem}\label{lemisbd}
For the independent set model on $G_n\lwc\mu$,
\[
\lim_{\lm\decto0}\limsup_{n\to\infty}\bigl|
\phi_n(\lm)\bigr|=0=\lim_{\lm\decto
0}\Phi\bigl(
\lm,h^\ddagger\bigr),\qquad \ddagger\in\set{0,1}.
\]
\end{lem}

\begin{pf}
The left limit follows from the trivial bounds $1\le Z_n\le(1+\lm)^n$.
Next, for any $h\in\cH$,
\[
\Phi^\vertex_T(\lm,h) = \log \biggl\{ 1 + \lm\prod
_{j\in\pd\rt} h_{j\to\rt}(0) \biggr\},
\]
so $\Phi^\vertex_T(\lm,h)\to0$ both $\mu$-a.s. and in $\mu
$-expectation as $\lm\decto0$, by bounded convergence. The same holds
for $\Phi^\edge_T(\lm,h^\ddagger)$, $\ddagger\in\set{0,1}$,
using the bound
$h^\ddagger_{x\to y}(0) \ge1/(1+\lm)$.
\end{pf}

\begin{pf*}{Proof of Theorem~\ref{tis}}
The independent set model \eqref{eis} is of form \eqref{efm} with
$\spins=\set{0,1}$, $\psi(\si,\si')=\Ind{\si\si'\ne1}$, and
$\vpsi(\si
)=\lm^\si\equiv e^{B\si}$, so \hypreg\ is clearly satisfied with
$\si
^\perm=0$ the permitted state. By definition of $\lm_c$, if $\lm<\lm
_c$, then $h^0=h^1\equiv h$ in $\cH$, and it then follows from the
recursive structure of the tree that $h\in\hstar_\mu(\lm)$. Since
$h^\ddagger_{x\to y}(0)\ge1/(1+\lm)$ as noted above, \hypintB\ is
satisfied on any compact interval of $\lm$.

For $T\in\treesrt$, as noted above the root occupation probability on
$\subt{}{s}$ for $s\ge2t-1$ with \emph{any} boundary conditions is
sandwiched between $\vh^{2t-1,0}_T(1)$ and $\vh^{2t-1,1}_T(1)$, with
the former increasing to $\vh^0_T(1)$ and the latter decreasing to
$\vh
^1_T(1)$. Since the $\vh^{t,\ddagger}_T$ are clearly continuous in
$\lm
$, it follows that $\vh^0_T(1)$ and $\vh^1_T(1)$ are, respectively, lower
and upper semi-continuous in $\lm$, so if they coincide, then their
common value $\vh_T(1)$ is continuous in $\lm$. Applying this with
$T=\subt{x\to y}{}$ gives the $\mu^\indir$-a.s. continuity of
$h^\ddagger_{x\to y}$ on $(0,\lm_c)$.

For $T\in\treesrt$, $\vh^\ddagger_T$ for $\ddagger\in\set{0,1}$
is a
function of $(h^{1-\ddagger}_{j\to\rt})_{j\in\pd\rt}$, so for
$\lm<\lm
_c$ we have that $\vh^0_T=\vh^1_T$, $\mu$-a.s. It then follows from the
preceding observations and Remark~\ref{rbp} that the boundary effect
vanishes and $\abs{\gibbs_T}=1$ $\mu$-a.s. Thus, we are in the setting
of Theorem~\ref{tfmunique}(b), and it remains only
to complete the verification of \hypdiffB, that is, the boundedness in
total variation of the messages $h_{x\to y}$:

\begin{longlist}[(a)]
\item[(a)] No verification is needed since boundedness in total
variation is simply assumed.

\item[(b)] For $T\in\treesrt$, $\vu^{t,\ddagger}_T\equiv\vh
^{t,\ddagger}_T(0)$ satisfies
\[
\log\vu^{2t+1,\ddagger}_T = -\log \biggl( 1 + \lm\prod
_{j\in\pd\rt} \f{1} {1 + \lm\prod_{k\in\pd
j\setminus\rt}
\vu^{2t-1,\ddagger}_{T_{k\to j}}} \biggr).
\]
Differentiating with respect to $\lm$, we find that $r^{t,\ddagger
}_T\equiv(1+\lm) \pd_\lm\log\vu^{t,\ddagger}_T$ satisfies
%
\begin{equation}
\label{eisderiv} \bigl|r^{2t+1,\ddagger}_T\bigr| \le1 + \f{\lm} {1+\lm}
D_\rt + \biggl(\f{\lm} {1+\lm} \biggr)^2 \sum
_{j\in\pd\rt} \sum_{k\in
\pd
j\setminus\rt}
\bigl|r^{2t-1,\ddagger}\bigr|_{T_{k\to j}}.
\end{equation}
Since $\vu_T^{1,1}=1/(1+\lm)$ for any $T\in\treesrt$, we find that
\[
\sup_{t\ge0} \bigl|r_T^{2t-1,1}\bigr| \le1+\sum
_{\ell\ge1} \biggl(\f {\lm} {1+\lm } \biggr)^{\ell} \bigl|\pd
\subt{} {\ell-1}\bigr|.
\]
If $\lm_1/(1+\lm_1)< 1/\br T$, then this is finite and uniformly
bounded on $[\lm_0,\lm_1]$ (see~\eqref{ebranch} or \cite{MR1062053}, Section~2), and consequently $\vu^1_T\equiv\lim_{t\to\infty
}\vu
^{2t-1,1}_T$ has deterministically bounded total variation on $[\lm
_0,\lm_1]$.
If $\lm_1<\lm_{c,\mu}$, then $h^{2t-1,1}_{x\to y} \to h_{x\to y}$ on
$[\lm_0,\lm_1]$, so if
$\br T_{x\to y}\le\De-1$ $\mu^\indir$-a.s. and $\lm_1/(1+\lm
_1)<1/(\De
-1)$ [i.e., $\lm_1(\De-2)<1$], then $h$ has deterministically bounded
total variation on $[\lm_0,\lm_1]$.

\item[(c)]
Since the limiting measure is supported on $\treereg{d}$, only
$h\equiv
h_{(\treereg{d},x\to y)}$ is of relevance,
and \eqref{ebpis} reduces to $\BP^\lm u = (1+\lm u^{d-1})^{-1}$. For
$\lm\le\lm_c=\lm_c(d)$
there is a unique fixed point (see \cite{MR844469}, Section~2), which is
then easily seen to be monotone in $\lm$.

Thus \hypdiffB\ is verified in parts (a)--(c). Also, $\phi(\lm_c)=\lim_{\lm\incto\lm_c}\phi(\lm)$
as an
immediate consequence of Lemma~\ref{lemlogZderiv}. The rest of the
theorem follows by applying Theorem~\ref{tfmunique} and then taking
$B_0 =\log\lm_0 \to-\infty$, relying on the boundary value given by
Lemma~\ref{lemisbd}.\quad\qed
\end{longlist}
\noqed\end{pf*}

\section{Bethe prediction as optimization over local polytope}\label
{sechloc}

Throughout this section we assume that $\vec{\psi}\equiv(\psi,\vpsi)$
satisfying \hypreg\ specifies a factor model \eqref{efm}, and that
$\mu
$ is a unimodular measure on $\treesrt$ with $\E_\mu[D_\rt]<\infty
$. We
study the Bethe prediction as the optimization of the Bethe free energy
functional $\Phi_\mu$ on $\hloc$ as defined by \eqref
{ebethehlocfirst}. We first verify that this agrees with the
previous definition~\eqref{ebethehstar} of $\Phi_\mu$ on $\hstar
_\mu
(\be,B)$, which we always regard as being embedded into $\hloc$ via~\eqref{eembed}. Recall from Definition~\ref{dhloc} that for $\bh
\in
\hloc$, the one-point marginals of $\bh_{xy}$ are denoted $\vh_x$ and
$\vh_y$, and are measurable functions $\treesrt\to\simplex$.

\begin{lem}\label{lembetheext}
The functional $\Phi_\mu$ on $\hloc$ given by \eqref
{ebethehlocrelent} agrees with the previous definition \eqref
{ebethehstar} on $\hstar_\mu(\be,B)$, subject to finiteness of $\E
_\mu
[\Phi^\edge_T]$.
\end{lem}

\begin{pf}
If $\bh$ corresponds to $h\in\hstar_\mu(\be,B)$, then \eqref{eembed}
and \eqref{ehrecurs} imply that
\begin{eqnarray*}
\bh_{xy}\bigl(\si,\si'\bigr) \exp\bigl\{
\Phi^{(xy)}_T(h)\bigr\} &=& \psi\bigl(\si,\si'
\bigr)h_{x\to y}(\si) h_{y\to x}\bigl(\si'\bigr),
\\
\vh_\rt(\si) \exp\bigl\{\Phi^\vertex_T(h)\bigr
\} &=& \vpsi(\si) \prod_{j\in\pd\rt} \biggl( \sum
_{\si_j} \psi(\si,\si_j) h_{j\to\rt}(
\si_j) \biggr).
\end{eqnarray*}
Letting $\Phi^{(i)}(\bh)$ ($1\le i\le3$) denote the three terms on the
right-hand side of \eqref{ebethehloc}, it follows from the above that
\begin{eqnarray*}
\Phi^{(1)}(\bh) &=& \E_\mu\bigl[\Phi^\vertex_T(h)
\bigr] - \E_\mu \biggl[ \sum_{j\in\pd\rt} \sum
_\si \vh_\rt(\si) \log \biggl( \sum
_{\si_j} \psi(\si,\si_j) h_{j\to\rt}(
\si_j) \biggr) \biggr],
\\
\Phi^{(2)}(\bh) &=& \E_\mu\bigl[\Phi^\edge_T(h)
\bigr] - \f{1} {2} \E_\mu \biggl[ \sum_{j\in\pd\rt}
\sum_{\si,\si_j} \bh_{\rt j}(\si,
\si_j) \log \bigl( h_{\rt\to j}(\si) h_{j\to\rt}(
\si_j) \bigr) \biggr]
\\
&=& \E_\mu\bigl[\Phi^\edge_T(h)\bigr] -
\E_\mu \biggl[ \sum_{j\in\pd\rt}\sum
_\si \vh_\rt(\si) \log h_{\rt\to j}(\si)
\biggr],
\\
\Phi^{(3)}(\bh) 
&=&\E_\mu \biggl[ \sum
_{j\in\pd\rt} \sum_\si
\vh_\rt(\si)\log \biggl( \sum_{\si_j} \psi(
\si,\si_j) h_{\rt\to j}(\si) h_{j\to\rt}(\si_j)
\biggr) \biggr]\\
&&{} - 2\E_\mu\bigl[\Phi^\edge_T(h)
\bigr],
\end{eqnarray*}
where unimodularity was used in the simplification of $\Phi^{(2)}$.
Adding these three identities gives
$\Phi_\mu(\bh)=\E_\mu[\Phi^\vertex_T(h)-\Phi^\edge_T(h)]$,
as claimed.
\end{pf}

As mentioned in Section~\ref{sssecintroopt}, our definition $\Phi
_\mu
$ of the Bethe free energy functional on $\hloc$ is an infinite-tree analogue
of the definition of \cite{MR2246363} for finite graphs. It is proved in
\cite{MR2246363}, Proposition~6, that when $\psi>0$,
all local maxima of the Bethe free energy lie in the \emph{interior} of
the local polytope. We now prove an analogous result for infinite
unimodular trees, assuming only permissivity of $\vec{\psi}$.

\begin{ppn}\label{pinter}
For permissive $\vec{\psi}$, if $\bh$ is a local maximizer of $\Phi
_\mu$
over $\hloc$, then $\bh\in\hlocopsi$.
\end{ppn}

\begin{pf}
Assume without loss that $\bh\in\hloc[\psi]$, since otherwise clearly
$\Phi_\mu(\bh)=-\infty$. If $\pr{u}\in\hloc[\psi]$, then it
follows by
convexity of $\hloc$ that
$\bh^\eta
\equiv\bh+\eta(\pr{u}-\bh)
\equiv\bh+\eta\prr{\de}$
belongs to $\hloc[\psi]$ for any $\eta\in(0,1]$. Letting
\[
R^\eta(\prr{\de}) \equiv\f{2} {\eta} \bigl[\Phi_\mu\bigl(
\bh^\eta\bigr) - \Phi _\mu(\bh )\bigr],\qquad \wh
R^\eta(\prr{\de})\equiv\f{R^\eta(\prr{\de})} {\abs{\log
\eta}},
\]
our claim will follow upon showing that if $h\notin\hlocopsi$, then
there exists such $\pr{u}$ for which
\[
\lim_{\eta\decto0} \wh R^\eta(\prr{\de}) = \wh
R^0(\prr{\de})>0.
\]
To this end, note that by an easy computation
$[H(\bh^\eta)-H(\bh)]/\eta
=\break -\anglb{\log\bh^\eta}_{\prr{\de}} - \anglb{f^\eta(\prr{\de
}/\bh)}_{\bh}$,
where $f^\eta(r)\equiv\eta^{-1}\log(1+\eta r)$ and $(\prr{\de}/\bh
)(\si
,\si')$ is defined to be $\prr{\de}(\si,\si')/\bh(\si,\si')$ if
$\bh(\si
,\si')>0$, zero otherwise; note that $\pr{u}=\bh+\prr{\de}\ge0$ implies
$\prr{\de}/\bh\ge-1$. Thus from \eqref{ebethehlocfirst} we obtain
$R^\eta(\prr{\de})= R^\eta_1(\prr{\de}) + R^\eta_2(\prr{\de})$ where
%
%
\begin{eqnarray}\label{eReta}
R^\eta_1(\prr{\de}) &\equiv&\E_\mu
\biggl[ 2\angl{\vxi}_{\vde_\rt} +2(D_\rt-1) \bigl\langle f^\eta(
\vde_\rt/\vh_\rt)\bigr\rangle_{\vh_\rt} \nonumber\\
&&\hspace*{19pt}{}+\sum
_{j\in\pd\rt} \bigl( \angl{\xi}_{\prr{\de}_{\rt j}} -
\bigl\langle f^\eta(\prr{\de}_{\rt j} /\bh_{\rt j})\bigr\rangle_{\bh_{\rt j}}
\bigr) \biggr],
\\
 R^\eta_2(\prr{\de}) &\equiv&
\E_\mu \biggl[ 2(D_\rt-1)\bigl\langle{\log\vh^\eta_\rt}\bigr\rangle_{\vde_\rt}
-\sum_{j\in\pd\rt} \bigl\langle\log\bh^\eta_{\rt j}\bigr\rangle_{\prr{\de}_{\rt j}}
\biggr].\nonumber
\end{eqnarray}
Since for $r\ge-1$ and $\eta\in(0,1)$, we have $\eta^{-1}\log
(1-\eta
)\le f^\eta(r)\le r$, it follows from dominated convergence (and the
boundedness of $\xi$ on $\supp\prr{\de}$) that $R^\eta_1(\prr{\de})$
converges to a finite limit as $\eta\decto0$, and so converges to zero
upon rescaling by $\abs{\log\eta}$. Again by dominated convergence,
$R^\eta_2(\prr{\de})/\abs{\log\eta}$ converges as $\eta\decto0$ to
%
\begin{eqnarray}
\label{eRzero} \wh R^0(\prr{\de})& =&\E_\mu \biggl[ (2-2
D_\rt) \vu_\rt\bigl(\bigl\{\si\dvtx \vh_\rt(
\si)=0\bigr\}\bigr)
\nonumber
\\[-8pt]
\\[-8pt]
\nonumber
&&\hspace*{18pt}{}+ \sum_{j\in\pd\rt} \pr{u}_{\rt j}
\bigl(\bigl\{\si,\si' \dvtx \bh_{\rt
j}\bigl(\si,\si
'\bigr)=0 \bigr\}\bigr) \biggr].
\end{eqnarray}
Let $A_\rt\equiv A_{(T,\rt)}\equiv\set{\si\in\spins\dvtx  \vh_\rt
(\si)=0}$.
Since $\bh_{xy}(\si,\si')=0$ whenever either $\si\in A_x$
or $\si' \in A_y$, we have by unimodularity of $\mu$ that
\begin{eqnarray*}
\wh R^0(\prr{\de}) &\ge& \E_\mu \biggl[ (2-2
D_\rt) \vu_\rt(A_\rt) + \sum
_{j\in\pd\rt} \bigl\{ \vu_\rt(A_\rt) +
\vu_j(A_j) - \pr{u}_{\rt j} (A_\rt
\times A_j) \bigr\} \biggr]
\nonumber
\\
&=& \E_\mu \biggl[2 \vu_\rt(A_\rt)- \sum
_{j\in\pd\rt} \pr {u}_{\rt j} (A_\rt
\times A_j) \biggr] = \E_\mu \biggl[ \sum
_{j\in\pd\rt} \wh R_{\rt\to j} \biggr],
\end{eqnarray*}
where $\wh R_{\rt\to j}\equiv
\Ind{D_\rt>0} [D_\rt^{-1} \vu_\rt(A_\rt) + D_j^{-1} \vu_j(A_j) -
\pr
{u}_{\rt j} (A_\rt\times A_j)]$ [by \eqref{evhdegzero}, necessarily
$A_\rt= \varnothing$ when $D_\rt=0$].

Noting that $A_\rt^c\ne\varnothing$, consider the measurable function
$\vu
\dvtx \treesed\to\simplex$ defined (up to $\mu$-equivalence) by
%
\begin{eqnarray}
\label{emarg} \vu_\rt(\si) &\equiv&\vu_{(T,\rt)}(\si)
\nonumber
\\[-8pt]
\\[-8pt]
\nonumber
& \equiv&
\cases{ \displaystyle\mathbf{1}\bigl\{\si=\si^\perm\bigr\}, & \quad $A_\rt^c=
\bigl\{\si^\perm\bigr\},$\vspace*{2pt}
\cr
\displaystyle \f12 \biggl(\mathbf{1}\bigl\{\si=
\si^\perm\bigr\} +\f{\Ind{\si\in A_\rt^c\setminus
\set{\si^\perm}}} {\abs{A_\rt^c\setminus\set{
\si^\perm}}} \biggr), &\quad $\mbox{else.}$ } %
\end{eqnarray}
Among those $\pr{u}\in\hloc$ with support contained in $\set{(\si
,\si
')\dvtx \si^\perm\in\set{\si,\si'}}$, there is a unique one with marginals
\eqref{emarg}. On the event $\set{D_\rt>0}$, we have the following:
\begin{itemize}
\item[--] If $\si^\perm\in A_\rt\cap A_j$, then
\[\!\!
\pr{u}_{\rt j}\bigl(\si,\si'\bigr) = \f12 \biggl( \mathbf{1}\bigl\{\si=
\si^\perm\bigr\} \f{\Ind{\si'\in A_j^c
\setminus\set{\si^\perm}}} {\abs{ A_j^c
\setminus\set{\si^\perm} }} + \mathbf{1}\bigl\{\si'=
\si^\perm\bigr\} \f{\Ind{\si\in A_\rt^c\setminus\set{
\si^\perm}}} {\abs{ A_\rt^c\setminus\set{
\si^\perm} }} \biggr),
\]
so $\wh R_{\rt\to j}=(2D_\rt)^{-1}+(2D_j)^{-1}$.
\item[--] If $\si^\perm\in A_\rt\cap A_j^c$, then $\vu_\rt(A_\rt
)\ge
1/2$ while $\vu_j(A_j)=0=\pr{u}_{\rt j}(A_\rt\times A_j)$, so $\wh
R_{\rt\to j}\ge(2D_\rt)^{-1}$. Symmetrically if $\si^\perm\in
A_\rt
^c\cap A_j$, then $\wh R_{\rt\to j}\ge(2D_j)^{-1}$.
\item[--] If $\si^\perm\notin A_\rt\cup A_j$, then $\wh R_{\rt\to j}=0$.
\end{itemize}
Thus $\wh R^0(\prr{\de})\ge0$, with strict inequality unless $\si
^\perm
\notin A_\rt\cup A_j$ $\mu^\indir$-a.s., in which case we take
$A_\rt
,A_j$ in place of $A_\rt^c,A_j^c$ in \eqref{emarg}. Then
\[
\wh R_{\rt\to j} = (2 D_\rt)^{-1}
\Ind{A_\rt\ne\varnothing} + (2 D_j)^{-1}
\Ind{A_j \ne\varnothing},
\]
so $\wh R^0(\prr{\de})>0$ unless $\mu(A_\rt=\varnothing)=1$. But in this
case taking $\pr{u}\in\hloc$ identically equal to the uniform measure
on $\supp\psi$ gives
\[
\wh R^0(\prr{\de}) = \f{1} {\abs{\supp\psi}} \E_\mu
\biggl[\sum_{j\in\pd\rt} \bigl|(\supp\psi)\setminus(\supp
\bh_{\rt j})\bigr| \biggr].
\]
If $\bh\notin\hlocopsi$, then this is positive, completing the proof of
our claim.
\end{pf}

Our main result in this section is the following infinite-tree analogue
of \cite{MR2246363}, Theorem~2, characterizing the interior stationary
points of $\Phi_\mu$ as fixed points of the Bethe recursion.

\begin{ppn}\label{popt}
For $\vec{\psi}$ permissive, any stationary point of $\Phi_\mu$ inside
$\hlocopsi$ belongs to $\hstar$.
\end{ppn}

\begin{pf}
Let $\hlocpmpsi$ denote the space of measurable functions $\prr{\de
}\dvtx \treespr\to\R^{\spins^2}$ (defined up to $\mu^\indir$-equivalence)
such that $\supp\prr{\de}_{xy}\subseteq\supp\psi$, $\prr{\de
}_{xy}(\si,\si
')=\prr{\de}_{yx} (\si',\si)$, the one-point marginals $\vde_x(\si
)\equiv\sum_{\si'} \prr{\de}_{xy}(\si,\si')$ do not depend on the
choice of $y\in\pd x$, and $\sum_\si\vde(\si)=\sum_{\si,\si
'}\prr{\de
}(\si,\si')\equiv0$.

\textit{Step} 1. We first show that if $\bh\in\hlocopsi$ is a stationary
point of $\Phi_\mu$, then there exists $\lm\dvtx \treespr\to\R^{\spins}$
measurable such that
%
\begin{equation}
\label{estatfactors} \bh_{xy}\bigl(\si,\si'\bigr) =\psi
\bigl(\si,\si'\bigr)\exp\bigl\{\lm_{x\to y}(\si) +
\lm_{y\to x}\bigl(\si'\bigr)\bigr\}, \qquad\mu^\indir
\mbox{-a.s.}
\end{equation}
Since $\bh\in\hlocopsi$, if $\prr{\de}\in\hlocpmpsi$ with $\abs
{\prr{\de
}} \le\bh$ $\mu^\indir$-a.s., then $\bh^\eta\equiv\bh+\eta\pr
{\de}$
belongs to $\hloc[\psi]$ for all $\abs{\eta}\le1$. Taking $\eta
\to0$ in
\eqref{eReta} gives (by stationarity of $\Phi_\mu$ at $\bh$)
\[
0=R^0(\prr{\de}) =\E_\mu \biggl[ 2\bigl\langle{
\vka_\rt'}\bigr\rangle_{\vde_\rt} +\sum
_{j\in\pd\rt} \bigl\langle{\prr{\ka}'}\bigr\rangle_{\prr{\de}}
\biggr],
\]
where $\vka'_x\equiv\vxi+(D_x-1)\log\vh_x$,
$\prr{\ka}'_{xy}\equiv(\xi-\log\bh_{xy})\I_{\supp\psi}$.

Consider now $\prr{\de}$ with one-point marginals $\vde\equiv0$, so
that the value of $\vka'$ becomes irrelevant: in this case the value of
$R^0(\prr{\de})$ is unchanged upon replacing $\prr{\ka}'$ by
\[
\prr{\ka}_{xy}\bigl(\si,\si'\bigr) \equiv
\I_{\supp\psi}\bigl(\si,\si'\bigr) \bigl[\prr{
\ka}'_{xy}\bigl(\si,\si'\bigr) +
\lm_{x\to y}(\si) + \lm_{y\to x}\bigl(\si'\bigr)
\bigr].
\]
We claim it is possible to choose $\lm$ such that $\prr{\ka}$ has
one-point marginals $\vka\equiv0$, $\mu^\indir$-a.s. This amounts to
solving the linear system
%
\begin{equation}
\label{elinsys} \pmatrix{ a_{x\to y} \vspace*{2pt}
\cr
a_{y\to x} }
= \pmatrix{ I & Q \vspace*{2pt}
\cr
Q & I } \pmatrix{\lm_{x\to y}
\vspace*{2pt}
\cr
\lm_{y\to x} } \equiv\mathbf{Q} \pmatrix{
\lm_{x\to y}\vspace*{2pt}
\cr
\lm_{y\to x}},
\end{equation}
where, writing $r(\si)\equiv\abs{\set{\si'\dvtx \psi(\si,\si')>0}}$,
\[
a_{x\to y}(\si) \equiv-\f{\sum_{\si'} \prr{
\ka}'_{xy}(\si,\si')} {r(\si)},\qquad Q
\bigl(\si,\si'\bigr)\equiv\f{\I_{\supp\psi}(\si,
\si')} {r(\si)}.
\]
For $\vec{\psi}$ permissive, the Markov kernel $Q$ is irreducible and
aperiodic, with stationary distribution $\vec r\equiv(r(\si))_\si$ (by
symmetry of $\psi$). By the Perron--Frobenius theorem, $Q,Q^2$ both
have unique left eigenvector $\vec r$ corresponding to eigenvalue $1$.
Therefore $\dim\ker(I-Q^2)=1$, from which it is easy to see that
$\ker
\mathbf{Q}^t=(\im\mathbf{Q})^\perp$ is the linear span of $(\vec
r,-\vec r)$. Since the assumed symmetry properties of $\psi$ and $\bh$
imply that
\[
\bigl\langle(\vec r,-\vec r), (a_{x\to y},a_{y\to x})\bigr\rangle =\sum
_{\si,\si'} \bigl(-\prr{\ka}'_{xy}\bigl(
\si,\si'\bigr) +\prr{\ka}'_{yx}\bigl(\si,
\si'\bigr)\bigr) =0, \qquad \mu^\indir\mbox{-a.s.,}
\]
there is a unique solution $(\lm_{x\to y},\lm_{y\to x})$ to the system
\eqref{elinsys} giving the required solution to \eqref{estatfactors}.

For this choice of $\prr{\ka}$, $\prr{\de}=c\prr{\ka}$ belongs to
$\hlocpmpsi$ for any measurable $c\dvtx \treespr\to\R_{>0}$ with
$c_{xy}=c_{yx}$. We can choose $c$ small enough so that $\abs{\prr{\de
}}<\abs{\bh}$ on $\supp\psi$ $\mu^\indir$-a.s. With this choice,
$0=R^0(\prr{\de})$ becomes the $\mu$-expectation of a (weighted) sum of
squares, so $\prr{\ka}\equiv0$, and rearranging gives \eqref{estatfactors}.

\textit{Step} 2. Returning now to general $\prr{\de}\in\hlocpmpsi$ with
$\abs{\prr{\de}} \le\bh$ $\mu^\indir$-a.s., we obtain from \eqref
{estatfactors} the simplification
%
\begin{eqnarray}
\label{eRzerosimp} 0&=&R^0(\prr{\de}) =\E_\mu \biggl[ 2
\bigl\langle\vka_\rt'\bigr\rangle_{\vde_\rt} -\sum
_{j\in\pd\rt} \bigl(\angl{\lm_{\rt\to j}}_{\vde_\rt} +\angl{
\lm_{j\to\rt}}_{\vde_j} \bigr) \biggr]
\nonumber
\\[-8pt]
\\[-8pt]
\nonumber
 &=&2 \E_\mu \biggl[
\biggl\langle\vka'_\rt-\sum
_{j\in\pd\rt} \lm_{\rt\to j} \biggr\rangle_{\vde_\rt}
\biggr],
\end{eqnarray}
using unimodularity of $\mu$ for the last identity. We claim that
%
\begin{eqnarray}
\label{eequalsavg} \vde'_x(\si) &\equiv&
\vka'_x(\si)-\sum_{y\in\pd x}
\lm_{x\to y}(\si) - \f{1} {\abs{\spins}} \sum_{\si'}
\biggl( \vka'_x\bigl(\si'\bigr)-\sum
_{y\in\pd x} \lm_{x\to y}\bigl(\si'
\bigr) \biggr)
\nonumber
\\[-8pt]
\\[-8pt]
\nonumber
&=&0,\qquad \mu^\indir\mbox{-a.s.}
\end{eqnarray}
Indeed, for any $\vde' \dvtx \treesrt\to\R^\spins$ measurable with
$\sum_\si
\vde'_\rt(\si)\equiv0$ $\mu$-a.s.,
\[
\prr{\de}'_{xy}\bigl(\si,\si'\bigr)
\equiv\vde'_x(\si)\mathbf{1}\bigl\{\si'=
\si^\perm\bigr\} +\vde'_y\bigl(\si'
\bigr)\mathbf{1}\bigl\{\si=\si^\perm\bigr\}
\]
defines an element of $\hlocpmpsi$. By considering \eqref
{eRzerosimp} with $\prr{\de}=c\prr{\de}'$ where $c_{xy}=c_{yx}$ is
small enough so that $\abs{c\prr{\de}'}<\abs{\bh}$, we obtain the claim
\eqref{eequalsavg}.

\textit{Step} 3.
Rearranging \eqref{eequalsavg} we find that $\bh$ satisfies $\mu
^\indir$-a.s.
%
%
\begin{eqnarray}
\label{ebarh-form} \bh_{\rt j}\bigl(\si,\si'\bigr) &=& \psi
\bigl(\si,\si'\bigr) \exp\bigl\{ \lm_{\rt\to
j}(\si)+\lm
_{j\to\rt}\bigl(\si'\bigr)\bigr\},
\\
\label{eh-form} \vh_\rt(\si) &\cong&\exp \biggl\{ \f{\sum
_{j\in\pd\rt} \lm_{\rt\to j}(\si) - \vxi(\si )}
{D_\rt-1} \biggr\}.
\end{eqnarray}
If we then re-parametrize
%
\begin{equation}
\label{elm-form} \lm_{\rt\to j} \equiv\vxi+ \sum
_{k \in\pd\rt\setminus j} \log \wh m_{k\to\rt}, \qquad\mu^\indir
\mbox{-a.s.}
\end{equation}
(well defined, for each $T$ and $\si\in\spins$, by invertibility of the
$D_\rt$-dimensional matrix $\I\I^t - I$), then formula \eqref
{eh-form} for $\vh_\rt$ becomes
\[
\vh_\rt(\si) \cong\vpsi(\si) \prod_{k\in\pd\rt}
\wh m_{k\to
\rt}(\si ),\qquad \mu^\indir\mbox{-a.s.}
\]
On the other hand, $\vh_\rt$ is the first marginal of $\bh_{\rt j}$,
and setting the above equal to the sum of \eqref{ebarh-form} over
$\si
'$ gives
[making use of \eqref{elm-form}]
\[
\wh m_{j\to\rt}(\si) \cong\sum_{\si'} \psi
\bigl(\si,\si'\bigr) e^{\lm
_{j\to\rt
}(\si')},\qquad \mu^{\indir}
\mbox{-a.s.}
\]
Thus, if we define $m\dvtx \treespr\to\simplex$, $m_{x\to y}(\si)\cong
e^{\lm
_{x\to y}(\si)}$, then \eqref{elm-form} can be written in terms of
$m$ as
\[
m_{\rt\to j}(\si) \cong\vpsi(\si) \prod_{k\in\pd\rt\setminus j}
\biggl( \sum_{\si_k} \psi(\si,\si_k)
m_{j\to\rt}(\si_k) \biggr),\qquad \mu^\indir\mbox{-a.s.,}
\]
that is, $m\in\hstar$. Then \eqref{ebarh-form} is precisely the
statement that $m$ maps to $\bh$ via \eqref{eembed}, which completes
the proof.
\end{pf}

\begin{pf*}{Proof of Theorem~\ref{tbetheopt}}
By \hypreg\ the set $\hlocf$ of $\bh\in\hloc$ for which $\Phi(\be
,B,\bh
)>-\infty$ is nonempty and does not depend on $(\be,B)$, so without
loss we will restrict to $\bh\in\hlocf$.

Again by \hypreg, the functions $(\be,B)\mapsto\Phi_\mu(\be,B,\bh)$
indexed by $\bh\in\hlocf$ are uniformly equicontinuous on compact
regions of $(\be,B)$: for any $\ep>0$
there exists $\de>0$ sufficiently small so that if $(\be,B)$ and
$(\be
',B')$ are within distance $\de$, then
$\abs{\Phi_\mu(\be,B,\bh)-\Phi_\mu(\be',B',\bh)}<\ep$
for all $\bh\in\hlocf$. Let $\bh\in\hlocf$ such that $\Phi_\mu
(\be,B,\bh
) \ge\wt\Phi_\mu(\be,B)-\ep$. Then
\[
\wt\Phi_\mu\bigl(\be',B'\bigr) \ge
\Phi_\mu\bigl(\be',B',\bh\bigr) \ge\wt
\Phi_\mu(\be,B)-2 \ep
\]
for all $(\be',B')$ within distance $\de$ of $(\be,B)$. Reversing the
roles of $(\be,B)$ and $(\be',B')$ completes the proof of part (a). The statement of part (b) is a
summary of the results of Lemma~\ref{lembetheext},
Propositions~\ref{pinter} and~\ref{popt}.
\end{pf*}

We supplement Proposition~\ref{popt} by computing the second
derivatives\break $\pd_\eta^2\Phi_\mu(\bh+\eta\bde)$ at interior stationary
points $\bh$, giving a criterion to verify that such points are local
maximizers.

\begin{ppn}\label{plocmax}
For permissive $\vec{\psi}$, let $\bh\in\hlocopsi$ be a stationary point
of $\Phi_\mu$, and let $\prr{\de}\in\hlocpmpsi$ with $\abs{\pr
{\de}}\le
\abs{\bh}$. Then $\bh$ is a local maximizer of $\Phi$ on the
one-dimensional space $\hloc\cap\set{\bh+\eta\prr{\de}\dvtx \eta\in
\R}$ if
and only if
%
\begin{eqnarray}
\label{elocmax} 4\pd_\eta^2\Phi_\mu(\bh+\eta
\bde)|_{\eta=0} &=&\E_\mu \biggl[ 2(D_\rt-1)\bigl\langle{(
\vde_\rt/\vh_\rt)^2}\bigr\rangle_{\vh_\rt} -\sum
_{j\in\pd\rt} \bigl\langle{(\prr{\de}_{\rt j} /
\bh_{\rt j})^2}\bigr\rangle_{\bh_{\rt j}} \biggr]
\nonumber
\\[-8pt]
\\[-8pt]
\nonumber
&\le&0,
\end{eqnarray}
or equivalently
%
\begin{eqnarray}
\label{elocmaxunim}\quad && \E_\mu\bigl[\bigl\langle{(\vde_\rt/
\vh_\rt)^2}\bigr\rangle_{\vh_\rt}\bigr]
\nonumber
\\[-8pt]
\\[-8pt]
\nonumber
&&\qquad\ge\f12
\E_\mu \biggl[ \sum_{j\in\pd\rt} \bigl(\bigl\langle {(
\vde_\rt/\vh_\rt)^2}\bigr\rangle_{\vh_\rt} +\bigl\langle{(
\vde_h/\vh_h)^2}\bigr\rangle_{\vh_h} -\bigl\langle{(\prr{
\de}_{\rt j} /\bh_{\rt j})^2}\bigr\rangle_{\bh_{\rt j}} \bigr)
\biggr].
\end{eqnarray}
It is a strict local maximizer if \eqref{elocmax} and \eqref
{elocmaxunim} hold with strict inequality.
\end{ppn}

\begin{pf}
For $\bh\in\hlocopsi$ and $\prr{\de}\in\hlocpmpsi$ with $\abs
{\prr{\de
}}\le\abs{\bh}$, arguing as in the proof of Proposition~\ref{popt} gives
\begin{eqnarray*}
&&2\pd_\eta\Phi_\mu(\bh+\eta\bde)|_{\eta=0}\\
&&\qquad=\lim
_{\eta\to0}R^\eta(\prr{\de}) =R^0(\de)
\\
&&\qquad\equiv\E_\mu \biggl[ 2\angl{\vxi}_{\vde_\rt}
+2(D_\rt-1)\angl{\log\vh_\rt}_{\vde_\rt} +\sum
_{j\in\pd\rt} \bigl( \angl{\xi}_{\prr{\de}_{\rt j}} -\angl{\log
\bh_{\rt j}}_{\prr{\de}_{\rt j}} \bigr) \biggr].
\end{eqnarray*}
If $\bh$ is further a stationary point of $\Phi_\mu$, then, for
$\eta<1$,
\begin{eqnarray*}
T^\eta(\prr{\de}) &\equiv&\f2\eta R^\eta(\prr{\de}) =\f2
\eta\bigl[R^\eta(\prr{\de})-R^0(\prr{\de})\bigr]
\\
&=& 2 \E_\mu \biggl[ 2(D_\rt-1) \bigl\langle{f^\eta(
\vde_\rt/\vh_\rt)}\bigr\rangle_{\vde_\rt} -\sum
_{j\in\pd\rt}\bigl\langle{f^\eta(\prr{\de}_{\rt j}/
\bh_{\rt
j})}\bigr\rangle_{\prr{\de
}_{\rt j}}
\\
&&\hspace*{22pt}{} +2(D_\rt-1) \bigl\langle{g^\eta(\vde_\rt/
\vh_\rt)}\bigr\rangle_{\vh_\rt} -\sum_{j\in\pd\rt}
\bigl\langle{g^\eta(\prr{\de}_{\rt j} /\bh_{\rt j})}\bigr\rangle_{\bh_{\rt j}}
\biggr],
\end{eqnarray*}
where $g^\eta(r)\equiv[f^\eta(r)-r]/\eta$, with $\lim_{\eta\to
0}g^\eta
(r)=-r^2/2$. Since $\abs{\prr{\de}/\bh}\le1$, it follows by dominated
convergence that
\begin{eqnarray*}
4\pd_\eta^2\Phi_\mu(\bh+\eta
\bde)|_{\eta=0}& =&\lim_{\eta\to0} T^\eta(\prr{\de})
=T^0(\prr{\de}) \\
&\equiv&\E_\mu \biggl[ 2(D_\rt-1)
\bigl\langle{(\vde_\rt/\vh_\rt)^2}\bigr\rangle_{\vh_\rt} -
\sum_{j\in\pd\rt} \bigl\langle{(\prr{\de}_{\rt j} /
\bh_{\rt j})^2}\bigr\rangle_{\bh_{\rt j}} \biggr].
\end{eqnarray*}
The stationary point $\bh$ is a local maximizer on $\hloc\cap\set
{\bh
+\eta\prr{\de}\dvtx \eta\in\R}$ if and only if $\pd_\eta^2\Phi_\mu
(\bh+\eta
\bde)|_{\eta=0}\le0$, which gives \eqref{elocmax}. Condition
\eqref
{elocmaxunim} is equivalent by an application of unimodularity.
\end{pf}

\section{Application to Ising and Potts models}\label{secisingpotts}

In this section we apply Theorem~\ref{tfm} to prove our results for
the ferromagnetic Ising and Potts models, Theorems~\ref{tising}--\ref
{tpottsreg}. Although both models have regimes of multiple fixed
points, monotonicity arguments allow us to restrict the space of fixed
points. In the Ising model we can restrict to a unique fixed point and
give a complete verification of the Bethe free energy prediction; in
the Potts model with $q>2$ there remain regimes of nonuniqueness where
we can only provide bounds.

\subsection{Ising model}\label{ssecising}

We first prove Theorem~\ref{tising}. Recall definition \eqref{eising}
for the Ising measure $\nu^{\be,B}_G$ for a finite graph $G=(V,E)$, and
more generally (from Definition~\ref{dbd}) the Ising measures $\nu
^{\free,\be,B}_{U,G}$ and $\nu^{+,\be,B}_{U,G}$ for a finite sub-graph
$U$ of a (possibly infinite) graph $G$ with free and $+$ boundary
conditions. We will make use of the following direct consequence of the
classical Griffiths's inequality; see, for example, \cite{MR2108619}, Theorem~IV.1.21.

\begin{lem}\label{lemgrif}
For the Ising model with parameters $\be,B\ge0$ on $U$ a finite
sub-graph of a graph $G$ with
boundary conditions $\ddagger\in\set{\free,+}$, the \emph{magnetization}
$\anglb{\si_v}^{\ddagger,\be,B}_{U,G}$ at vertex $v\in U$ is nonnegative,
nondecreasing in $\be,B$, nondecreasing in $U$ for $\ddagger=\free$
and nonincreasing in $U$ for $\ddagger=+$.
\end{lem}

Recall from Section~\ref{sssecintroising} the definitions of $\vh
^{t,\ddagger}_T$ for $\ddagger\in\set{\free,+}$; the measure $\vh
^{t,\ddagger}_T$ is parametrized by the corresponding magnetization
$\vm
^{t,\ddagger}_T\equiv\vh^{t,\ddagger}_T(+)-\vh^{t,\ddagger
}_T(-)$. By
Lemma~\ref{lemgrif}, $\vm^{t,\free}_T$ is nondecreasing in $t$ while
$\vm^{t,+}_T$ is nonincreasing, so there exist well-defined limits
$\vm^\ddagger_T(\be,B)\equiv\lim_{t\to\infty} \vm^{t,\ddagger
}_T(\be
,B)$. The following result from~\cite{MR2733399}, an extension of
\cite{MR2650042}, Lemma~4.3, shows that these limits agree on \emph{any}
$T\in\treesrt$.

\begin{lem}[({\cite{MR2733399}, Lemma~3.1})]\label{lemisingunique}
For the Ising model \eqref{eising} on an infinite tree~$T$ with $\be
,B>0$, there exists a constant $C\equiv C(\be,B)$ such that
\[
\vm^{t,+}_T-\vm^{t,\free}_T \le C/t\qquad
\forall t \ge1.
\]
\end{lem}

By this result we can define $h\in\cH$ by
$h_{x\to y}
=\vh^\free_{\subt{x\to y}{}}
=\vh^+_{\subt{x\to y}{}}$,
and we now proceed to verify the Bethe prediction $\phi(\be,B)=\Phi
_\mu
(\be,B,h)$.

\begin{pf*}{Proof of Theorem~\ref{tising}}
The Ising model \eqref{eising} is of form \eqref{efm} with $\spins
=\set{\pm1}$, $\xi(\si,\si')=\be\si\si'$ and $\vxi(\si)=B\si
$, so
\hypreg\ and \hypint\ are
clearly satisfied (with no additional moment conditions on $D_\rt$,
since $\psi>0$).
It follows directly from the recursive structure of the tree that $h\in
\hstar$. It will be shown in Lemma~\ref{lpottsbd}
that for $\be\ge0$ fixed,
\[
\lim_{B\to\infty}\limsup_{n\to\infty}\bigl|
\phi_n(\be,B)-\Phi _\mu(\be,B,h)\bigr|=0,
\]
so to prove the theorem we will interpolate from $(\be,B)$ to $(\be
,B_1)$, then take $B_1 \to\infty$.

It follows from Lemmas~\ref{lemgrif} and~\ref{lemisingunique}
that for $T\in\treesrt$, $\vm^{\free}_T(\be,B)=\vm^+_T(\be
,B)\equiv
m_T(\be,B)$ is the increasing limit of $\vm^{t,\free}_T(\be,B)$ and the
decreasing
limit of $\vm^{t,+}_T(\be,B)$. The $\vm^{t,\ddagger}(\be,B)$ are
continuous and nondecreasing in $\be,B$, so $m$ inherits these
properties by the same argument as in the proof of Theorem~\ref{tis},
and so (since it takes values in $[-1,1]$) is of uniformly bounded
total variation. This verifies both \hypdiffbeta\ and \hypdiffB\
(though we will use only the latter).

We conclude by showing [cf. \eqref{edBassump}] that
\[
\lim_{n\to\infty} \E_n \bigl[ \bigl\langle\pd\vxi(
\si_{I_n})\bigr\rangle^{\be,B}_n \bigr] = a^\vertex(
\be,B) \equiv\E_\mu\bigl[\bigl [\!\bigl[\pd_B\vxi(
\si_\rt)\bigr]\!\bigr]^{h,\be,B}_T \bigr].
\]
Here $\pd_B \vxi(\si) = \si$, and it follows from Lemma~\ref{lemgrif},
our assumption of $G_n \lwc\mu$ and Fatou's lemma that
\begin{eqnarray*}
\E_\mu\bigl[ \danglb{\si_\rt}^{h^\free,\be,B}_T
\bigr] &\le&\liminf_{t\to\infty} \E_\mu\bigl[ \anglb{
\si_\rt}^{\free,\be
,B}_{\subt
{}{t},T} \bigr] \le\liminf
_{n\to\infty} \E_n \bigl[ \anglb{\si_{I_n}}^{\be,B}_n
\bigr]
\\
&\le&\limsup_{n\to\infty} \E_n \bigl[\anglb{
\si_{I_n}}^{\be,B}_n\bigr] \le\limsup
_{t\to\infty} \E_\mu\bigl[ \anglb{\si_\rt}^{+,\be
,B}_{\subt
{}{t},T}
\bigr] \\
&\le&\E_\mu\bigl[ \danglb{\si_\rt}^{h^+,\be,B}_T
\bigr].
\end{eqnarray*}
The left-most and right-most expressions coincide by Lemma~\ref
{lemisingunique} so equality holds throughout.

By Theorem~\ref{tfm}(b), $\phi(\be,B)=\Phi(\be
,B,h^+)=\Phi
(\be,B,h^\free)$ for $\be\ge0$, \mbox{$B>0$}.
Since $\phi_n$ is symmetric in $B$ and continuous at $B=0$ (uniformly
in $n$),
we have $\phi(\be,B)=\phi(\be,-B)$ and $\phi(\be,0) = \lim_{B\to
0} \phi
(\be,B)$.
\end{pf*}

\subsection{Potts model}\label{ssecpotts}

We now apply Theorem~\ref{tfm} to deduce our result (Theorem~\ref
{tpotts}) for the Potts model \eqref{epotts} with $\be,B\ge0$. From
now on we let $\spins\equiv[q]$ with $q\ge2$. It will be convenient to
generalize \eqref{epotts} to the \emph{inhomogeneous} Potts model
\[
\nu^{\vec{\be},\vec B}_G(\vec{\si}) \cong\exp \biggl\{ \sum
_{(ij)\in E} \be_{ij}\cdot\Ind{\si_i=
\si_j} + \sum_{i\in V} B_i\cdot
\Ind{\si_i=1} \biggr\},\qquad \vec{\si}\in\spins^V.
\]

We now introduce the coupling of the Potts model with a random-cluster
model which we use to obtain monotonicity properties. The following
representation is as in \cite{arXiv09011625}; see also \cite{MR1757955}. If
$G=(V,E)$ is a
finite graph, let $G^\ghost$ be the graph formed by adding an edge from
every $v\in V$ to a ``ghost vertex'' $\vghost$, that is, $G^\ghost
=(V^\ghost,E^\ghost)$ where $V^\ghost= V\cup\set{\vghost}$ and
$E^\ghost=E\cup\set{(v,\vghost) \dvtx v\in V}$. Writing $\vec{\si}$ for
elements of $\spins^{V^\ghost}$ and $\vec{\eta}$ for elements of
$\set
{0,1}^{E^\ghost}$ (\emph{bond configurations}), consider the
probability measure on pairs $(\vec{\si},\vec{\eta})$ defined by
%
\begin{eqnarray}
\label{ees}\qquad  &&\es^{\vec{\be},\vec B}_G(\vec{\si},\vec{\eta})
\nonumber
\\[-8pt]
\\[-8pt]
\nonumber
&&\qquad \cong\Ind{
\si_\vghost=1} \prod_{\eta_{ij}=1} \bigl\{
\bigl(e^{\be_{ij}\cdot\Ind{\si_i=\si_j}}-1\bigr) \bigr\} \prod_{\eta_i=1}
\bigl\{ \bigl(e^{B_i\cdot\Ind{\si_i=\si_\vghost}}-1\bigr) \bigr\}.
\end{eqnarray}
The marginal on $\vec{\si}_V$ is the inhomogeneous Potts measure $\nu
^{\vec{\be},\vec B}_G$, while the marginal on $\vec{\eta}$ is the
(\emph
{inhomogeneous}) \emph{random-cluster measure}
%
\begin{equation}
\label{erc} \rc^{\vec{\be},\vec B}_G(\vec{\eta}) \cong\prod
_{e \in E^\ghost} p_e^{\eta_e} (1-p_e)^{1-\eta_e}
\prod_{C\in\eta} \Thet(C),
\end{equation}
where $p_{ij} \equiv1-e^{-\be_{ij}}$ for $(i,j)\in E$ and $p_{i
\vghost
} \equiv1-e^{-B_i}$ for $i\in V$,
and the last product is taken over connected components $C$ of $\vec
{\eta}
$, with $\Thet(C) = q$ unless $\vghost\in C$ in which case $\Thet
(C)=1$. Given a configuration $\vec{\eta}$, a realization of the
conditional law $\es^{\be,B}_G(\vec{\si}=\cdot|\vec{\eta})$ is obtained
by choosing a constant spin on each connected component $C$ of $\vec
{\eta}
$ independently and uniformly over $[q]$, except for $C$ containing
$\vghost$ which is given spin $1$.

For a detailed account the random-cluster model, see \cite{MR2243761}; we will use only the following basic properties:

\begin{ppn}\label{ppnrcmonot}
The random-cluster measure $\rc^{\vec{\be},\vec B}_G$ is FKG. It is also
increasing, in the sense of stochastic domination, in $(\vec{\be
},\vec B)$.
\end{ppn}

\begin{pf}
The FKG property follows by a straightforward modification of the
proof of \cite{MR1757955}, Theorem~III.1(i). Monotonicity in $(\vec
{\be
},\vec
B)$ follows by modifying the proof of \cite{MR2243761}, Theorem~3.21.
\end{pf}

Recalling Definition~\ref{dbd}, for $U$, a finite sub-graph of a graph
$G$ and $\ddagger\in\set{\free}\cup[q]$ (with $\free=$ free), let
$\nu
^{\ddagger,\be,B}_{U,G}$
denote the Potts model on $U$ with $\ddagger$ boundary conditions.
\begin{cor}\label{corpottsmonot}
For the Potts model with parameters $\be,B\ge0$ on $U$ a finite
sub-graph of a graph $G$ with boundary conditions $\ddagger\in\set
{\free
,1}$, and for any vertices $v,w\in U$, the quantities
\[
\nu^{\ddagger,\be,B}_{U,G}(\si_v=1), \qquad\nu^{\ddagger,\be
,B}_{U,G}(
\si _v=\si_w)
\]
are nondecreasing in $\be$ and $B$, nonincreasing in $U$ for
$\ddagger=1$
and nondecreasing in $U$ for $\ddagger=\free$.
\end{cor}

\begin{pf} Note that $\nu^{\free,\be,B}_{U,G}$ is the marginal on $\vec
{\si}_U$
of the measure $\es^{\vec{\be},\vec B}_G$ with
\[
B_i = B\qquad \forall i\in V, \qquad\be_e = \be\cdot\Ind{e \in
E_U}.
\]
Similarly, $\nu^{1,\be,B}_{U,G}$ is the marginal on $\vec{\si}_U$
of the
measure $\es^{\vec{\be}',\vec B'}_G$ with
\[
\be_e'=\be\qquad\forall e\in E,\qquad B_i'
= B\cdot\Ind{i\in V_U} + \infty\cdot\Ind{i \notin
V_U}.
\]
Clearly, $(\vec{\be},\vec B)$ is nondecreasing in $U$ while $(\vec
{\be}
',\vec B')$ is nonincreasing, and both are nondecreasing in $\be,B$.
The result therefore follows from Proposition~\ref{ppnrcmonot} by
showing that for any $(\vec{\be},\vec B)$,
the conditional probabilities $\es^{\vec{\be},\vec B}_G(\si
_v=1|\vec{\eta}
)$ and $\es^{\vec{\be},\vec B}_G(\si_v=\si_w|\vec{\eta})$
are monotone functions of $\vec{\eta}$. Indeed, letting $\es\equiv
\es
^{\vec{\be},\vec B}_G$ and writing $v\conn w$ to indicate that $v,w$
belong to the same connected component of~$\vec{\eta}$, we have
\begin{eqnarray*}
\es(\si_v=1|\vec{\eta}) &=& \mathbf{1}\bigl\{v\conn\vghost\bigr\} + \f{1-\Ind{v\conn
\vghost}} {q},
\\
\es(\si_v=\si_v|\vec{\eta}) &=& \Ind{v\conn w} +
\f{1-\Ind{v\conn w}} {q}.
\end{eqnarray*}
These are increasing functions of $\vec{\eta}$ so the proof is complete.
\end{pf}

Under the measures with $\ddagger\in\set{\free,1}$, any one-vertex
marginal must be uniform on the spins $\ne1$, and so is characterized
by the probability given to spin $1$. In particular, recall from
Section~\ref{sssecintropotts} the definitions of $\vh^{t,\ddagger}_T$
for $\ddagger\in\set{\free,1}$; existence of the $t\to\infty$ limits
$\vh^\ddagger_T$ is now justified by Corollary~\ref{corpottsmonot},
so we can define $h^\ddagger\in\cH$ by $h^\ddagger_{x\to y}=\vh
^\ddagger
_{\subt{x\to y}{}}$. The following lemma gives the boundary values for
the interpolation in $(\be,B)$ using $h^\ddagger$:

\begin{lem}\label{lpottsbd}
For the Potts model on $G_n\lwc\mu$, let
\begin{eqnarray*}
\wt\Phi_\mu(\be,B)&\equiv& B+\be\E_\mu[D_\rt]/2
+ \E_\mu\bigl[\bar \vph\bigl(\abs {T}\bigr)\bigr],\\
 \bar\vph(n) &\equiv&\bar
\vph^B(n) \equiv n^{-1}\log\bigl(1 + (q-1) e^{-Bn}
\bigr).
\end{eqnarray*}\vspace*{-12pt}
\begin{longlist}[(a)]
\item[(a)]
For all $B\in\R$ and any $h \in\cH$, $\phi(0,B)=\log
(e^B+q-1)=\Phi_\mu
(0,B,h)$.

\item[(b)]
For $\be\ge0$ and $h\in\hstar$,
\[
\lim_{B\to\infty} \limsup_{n\to\infty} \bigl|
\phi_n(\be,B)-\wt\Phi_\mu(\be,B)\bigr| =0=\lim
_{B\to\infty}z\bigl|\Phi_\mu(\be,B,h)-\wt
\Phi_\mu(\be,B)\bigr|.
\]

\item[(c)]
For $B\ge0$,
$\lim_{\be\to\infty} \limsup_{n\to\infty}
\abs{\phi_n(\be,B)-\wt\Phi_\mu(\be,B)}
=0$.

\item[(d)]
For $B>0$ and $\ddagger\in\set{\free,1}$,
$\lim_{\be\to\infty}\abs{\Phi_\mu(\be,B,h^\ddagger)-\wt\Phi
_\mu(\be,B)}=0$.
\end{longlist}
\end{lem}

\begin{pf}
(a) At $\be=0$, $\psi\equiv1$ so the spins
are independent. Thus, for all $n\ge1$, $h\in\cH$ and $T\in\treesrt$,
\[
\phi_n(0,B)=\log\bigl(e^B+q-1\bigr)=
\Phi_T^\vertex(0,B,h)=\Phi_T(0,B,h),
\]
since $\Phi^{(\rt j)}_T \equiv0$ for all $j \in\pd\rt$.

(b) The value of $Z_n(\be,B)$ is bounded below
by considering only the ground state $\vec{\si}\equiv1$,
and bounded above by decomposing $\spins^V$ according to the subset of
$k$ vertices where the spin is not $1$.
For $\be\ge0$ this gives
\[
1 \le Z_n(\be,B) e^{- B n - \be|E_n|} \le \sum
_{k=0}^n \pmatrix{n\cr k} (q-1)^k
e^{-Bk} = \bigl(1+(q-1)e^{-B}\bigr)^n,
\]
so if we define $\bar\phi_n (\be,B) \equiv\phi_n(\be,B)- B - \be
\E
_n[\abs{E_n}]/n$, then\break $\lim_{B\to\infty}\limsup_{n\to\infty
}\abs{\bar
\phi_n(\be,B)}=0$.
Recalling \eqref{eedgefreq}, this proves the left identity in (b).

We next define
\begin{eqnarray*}
\bar\Phi^\vertex_T&\equiv&\Phi^\vertex_T-B-
\be D_\rt,\qquad \bar\Phi^\edge_T\equiv
\Phi^\edge_T-\be D_\rt/2, \qquad\bar
\Phi_T \equiv\bar\Phi^\vertex_T-\bar
\Phi^\edge_T,\\
 \bar\Phi_\mu&\equiv&
\E_\mu\bar\Phi_T,
\end{eqnarray*}
so that to prove the right identity in (b) it
suffices to show $\lim_{B\to\infty}\bar\Phi_\mu(\be,\break B,h)=0$
for any $h\in\hstar$. Indeed,
\eqref{ehrecurs} gives that $\mu$-a.s., $\lim_{B\to\infty}h^{\be
,B}_{\rt\to j} (\si)=\Ind{\si=1}$ for all $j\in\pd\rt$,
hence also $\lim_{B\to\infty}h^{\be,B}_{j \to\rt}(\si) =\Ind
{\si=1}$
for all $j\in\pd\rt$ by equivalence of $\mu^\indir$ and $\mu
^\outdir$. Thus
\[
\lim_{B\to\infty}\bar\Phi^\vertex_T(\be,B,h) =
0= \lim_{B\to\infty}\bar\Phi^\edge_T(
\be,B,h),\qquad \mu\mbox{-a.s.}
\]
It is easily verified that
%
\begin{equation}
\label{ebarphibds} -\be D_\rt\le\bar\Phi^\vertex_T(
\be,B,h) \le\log q,\qquad -\be D_\rt/2 \le\bar\Phi^\edge_T(
\be,B,h) \le0,
\end{equation}
so $\bar\Phi_\mu(\be,B,h)\to0$ by dominated convergence.

(c)
Suppose first that $G_n$ is connected. Then $Z_n(\be,B)$ is bounded
below by considering only the $q$ constant-spin configurations, and
bounded above by decomposing $\spins^V$ according to the subset of
$\ell
$ edges across which the spins disagree. Since $G_n$ is connected,
removing $\ell$ edges leaves at most $\ell+1$ connected components, of
sizes $k_0,\ldots,k_\ell$ summing to $n$. Therefore, with $\vph
(n)\equiv
\vph^B(n)\equiv n\bar\vph^B(n)$, we have
\[
e^{\vph(n)}\le Z_n(\be,B) e^{-Bn-\be\abs{E_n}} \le\sum
_{\ell=0}^{\abs{E_n}} \pmatrix{\abs{E_n} \cr\ell}
e^{-\be
\ell} \max_{k_0,\ldots,k_\ell} \Biggl\{ \exp \Biggl\{ \sum
_{r=0}^\ell\vph(k_r) \Biggr\}
\Biggr\},
\]
where the maximum is taken over $k_0,\ldots,k_\ell\in\Z_{\ge0}$ summing
to $n$. By convexity of $\vph$ this maximum is achieved with $k_r=n$
for some $r$, so
%
\begin{eqnarray}
\label{elogZubd} \vph(n)& \le& n\bar\phi_n(\be,B) \le\vph(n)+
\E_n \Biggl[ \log \Biggl\{ \sum_{\ell=0}^{\abs{E_n}}
\pmatrix{\abs{E_n} \cr\ell } e^{-\be\ell} q^\ell \Biggr
\} \Biggr]
\nonumber
\\[-8pt]
\\[-8pt]
\nonumber
&=& \vph(n)+\E_n\bigl[\abs{E_n}\bigr]\log
\bigl(1+qe^{-\be}\bigr).
\end{eqnarray}
If $G_n$ has connected components $C^j=(V^j,E^j)$, $j \ge1$, with
$\abs
{V^j}=n^j$, then clearly $Z_n(\be,B) = \prod_j Z_{C^j}(\be,B)$, so
%
\begin{equation}
\label{ebdcnctdg} 0\le\bar\phi_n(\be,B)-\f{1} {n}\E_n
\biggl[\sum_j \vph\bigl(n^j\bigr)
\biggr] \le\f{1} {n} \E_n\bigl[\abs{E_n}\bigr] \log
\bigl(1+qe^{-\be}\bigr).
\end{equation}
With $j(i)$ denoting the index of the connected component of $G_n$
containing vertex~$i$, we have $n^{-1}\E_n[\sum_j \vph(n^j)]=\E
_n[\bar
\vph(n^{j(I_n)})]$. Then, since $\bar\vph'(n)\le0$,
\[
\E_n\bigl[\bar\vph\bigl(\bigl|\ball{t} {I_n}\bigr|\bigr) \cdot\mathbf{1}\bigl\{
\ball{t} {I_n}=C^{j(I_n)}\bigr\}\bigr] \le\E_n\bigl[
\bar\vph\bigl(n^{j(I_n)}\bigr)\bigr] \le\E_n\bigl[\bar\vph\bigl(
\bigl|\ball{t} {I_n}\bigr|\bigr)\bigr].
\]
Since $G_n \lwc\mu$, letting $n\to\infty$ followed by $t\to\infty
$ in
the above inequalities gives $\E_n[\bar\vph(n^{j(I_n)})] \to\E_\mu
[\bar
\vph(\abs{T})]$, and so (c) follows from
\eqref{ebdcnctdg} by taking first $n\to\infty$ and then $\be\to
\infty$.

(d) Clearly $h_T^\free=h_T^1$ for any
finite $T\in\treesrt$ (as $\pd\subt{}{t}=\varnothing$ for large enough
$t$). In the $\be\to\infty$ limit only the constant-spin configurations
contribute, so
%
\begin{equation}
\label{ehTr} \lim_{\be\to\infty} h^{\ddagger,\be,B}_T(\si)
= e^{-\vph(\abs{T})-B\abs{T}(1-\Ind{\si=1})},\qquad \ddagger\in \set{\free,1}.
\end{equation}
For $T$ infinite, recall from Corollary~\ref{corpottsmonot} that
$h^\free_{\subt{}{t}}(1)
\le h^\free_T(1) \le h^1_T(1)$, so if $B>0$, then
\[
1=\lim_{t\to\infty} \lim_{\be\to\infty}
h^{t,\free,\be,B}(1) \le\lim_{\be\to\infty} h^{\free,\be,B}(1) \le\lim
_{\be\to\infty} h^{1,\be,B}(1),
\]
so that \eqref{ehTr} again holds for $T$ infinite.
We then compute
\begin{eqnarray*}
\lim_{\be\to\infty}\bar\Phi^\vertex_T\bigl(
\be,B,h^\ddagger\bigr) &=& -\sum_{j\in\pd\rt} \vph\bigl(
\abs{\subt{j\to\rt} {}}\bigr) + \vph\bigl(\abs{T}\bigr),
\\
\lim_{\be\to\infty}\bar\Phi^\edge_T\bigl(
\be,B,h^\ddagger\bigr) &=&-\f{1} {2}\sum_{j\in\pd\rt}
\vph\bigl(\abs{\subt{j\to\rt} {}}\bigr) -\f{1} {2}\sum_{j\in\pd\rt}
\vph\bigl(\abs{\subt{\rt\to j} {}}\bigr) +\f{D_\rt} {2} \vph\bigl(\abs{T}\bigr),
\end{eqnarray*}
$\mu$-a.s.,
where the first identity uses $\abs{T}=1+\sum_{j\in\pd\rt}\abs
{\subt
{j\to\rt}{}}$ and the second uses $\abs{T}=\abs{\subt{\rt\to
j}{}}+\abs
{\subt{j\to\rt}{}}$. Convergence also holds in $\mu$-expectation, using
the upper bounds in \eqref{ebarphibds} together with
\begin{eqnarray*}
\bar\Phi^\vertex_T\bigl(\be,B,h^\ddagger\bigr)& \ge&
\sum_{j\in\pd\rt} \log h^{\free,\be,B}_{j\to\rt}(1),\\
\bar\Phi^\edge_T\bigl(\be,B,h^\ddagger\bigr) &\ge&
\f{1} {2} \sum_{j\in\pd\rt} \log h^{\free,\be,B}_{j\to\rt}(1)
+\f{1} {2} \sum_{j\in\pd\rt} \log h^{\free,\be,B}_{\rt\to j}(1),
\end{eqnarray*}
and the fact that $h^{\free,\be,B}_{x\to y}(1)\ge1/q$ for $\be,B\ge0$
(by Corollary~\ref{corpottsmonot}). Thus, using unimodularity of
$\mu
$, we have
\[
\lim_{\be\to\infty}\bar\Phi_\mu\bigl(\be,B,h^\ddagger
\bigr) = \E_\mu\bigl[ (1-D_\rt/2) \vph\bigl(\abs{T}\bigr) \bigr],
\]
and we conclude by showing that this coincides with $\E_\mu[\bar\vph
(\abs{T})]$. The case $\abs{T}=\infty$ is trivial; otherwise, another
application of unimodularity gives
\begin{eqnarray*}
\f{1} {2} \E_\mu\bigl[ D_\rt\vph\bigl(\abs{T}\bigr) \bigr] &=&
\f{1} {2} \E_\mu \biggl[ D_\rt\sum
_{x\in T} \bar\vph\bigl(\abs{T}\bigr) \biggr] = \f{1} {2}
\E_\mu \biggl[ \sum_{x\in T}
D_x \bar\vph\bigl(\abs{T}\bigr) \biggr]
\\
&=& \E_\mu\bigl[\bar\vph\bigl(\abs{T}\bigr) \abs{E_T}\bigr] =
\E_\mu\bigl[\vph\bigl(\abs{T}\bigr)\bigr]-\E_\mu\bigl[\bar\vph\bigl(
\abs{T}\bigr)\bigr].
\end{eqnarray*}
Therefore, $\lim_{\be\to\infty}\bar\Phi_\mu(\be,B,h)=\E_\mu
[\bar\vph
(\abs{T})]$ which concludes the proof.
\end{pf}

\begin{pf*}{Proof of Theorem~\ref{tpotts}}
The Potts model \eqref{epotts} is of form \eqref{efm} with $\spins
=[q]$, $\xi(\si,\si') = \be\cdot\Ind{\si=\si'}$, and $\vxi(\si
) = B\cdot
\Ind{\si=1}$, so \hypreg\ and \hypint\ are clearly satisfied. It
follows from the recursive structure of the tree that $h^\ddagger\in
\hstar$ for $\ddagger\in\{\free,1\}$. For part (a),
along any interpolation path contained in $\cR_\mu$, both
\hypdiffbeta\ and \hypdiffB\ are satisfied by Corollary~\ref
{corpottsmonot} and the same argument used in the proof of
Theorem~\ref{tis}. For part (b), \hypdiffbeta\ and
\hypdiffB\ are satisfied by the additional hypothesis of continuity.

The inequalities in part (b) then follow from Theorem
\ref{tfm} once we verify [cf.~\eqref{edbetaassump}, \eqref{edBassump}]
\begin{eqnarray*}
a^\vertex\bigl(\be,B,h^\free\bigr) &\le&\liminf
_{n\to\infty} a^\vertex_n(\be,B) \le\limsup
_{n\to\infty} a^\vertex_n(\be,B) \le
a^\vertex\bigl(\be,B,h^1\bigr),
\\
a^\edge\bigl(\be,B,h^\free\bigr) &\le&\liminf
_{n\to\infty} a^\edge_n(\be,B) \le\limsup
_{n\to\infty} a^\edge_n(\be,B) \le
a^\edge\bigl(\be,B,h^1\bigr),
\end{eqnarray*}
where $a^\vertex_n(\be,B) = \E_n [\anglb{\Ind{\si_{I_n}=1}}^{\be,B}_n]$
and $a^\edge_n(\be,B) = \f{1}{2} \E_n  [ \sum_{j\in\pd I_n}
\anglb
{\Ind{\si_{I_n}=\si_j}}^{\be,B}_{n}  ]$. Indeed, by Vitali's
convergence theorem, the assumption $G_n\lwc\mu$ and Corollary~\ref
{corpottsmonot} [with $U=B_t (I_n) \subseteq G_n$], we have
\[
a^\edge\bigl(\be,B,h^\free\bigr) = \liminf
_{t \to\infty} \f{1} {2} \E_\mu \biggl[ \sum
_{j\in\pd\rt} \bigl\langle{\Ind{\si_\rt=\si_j}}\bigr\rangle^{\free,\be,B}_{\subt{}{t},T}
\biggr] \le\liminf_{n\to\infty} a^\edge_n(
\be,B),
\]
and the other inequalities are proved similarly. Together these
inequalities imply that
\[
\lim_{n \to\infty} \bigl[\phi_n\bigl(
\be',B'\bigr)-\phi_n(\be,B)\bigr] = \Phi
\bigl(\be ',B',h^\ddagger\bigr)-\Phi\bigl(
\be,B,h^\ddagger\bigr)
\]
for any $(\be,B)$ and $(\be',B')$ joined by an interpolation path
contained in $\cR_\mu$. The result of part (a)
then follows by letting $(\be',B')$ approach $\cR_\infty$ and
applying Lemma~\ref{lpottsbd}.
\end{pf*}

\subsection{Potts model with $d$-regular limiting tree}\label{ssecpottsreg}

In this section we prove Theorem~\ref{tpottsreg}, which amounts to
determining the shape of $\cR_{\ne}$ and establishing continuity of
$h^\free$ and $h^1$ in certain regimes.

Since the limiting measure is supported on $\treereg{d}$, only
$h\equiv
h_{(\treereg{d},x\to y)}$ is of relevance. Further, $h^\ddagger$ is
symmetric among the spins $\ne1$ for $\ddagger\in\set{\free,1}$, so
determination of $h^\ddagger$ reduces to solving a univariate recursion
for $h^\ddagger(1)$,
\[
h\mapsto\f{e^B [ e^\be h + (1-h)
]^{d-1}} {e^B [e^\be h + (1-h)
]^{d-1} + (q-1) [ h + ({(1-h)}/ {(q-1)})(e^\be+ q-2
) ]^{d-1}}.
\]
Our result follows from analysis of the fixed points of this mapping;
similar computations have appeared, for example, in \cite{MR0378152,MR714953} so some overlap among the analyses may occur.

A convenient parametrization is given by the log likelihood ratio $r
\equiv\log h - \log[(1-h)/(q-1)]$, in terms of which the recursion becomes
\[
r\mapsto f(r) \equiv f(r;\be,B) = B + (d-1)\log \biggl( \f{e^{\be+r} +
q-1} {e^r+e^\be+q-2} \biggr).
\]
With $f^{(t)}$ the $t$-fold iteration of $f$, let $r^\free$ denote the
increasing limit of $f^{(t)}(0)$ and $r^1$ the decreasing limit of
$f^{(t)}(\infty)$, as $t\to\infty$. The region $\cR_{\ne}$
corresponds to those $\be,B\ge0$ for which $r^\free\ne r^1$.

\begin{lem}\label{lempottsregregion}
There exists $\be_- >0$ such that for $\be\le\be_-$ the map $f$ has
exactly one fixed point for any $B \in\R$. For $\be>\be_-$ there exist
real-valued $B_-(\be)<B_+(\be)$ (smooth in $\be$) such that $f$ has
one, two or three fixed points depending on whether $B$ is in
$[B_-,B_+]^c$, $\set{B_-,B_+}$ or $(B_-,B_+)$. The curves extend
continuously to $B_-(\be_-)=B_+(\be_-)$.
\end{lem}

\begin{pf}
We have
%
\begin{equation}
\label{edfdr} f'(r)=\frac{(d-1) e^r  (e^{\beta}-1 )  (q+e^{\beta
}-1 )} {
(q+e^r+e^{\beta}-2 ) (q+e^{r+\beta}-1 )}
\end{equation}
so $f$ is increasing in $r$ with $f'(r)\to0$ as $r\to\pm\infty$. Since
$f(r;\be,B)=f(\be;r,B)$, it easily follows from \eqref{edfdr} that
$\pd
_\be f(r)$ has the same sign as $r$ while $\pd_\be[f'(r)]>0$. Further
\begin{eqnarray*}
f''(r) &=&-\frac{(d-1) e^{r+\be}  (e^{\beta}-1 )  (q+e^{\beta
}-1 )
(e^{2r}-\al )}{ (q+e^r+e^{\beta}-2 )^2  (q+e^{r+\beta}-1 )^2},\\
 \al&\equiv&(q-1) \bigl(1+
(q-2) e^{-\be}\bigr),
\end{eqnarray*}
with $\al>0$ since $q>1$. Notice that $f''(r)>0$ for $r$ sufficiently
negative and $f''(r)<0$ for $r$ sufficiently positive, with a single
sign change occurring at $(\log\al)/2$ which is zero for $q=2$ and
strictly positive for $q > 2$. This proves that $f$ has between one and
three fixed points. When $B=0$, one fixed point is always given by
$r^\free(\be,0)=0$. Further $f(r;0,0)\equiv0$, so (by monotonicity of
$f'$ in $\be$) there exists $\infty\ge\be_-\ge0$ such that $f'\le1$
everywhere for $\be\le\be_-$, and $f'$ exceeds~$1$ somewhere for
$\be
>\be_-$.

Solving the equation $f'(r)=1$ in terms of $t\equiv e^r$ yields solutions
\[
t_\pm(\be)=-\gam\pm\sqrt{\gam^2-\al},\qquad \gam\equiv
e^\be+ q- 2 - \f{d} {2}\bigl(1-e^{-\be}\bigr)
\bigl(e^\be+q-1\bigr).
\]
Since $\al>0$, $t_\pm(\be)$ are not positive if $\gam>-\sqrt{\al}$,
equal to $\sqrt{\al}>0$ if $\gam=-\sqrt{\al}$, and positive but not
equal if $\gam<-\sqrt{\al}$. If $d \ge2$, it is easy to check that
both $\al$ and $\gam$ decrease smoothly in $\be$, starting at $\gam
|_{\be=0}=q-1$ and $\al|_{\be=0}=(q-1)^2$, so there is a unique value
$\be=\be_->0$ at which $\gam=-\sqrt{\al}$: if $d=2$, then $\be
_-=\infty
$, and if $d>2$, then $\be_-$ is the logarithm of the unique finite
positive root $b_-$ of
%
\begin{equation}
\label{ebetam} (d-2)^2 b^2 + (d-2)^2(q-2)b -
d^2(q-1)=0.
\end{equation}
Hence, the equation $f'(r)=1$ has no solutions for $\be<\be_-$, and it
has solutions $\rho_\pm(\be)\equiv\log t_\pm(\be)$ for $\be\ge
\be_-$,
with $\rho_-(\be_-)=\rho_+(\be_-)$ and $\rho_-(\be)<\rho_+(\be
)$ for
$\be>\be_-$. The values of $B_-(\be)$, $B_+(\be)$ are then given
explicitly by
%
\begin{equation}
\label{epottsregbdB} B_\pm(\be)=\rho_\mp(\be)-f\bigl(
\rho_\mp(\be);\be,0\bigr),
\end{equation}
which clearly meet at $\be=\be_-$ and are smooth for $\be>\be_-$.
\end{pf}

Considering hereafter only $d>2$ (so that $\be_-<\infty$), suppose
$\be
>\be_-$, so that the functions $\rho_\pm$ are defined. Since $\pd
_\be
[f'(r)]>0$, $\rho_-$ and $\rho_+$ must be, respectively, decreasing and
increasing in $\be$. Further, since $f$ has a unique inflection point
at $(\log\al)/2$, we must have $\rho_-(\be)\le(\log\al)/2\le
\rho_+(\be
)$, with strict inequalities unless $\rho_-(\be)=\rho_+(\be)$. For
$q=2$ (Ising), this implies $\rho_-\le0\le\rho_+$ from which it is easy
to see that whenever $B>0$ we have $r^\free(\be,B)=r^1(\be,B)$, which
is then continuous in $(\be,B)$ by the same argument as in the proof of
Theorem~\ref{tis}. When $B=0$, $r^\free(\be,0)$ is zero for all
$\be$,
while $r^1(\be,B)$ is zero for $\be\le\be_-$ and strictly positive for
$\be>\be_-$.

For $q>2$ (Potts), this implies that $\rho_+(\be,B)>0$ while $\rho
_-(\be
,B)\ge0$ if and only if $f'(0;\be,B)\le1$. From the calculations above,
$f'(0)$ is zero at $\be=0$ and increases in $\be$. We therefore define
%
%
\begin{eqnarray}
\label{ebetap} \be_\free&\equiv&\inf\bigl\{\be\ge0\dvtx f(r;\be,0)=r
\mbox{ for some } r>0\bigr\},
\nonumber
\\
\be_+ &\equiv &\inf\bigl\{\be\ge0\dvtx \rho_-(\be)\le0\bigr\} = \inf\bigl\{\be\ge0\dvtx
f'(0;\be,0)\ge1\bigr\} \\
&=& \log\biggl (1 + \f {q} {d-2} \biggr)\nonumber
\end{eqnarray}
[where the formula for $\be_+$ comes from \eqref{edfdr}]. Clearly
$\be
_-\le\be_\free\le\be_+$, and in fact these inequalities are
strict: at
$\be_\free$, $f'$ must exceed one between zero and the positive fixed
point, so $\be_-<\be_\free$.\footnote{Note that $r^1(\be_\free,0)>0$,
that is, the 1-biased fixed point ``arises discontinuously.''}
Likewise, if $f'(0)\ge1$ at $\be=\be_\free$, the concavity of
$f(r)$ at
$r=0$ would imply the existence of a positive fixed point at\vadjust{\goodbreak} some $\be$
below $\be_\free$ which is a contradiction, so $\be_\free< \be_+$. We
refer again to Figure~\ref{fbp} which shows the maps $f(r;\be,B)$ for
the Ising and Potts models at several values of $\be$ while holding
$B=0$. Figure~\ref{fpotts}(b) shows the regime of $(\be,B)$ values
delineated by the curves $B_\pm(\be)$.

\begin{pf*}{Proof of Theorem~\ref{tpottsreg}}
(a) We found above that $\cR_{\ne
}=\varnothing$
for $d=2$ and $\cR_{\ne}=(\be_-,\infty)$ for $q=2$, so suppose
$d,q>2$. If $B>0$, $r^\free=r^1$ holds for all $\be\ge0$ with $\be
\notin
(\be_-,\be_+)$. For $\be\in(\be_-,\be_+)$ there is a closed interval
$[B_-(\be)\vee0,B_+(\be)]$ of $B$ values for which $r^\free<r^1$: this
interval is strictly positive for $\be<\be_\free$ and includes zero for
$\be\ge\be_\free$. If $B=0$, $r^\free=r^1$ for $0\le\be<\be
_\free$ and
$r^\free<r^1$ for $\be\ge\be_\free$. Recalling \eqref{epottsregbdB},
\begin{eqnarray*}
\pd_\be B_\pm(\be) &=& \pd_\be\bigl[
\rho_\mp(\be) - f\bigl(\rho_\mp(\be)\bigr)\bigr] = \bigl[1
- f'\bigl(\rho_\mp(\be)\bigr)\bigr]\pd_\be
\rho_\mp(\be) - (\pd_\be f) \bigl(\rho_\mp (
\be)\bigr)\\
& =& - (\pd_\be f) \bigl(\rho_\mp(\be)\bigr).
\end{eqnarray*}
This has the same sign as $-\rho_\mp(\be)$, which are both negative for
$0\le\be<\be_+$, so the curves $B_\pm(\be)$ are decreasing. Inverting
them gives the curves $\be_\free(B),\be_+(B)$ which delineate the
region $\cR_{\ne}$ as described in the theorem statement, with
$\be
_\free(0)=\be_\free$ and $\be_+(0)=\be_+$.

(b) Away from the boundary of $\cR_{\ne}$,
$h^\free$ and $h^1$ correspond to isolated zeros of a smooth function,
and so are continuous by the implicit function theorem. From part (a), any point of $\cR$ is connected to $\cR
_\infty$ by an interpolation path contained in $\cR$, so applying
Theorem~\ref{tpotts}(a) verifies the Bethe
prediction for $(\be,B)\notin\cR_{\ne}$.

Since changing $B$ only translates $f(r;\be,B)$, it is not difficult to
see that when $\be\in(\be_-,\be_+)$, the function $h^\free(\be
,B)$ is
continuous in $B$ for $B\in[0,B_+(\be)]$ while $h^1(\be,B)$ is
continuous for $B \in[B_-(\be)\vee0,\infty)$. It follows by
Lemma~\ref
{lemlogZderiv} that for $(\be,B)\in\pd\cR_{\ne}$ with $\be
=\be
_\free(B)$, $\phi(\be,B)=\Phi(\be,B,h^\free)$, while for $(\be
,B)\in\pd
\cR_{\ne}$ with $\be\ge\be_+(B)$, $\phi(\be,B)=\Phi(\be,B,h^1)$.

Recall our convention that $\be_0\le\be_1,B_0\le B_1$: by
Theorem~\ref
{tpotts}(b) we may interpolate in $B$ from $(\be
,B_0)\in\cR_{\ne}^\circ$ to $(\be,B_1)\in\cR$ using the message
$h^1$, yielding $\liminf_{n\to\infty}\phi_n(\be,B)\ge\Phi(\be,B,h^1)$
for $(\be,B)\in\cR_{\ne}^\circ$. Likewise, we may interpolate in
$B$ from $(\be,B_0)\in\cR$ to $(\be,B_1)\in\cR_{\ne}$ using
$h^\free$ (and once inside $\cR_{\ne}$ we may also interpolate in
$\be$ using $h^\free$), which gives $\liminf_{n\to\infty}\phi
_n(\be
,B)\ge\Phi(\be,B,h^\free)$ for $(\be,B)\in\cR_{\ne}^\circ$.

Next, since $h^\free(\be,B)$ and $h^1(\be,B)$ are lower and upper
semi-continuous, respectively, in $\be$, and both are nondecreasing in
$\be$, for $0<B<B_+$ we have that $h^\free(\be,B)\incto h^\free(\be
_+(B),B)$ as $\be\incto\be_+(B)$ and $h^1(\be,B)\decto h^1(\be
_\free
(B),B)$ as $\be\decto\be_\free(B)$. Again by Theorem~\ref
{tpotts}(b), we may interpolate in $\be$ from $(\be_0,B)=(\be
_\free
(B),B)\in\pd\cR_{\ne}$ to $(\be_1,B)\in\cR_{\ne}^\circ
$ using
$h^1$, and from $(\be_0,B)\in\cR_{\ne}^\circ$ to $(\be
_1,B)=(\be
_+(B),B)\in\pd\cR_{\ne}$ using $h^\free$, giving
\[
\limsup_{n\to\infty}\phi_n(\be,B) \le \min\bigl\{\wt
\Phi^\free(\be,B),\wt\Phi^1(\be,B)\bigr\}, \qquad(\be,B)\in
\cR_{\ne}^\circ,
\]
which completes the proof.
\end{pf*}\eject

\section*{Acknowledgments}
We thank Allan Sly and Ofer Zeitouni for many helpful conversations.
A. Dembo and N. Sun thank the Microsoft Research Theory Group for
supporting a visit during which part of this work was completed.

%
%



\printaddresses

\end{document}